\documentstyle[12pt]{article}
\begin{document}
\author{S.V. Ludkovsky.}
\title{Semidirect products of groups of loops and
groups of diffeomorphisms of real, complex and quaternion manifolds,
their representations.}
\date{29 June 2004}
\maketitle
\begin{abstract}
This article is devoted to the investigation of
semidirect products of groups of loops and groups
of homeomorphisms or groups of diffeomorphisms of
finite and infinite dimensional real, complex and quaternion manifolds.
Necessary statements about quaternion manifolds with quaternion holomorphic
transition mappings between charts of atlases are proved.
It is shown, that these groups exist and have the infinite dimensional
Lie groups structure, that is, they are continuous or differentiable
manifolds and the composition $(f,g)\mapsto f^{-1}g$ is continuous or
differentiable depending on a class of smoothness of groups.
Moreover, it is demonstrated that in the cases of complex and
quaternion manifolds these groups
have structures of complex and quaternion manifolds respectively.
Nevertheless, it is proved that these groups does not necessarily satisfy
the Campbell-Hausdorff formula even locally besides the exceptional case
of a group of holomorphic diffeomorphisms of a compact complex manifold.
Unitary representations of these groups $G'$ including irreducible
are constructed with the help of quasi-invariant measures
on groups $G$ relative to dense subgroups $G'$.
It is proved, that this procedure provides a family of the cardinality
$card ({\bf R})$ of pairwise nonequivalent irreducible unitary
representations. A differentiabilty of such representations is studied.
\end{abstract}

\section{Introduction.}
Gaussian quasi-invariant measures on loop groups
and diffeomorphism groups of Riemannian manifolds were investigated
earlier (by the author) in \cite{ludan,lurim1,lurim2}.
Traditionally geometric loop groups are considered
as families of mappings $f$ from the unit circle
$S^1$ to a manifold $N$ which map a marked point $s_0\in S^1$
to a marked point $y_0$ of $N$. In fact these mappings are defined
as equivalence classes in the group of $N$-valued
$C^{\infty }$-diffeomorphisms on $S^1$, where two mappings
$f_1$ and $f_2$ are called equivalent if there exists a
$C^{\infty }$-diffeomorphism $\varphi : S^1\rightarrow S^1$
preserving the orientation such that $f_1=f_2\circ \varphi $
and equivalence classes are closures of such families of
equivalent mappings. This concept was generalized and investigated in
\cite{ludan} to include loop groups
as families of $C^{\infty }$-mappings from one Riemannian
manifold $M$ to another $N$, which "preserve marked points"
$s_0\in M$ and $y_0\in N$. Here again equivalence classes
under appropriate $C^{\infty }$-diffeomorphisms are used.
This is done via the construction of an Abelian group from a
commutative monoid with a
unit together with the cancellation property which is available
under rather mild conditions on the finite or infinite dimensional
manifolds $M$ and $N$. In this paper besides Riemannian
manifolds also groups of loops of complex and quaternion manifolds are
investigated.
\par Groups of loops and groups of diffeomorphisms of quaternion
manifolds are defined and investigated here for the first time.
It is necessary to note, that semidirect products of these groups
also were not earlier considered.
\par Holomorphic functions of quaternion
variables were investigated in \cite{luoyst,luoyst2}.
There specific definition of superdifferentiability was considered,
because the quaternion skew field has the graded algebra structure.
This definition of superdifferentiability does not impose the condition
of right or left superlinearity of a superdifferential, since it leads
to narrow class of functions. There are some articles on quaternion
manifolds, but practically they undermine a complex manifold
with additional quaternion structure of its tangent space
(see, for example, \cite{museya,yano} and references therein).
Therefore, quaternion manifolds as they are defined below were not
considered earlier by others authors (see also \cite{lufsqv}).
Applications of quaternions in mathematics and physics can be found in
\cite{emch,guetze,hamilt,lawmich}.
\par If we consider the composition of two nontrivial
$C^n$-loops pinned in the marked point
$s_0$, where $n\ge 1$, then the resulting loop is continuous, 
but in general not of class $C^n$. This can easily be seen from
examples of loops on the unit circle
$S^1$ and the unit sphere $S^2$. If e.g. $f: S^1\rightarrow S^2$
is a $C^n$-loop, with $n\ge 1$. Then $f$ and $f'$
are continuous functions, which may be considered to be defined
on the interval. Then $f(0)=f(1)$ and
$\lim_{\theta \downarrow 0}f'(\theta )=
\lim_{\theta \uparrow 1}f'(\theta )=:f'(0)$. There exists another
$C^n$-loop $g$ such that $g(0)=f(0)$, but with $g'(0)\ne f'(0)$.
As an example of a non-$C^1$ continuous loop $h$ we concatenate
$f$ with $g$ to obtain the loop
$h(\theta ):=f(2\theta )$ for each $\theta \in [0,1/2)$ and
$h(\theta ):=g(2\theta -1)$ for each $\theta \in [1/2,1]$. 
Such mappings $h$ are generally speaking only piecewise
$C^n$-continuous and continuous.
Another reason that, although $S^n\vee S^n$ is a retract
of $S^n$, the manifolds $S^n\vee S^n$ and $S^n$ are not diffeomorphic.
In addition, there exists a continuous mapping from $S^1\times S^n$
onto $S^{n+1}$, but this mapping is not a diffeomorphism (see \cite{swit}).
Naturally, smooth compositions of mappings between loop monoids
and loop groups manifolds with corners (with the corresponding atlases)
are being used. Such compositions and mappings permit us to define
topological loop monoids and
topological loop groups. Below the reader will find two other reasons
to consider manifolds with corners.
\par A commutative monoid is not a free (pinned) loop space, because
it is obtained from the latter by factorization.
In order to construct loop groups, the underlying manifold
has to satisfy some mild conditions.
Finite dimensional manifolds are supposed to be compact.
This condition is not very restrictive, because each locally compact
space can be embedded into its Alexandroff (one-point) compactification
(see Theorem 3.5.11 in \cite{eng}).
If the manifold $M$ is infinite dimensional over $\bf C$, then it is
assumed, that $M$ is embedded as a closed bounded subset in a 
corresponding Banach space $X_M$. In order to define a group structure on a 
quotient space of a free loop space, such an embedding is required. 
\par Let $M$ and $N$ be complex or real manifolds with marked points
$s_0\in M$ and $y_0\in N$ respectively.
By definition the free loop space of these pointed manifolds
consists of those continuous mappings $f: M\to N$ which are (piecewise) 
holomorphic in the complex case or $C^{\infty }$ (in the real case)
on $M\setminus M'$, where $M'$ is a submanifold of $M$ (which may depend
on $f$) of codimension $1$ in $M$, and which have the property
that $f(s_0)=y_0$. There are at least two reasons to consider
such a class of mappings. The first one being the fact that
compositions of elements of a loop group can be defined correctly.
The second reason is the existence of an isoperimetric inequality
(for holomorphic loops) which causes a loop, which is close
to a constant loop $w_0: M \rightarrow \{ y_0 \} $, to be
constant in a neighbourhood of $s_0$ (see Remark 3.2 in \cite{hum}).
\par In this article loop groups of different classes of smoothness are
considered. Classes analogous to Gevrey classes
and also with the usage of Sobolev classes of $f: M\setminus M' \to N$  
are considered for the construction of dense loop 
subgroups and quasi-invariant measures. Henceforth, we consider
not only orientable manifolds $M$ and $N$, but also nonorientable
manifolds.
Loop commutative monoids with the cancellation property
are quotients of families of mappings $f$ from $M$ into a
manifold $N$ with $f(s_0)=y_0$ by the corresponding equivalence
relation. 
For the definition of the equivalence relation 
groups of holomorphic diffeomorphisms are not used here, because of strong
restrictions on their structure caused by 
holomorphicity (see Theorems 1 and 2 in \cite{bomon}).
Groups are constructed from monoids via the algebraic
procedure, which was possibly first described by A. Grothendieck
in an abstract context used in algebraic topology not meaning
concrete groups related with those investigated in this paper.
\par Loop groups are Abelian, non-locally compact for $dim_{\bf R}N>1$
and for them the Campbell-Hausdorff formula is not valid (in 
an open local subgroup). Apart from them, finite dimensional Lie groups 
satisfy locally the Campbell-Hausdorff formula.
This is guarantied, if impose on a locally compact topological
Hausdorff group $G$ two conditions: it is a $C^{\infty }$-manifold
and the mapping $(f,g)\mapsto f\circ g^{-1}$
from $G\times G$ into $G$ is of class $C^{\infty }$.
But for infinite dimensional $G$ the Campbell-Hausdorff 
formula does not follow from these conditions.
Frequently topological Hausdorff groups satisfying these 
two conditions are also called Lie groups, though they 
can not have all properties of finite dimensional 
Lie groups, so that the Lie algebras for them do not play 
the same role as in the finite dimensional case and therefore
Lie algebras are not so helpful.
If $G$ is a Lie group and its tangent space $T_eG$ is a Banach space,
then it is called a Banach-Lie group, sometimes it is undermined,
that they satisfy the Campbell-Hausdorff formula locally for a Banach-Lie 
algebra $T_eG$. In some papers the Lie group terminology
undermines, that it is finite dimensional.
It is worthwhile to call Lie groups satisfying the Campbell-Hausdorff
formula locally (in an open local subgroup) by Lie groups in the 
narrow sense;
in the contrary case to call them by Lie groups in the broad sense.
\par The problem to investigate unitary representations of nonlocally
compact groups with the help of quasi-invariant measures
was formulated in sixties of the 20-th century by I.M. Gelfand
\cite{shav}. For construction of quasi-invariant measures
on groups it can be used a local diffeomorphism of a Lie
group as a manifold with its tangent space, where diffeomorphism
need not be preserving group structure and may be other than
the exponential mapping of a manifold. On tangent space
it is possible to take a quasi-invariant measure, for example,
Gaussian or a transition measure of a stochastic process, for example,
Brownian (see \cite{shav,lurim2} and references therein).
General theorems about quasi-invariance 
and differentiabilty of transition probabilities on the Lie group $G$
relative to a dense subgroup $G'$ were given in \cite{beldal,dalschn}, 
but they permit a finding of $G'$ only abstractly and when a local
subgroup of $G$
satisfies the Campbell-Hausdorff formula. For Lie groups which 
do not satisfy the Campbell-Hausdorff formula even locally
this question has been remaining open, as it was pointed out by Belopolskaya 
and Dalecky in Chapter 6. They have proposed in such cases to investigate
concrete Lie groups that to find pairs $G$ and $G'$.
On the other hand, the groups considered in the present
article do not satisfy the Campbell-Hausdorff formula.
\par Below loop groups and diffeomorphism groups and their semidirect
products are considered
not only for finite dimensional, but also for infinite dimensional 
manifolds.
\par In particular, loop and diffeomorphism groups are important for the 
development of the representation theory of non-locally compact groups. 
Their representation theory has many 
differences with the traditional representation theory of 
locally compact groups and finite dimensional Lie groups, 
because non-locally 
compact groups have not $C^*$-algebras associated with the Haar 
measures and they have not underlying Lie algebras and relations between 
representations of groups and underlying algebras (see also \cite{lubp}). 
\par In view of the A. Weil 
theorem if a topological Hausdorff group $G$ has a quasi-invariant measure 
relative to the entire $G$, then $G$ is locally compact.
Since loop groups $(L^MN)_{\xi }$ are not locally compact, they can not have 
quasi-invariant measures relative to the entire group, but only relative to 
proper subgroups $G'$ which can be chosen dense in  $(L^MN)_{\xi }$,
where an index $\xi $ indicates on a class of smoothness.
The same is true for diffeomorphism groups (besides holomorphic 
diffeomorphism groups of compact complex manifolds).
Diffeomorphism groups of compact complex manifolds
are finite dimensional Lie groups (see \cite{kobtg}
and references therein).
It is necessary to note that there are quite another groups
with the same name loop groups, but they are infinite dimensional
Banach-Lie groups of $C^{\infty }$-mappings $f: M\to H$
into a finite dimensional Lie group $H$ with the pointwise group
multiplication of mappings with values in $H$
such that $C^{\infty }(M,H)$ satisfies locally the Campbell-Hausdorff
formula.
\par The traditional geometric loop groups and free loop spaces are important both 
in mathematics and in modern physical theories. Moreover,
generalized geometric loop groups also can be used in the same fields
of sciences
and open new opportunities. In the cohomology theory and physical
applications stochastic processes on the free loop spaces are used
\cite{eljm,jole,malb}.
In these papers were considered only particular cases of 
real free loop spaces and groups for finite dimensional manifolds $N$
for mappings from $S^1$ into $N$.
No any applications to the representation theory were given there.
\par On the other hand, representation theory of non-locally compact groups
is little developed apart from the case of locally compact groups
(see \cite{barut,fell} and references therein).
In particular, geometric loop and diffeomoprphism groups have important 
applications in modern physical theories (see \cite{ish,mensk}
and references therein).
Groups of loops and groups of diffemorphisms are also intesively
used in gauge theory.
Loop groups defined below with the help of
families of mappings from a manifold $M$ into another manifold $N$
with a dimension $dim (M) >1$ can be used in the membrane theory
which is the generalization of the string (superstring) theory.
\par One of the main tools in the investigation of unitary represenations 
of nonlocally compact groups are quasi-invariant measures.
In previous works of the author \cite{lurim1,lurim2} 
Gaussian quasi-invariant measures were constructed on diffeomorphism groups
with some conditions on real manifolds. For example, compact manifolds 
without boundary were not considered, as well as infinite dimensional 
manifolds with boundary. In this article new Gevrey-Sobolev classes 
of smoothness for diffeomorphism groups of infinite dimensional 
real, complex and quaternion manifolds 
are defined and investigated. This permits to define on them the Hilbert 
manifold structure. This in its turn simplifies the construction of 
stochastic processes and transition quasi-invariant probabilities on them.
Transition probabilities are constructed below, which
are quasi-invariant, for
wider classes of manifolds. Pairs of topological groups $G$ and their dense
subgroups $G'$ are described precisely.
Then measures are used for the study of associated unitary 
regular and induced representations of dense subgroups $G'$. 
\par Section 2 is devoted to the definitions of topological and 
manifold structures of loop groups and diffeomorphism groups
and their semidirect products and their dense subgroups.
For this necessary statements about structures of quaternion
manifolds are proved in Propositions 2.1.3.3, 4, Theorems 2.1.3.6, 7, 9,
Lemmas 2.1.6, 2.1.6.2.
The existence of these groups is proved and that they are
infinite dimensional Lie groups not satisfying even locally
the Campbell-Hausdorff formula besides the degenerate case of
groups of holomorphic diffeomorphisms of compact complex manifolds
(see Theorems 2.1.4.1, 2.1.7, 2.2.1, 2.8, 2.9, 2.10, 2.11,
Lemmas 2.5, 2.5.1, 2.6.2, 2.7). In the cases of complex
and quaternion manifolds it is proved that they have structures
of complex and quaternion manifolds respectively.
Their structure as manifolds and groups is studied not only for
orientable manifolds, but also for nonorientable manifolds
$M$ or $N$ over $\bf R$ or $\bf H$ (see Theorem 2.1.8).
In Section 3 transition quasi-invariant
differentiable probabilities are studied (see Theorem 3.3). 
Unitary representations of dense subgroups $G'$ 
founded in Sections 2 and 3 are investigated in Section 4.
All objects given in Sections 2-3 were
not considered by other authors, besides very specific particular 
cases of the diffeomorphism group of real and complex finite
dimensional manifolds and loop groups for $M=S^1$ outlined above.
In Section 4 unitary representations including topologically
irreducible of semidirect products and constituing them subgroups
are investigated.
It is proved, that this procedure provides a family of the cardinality
$card ({\bf R})$ of pairwise nonequivalent irreducible unitary
representations. A differentiability of such representations is studied.
Then a usefulness of differentiable representations is illustrated
for construction of representations of the corresponding algebras.
\par All results of this paper are obtained for the first time.

\section{Semidirect products of groups of loops and
groups of diffeomorphisms of finite and infinite dimensional manifolds.}
To avoid misunderstandings we first give our definitions
of manifolds considered here and then of loop and diffeomorphism groups.
\par {\bf 2.1.1. Remark.} An atlas $At(M)=\{ (U_j,\phi _j): j \} $ 
of a manifold $M$ on a Banach space $X$ over $\bf R$
is called uniform, if its charts satisfy 
the following conditions: \\
$(U1)$ for each $x\in G$ there exist
neighbourhoods $U_x^2\subset U_x^1\subset U_j$
such that for each $y\in U_x^2$ there is the inclusion
$U_x^2\subset U_y^1$; \\
$(U2)$ the image $\phi _j(U_x^2)\subset X$
contains a ball of the fixed positive radius
$\phi _j(U_x^2)\supset B(X,0,r):=\{ y: y\in X, \| y\| \le r \} ;$ \\
$(U3)$ for each pair of intersecting charts
$(U_1,\phi _1)$ and $(U_2,\phi _2)$ connecting mappings
$F_{\phi _2,\phi _1}=\phi _2\circ \phi _1^{-1}$
are such that $\sup_x \| F'_{\phi _2,\phi _1}(x) \|
\le  C$ and $\sup_x \| F'_{\phi _1,\phi _2}(x) \|
\le  C$, where $C=const >0$ does not depend on  
$\phi _1$ and $\phi _2$. For the diffeomorphism group
$Diff^t_{\beta ,\gamma }(M)$ and loop groups $(L^MN)_{\xi }$
we also suppose that manifolds satisfy conditions
of \cite{ludan,lurim1,lurim2} such that these groups are separable,
but here let $M$ and $N$ may be with a boundary,
where \\
$(N1)$ $N$ is of class not less, than (strongly) $C^{\infty }$ and such that
$sup_{x\in S_{j,l}} \| F^{(n)}_{\psi _j,\psi _l}(x) \| \le C_n$ for each 
$0\le n\in \bf Z$, when $V_{j,l}\ne \emptyset $, $C_n>0$ are constants,
$At(N):= \{ (V_j,\psi _j): j \} $ denotes an atlas of $N$, $V_{j,l}:=
V_j\cap V_l$ are intersections of charts, $S_{j,l}:=\psi _l(V_{j,l})$,
$\bigcup_jV_j=N$. 
\par Conditions $(U1-U3,N1)$ are supposed to be satisfied for the 
manifold $N$ for each loop group, as well as for the manifold $M$ for 
each diffeomorphism group. Certainly, 
classes of smoothness of manifolds are supposed to be not less 
than that of groups.
\par {\bf 2.1.2.1. Definition.} 
~ A canonical closed subset $Q$ of the Euclidean space
$X=\bf R^n$ or of the standard separable
Hilbert space $X=l_2({\bf R})$ over $\bf R$
is called a quadrant if it can be given by the condition $Q:=\{ x\in X:
q_j(x)\ge 0 \} $, where $(q_j: j\in \Lambda _Q)$ 
are linearly independent elements of the topologically adjoint space $X^*$.
Here $\Lambda _Q\subset \bf N$ (with $card (\Lambda _Q)=k\le n$
when $X=\bf R^n$) and $k$ is called the index of $Q$.
If $x\in Q$ and exactly $j$ of the $q_i$'s
satisfy $q_i(x)=0$ then $x$ is called a corner of index $j$.  
\par If $X$ is an additive group and also left and right module over
$\bf H$ with the corresponding associativity and distributivity
laws then it is called the vector space over $\bf H$.
In particular $l_2 ({\bf H})$ consisting of all sequences
$x = \{ x_n\in {\bf H}: n \in {\bf N} \} $ with the finite norm
$\| x \| <\infty $ and scalar product $(x,y):=\sum_{n=1}^{\infty }
x_ny_n^*$ with $\| x \| := (x,x)^{1/2}$ is called the Hilbert
space (of separable type) over $\bf H$, where $z^*$ denotes
the conjugated quaternion, $zz^* =: |z|^2$, $z\in \bf H$.
Since the unitary space $X=\bf C^n$ or the separable Hilbert space
$l_2({\bf C})$ over $\bf C$ or the quaternion space $X=\bf H^n$ or the
separable Hilbert space $l_2({\bf H})$ over $\bf H$ while considered
over the field $\bf R$
is isomorphic with $X_{\bf R}:=\bf R^{2n}$ or $l_2({\bf R})$
or $\bf R^{4n}$ or $l_2 ({\bf R})$ respectively,
then the above definition also
describes quadrants in $\bf C^n$ and $l_2({\bf C})$ and $\bf H^n$
and in $l_2 ({\bf H})$.
In the latter case we also consider generalized quadrants
as canonical closed subsets which can be given by
$Q:=\{ x\in X_{\bf R}:$ $q_j(x+a_j)\ge 0, a_j\in X_{\bf R},
j\in \Lambda _Q \} ,$ where $\Lambda _Q\subset \bf N$
($card(\Lambda _Q)=k\in \bf N$ when $dim_{\bf R}X_{\bf R}<\infty $).
\par {\bf 2.1.2.2. Notation.} If for each open subset $U\subset Q\subset X$ 
a function $f:  Q\to Y$ for Banach spaces $X$ and $Y$ over $\bf R$
has continuous Frech\'et differentials $D^{\alpha }f|_U$ on $U$ 
with $\sup_{x\in U} \| D^{\alpha }f(x) \|_{L(X^{\alpha },Y)} <\infty $ 
for each $0\le \alpha \le r$ for an integer $0\le r$ or $r=\infty $, 
then $f$ belongs to the class of smoothness $C^r(Q,Y)$,
where $0\le r\le \infty $,
$L(X^k,Y)$ denotes the Banach space of bounded $k$-linear
operators from $X$ into $Y$.
\par {\bf 2.1.2.3. Definition.} 
A differentiable mapping $f:  U\to U'$ is called a diffeomorphism if 
\par $(i)$ $f$ is bijective and there exist continuous mappings
$f'$ and $(f^{-1})'$, where $U$ and $U'$ 
are interiors of quadrants $Q$ and $Q'$ in $X$. 
\par In the complex case and in the quaternion case we consider
bounded generalized quadrants $Q$ and $Q'$ in $\bf C^n$ or $l_2({\bf C})$
and in $\bf H^n$ or $l_2({\bf H})$ respectively
such that they are domains with piecewise 
$C^{\infty }$-boundaries. We impose additional conditions on
the diffeomorphism $f$ in the complex case:
\par $(ii)$ ${\bar \partial }f=0$ on $U$, 
\par $(iii)$ $f$ and all its strong (Frech\'et) differentials (as multilinear
operators) are bounded on $U$, where $\partial f$ and ${\bar \partial }f$
are differential $(1,0)$ and $(0,1)$ forms respectively,
$d=\partial +{\bar \partial }$ is an exterior derivative.
In particular, for 
$z=(z^1,...,z^n)\in \bf C^n$, $z^j\in \bf C$, $z^j=x^{2j-1}+ix^{2j}$
and $x^{2j-1}, x^{2j}\in \bf R$ for each $j=1,...,n,$
$i=(-1)^{1/2}$, there are expressions:
$\partial f:=\sum_{j=1}^n(\partial f/\partial z^j)dz^j$,
${\bar \partial }f:=\sum_{j=1}^n(\partial f/\partial {\bar z}^j)d{\bar z}^j$.
In the infinite dimensional case there are equations:
$(\partial f)(e_j)=\partial f/\partial z^j$
and $({\bar \partial }f)(e_j)=\partial f/\partial {\bar z}^j$, 
where $\{ e_j: {j\in \bf N} \} $ is the standard orthonormal base
in $l_2({\bf C})$, $\partial f/\partial z^j=(\partial f/\partial x^{2j-1}
-i\partial f/\partial x^{2j})/2$,
$\partial f/\partial {\bar z}^j=(\partial f/\partial x^{2j-1}
+i\partial f/\partial x^{2j})/2$.
\par In the quaternion case consider quaternion holomorphic
diffeomorphisms $f$:
\par $(iv)$ ${\tilde \partial }f=0$ on $U$,
\par $(v)$ $f$ and all its superdifferentials (as
$\bf R$-multilinear $\bf H$-additive operators) are bounded on $U$,
where $\partial f$ and ${\tilde \partial }f$
are differential $(1,0)$ and $(0,1)$ forms respectively,
$d=\partial +{\tilde \partial }$ is an exterior derivative,
$\partial $ corresponds to superdifferentiation by $z$ and
${\tilde \partial }$ corresponds to superdifferentiation by
${\tilde z}:=z^*$, $z\in U$ (see \cite{luoyst}).
\par The Cauchy-Riemann Condition $(ii)$ in the complex case
or $(iv)$ in the quaternion case means that $f$ on $U$
is the complex holomorphic or quaternion holomorphic mapping
respectively.
\par {\bf 2.1.2.4. Definition and notation.}
A complex manifold or a quaternion manifold $M$ with corners is
defined in the usual way:  it is a metric separable space
modelled on $X=\bf C^n$ or $X=l_2({\bf C})$ or on $X=\bf H^n$
or $X=l_2({\bf H})$ respectively
and is supposed to be of class $C^{\infty }$.  Charts on $M$ are
denoted $(U_l, u_l, Q_l)$, that is, $u_l:  U_l\to u_l(U_l) \subset Q_l$ is
a $C^{\infty }$-diffeomorphism for each $l$,
$U_l$ is open in $M$, $u_l\circ {u_j}^{-1}$ is biholomorphic from the
domain $u_j(U_l\cap U_j)\ne \emptyset $ onto $u_l(U_l\cap U_j)$ 
(that is, $u_j\circ u_l^{-1}$ and 
$u_l\circ u_j^{-1}$ are holomorphic and bijective)
and $u_l\circ u_j^{-1}$ 
satisfy conditions $(i-iii)$ or $(i,iv,v)$ correspondingly
from \S 2.1.2.3, $\bigcup_jU_j=M$.
\par A point $x\in M$ is called a corner of index $j$
if there exists a chart $(U,u,Q)$ of $M$ with $x\in U$ and $u(x)$ is of index
$ind_M(x)=j$ in $u(U)\subset Q$. A set of all corners of index $j\ge 1$ is
called a border $\partial M$ of $M$, $x$ is called an inner point of $M$ if
$ind_M(x)=0$, so $\partial M=\bigcup_{j\ge 1}\partial ^jM$, where
$\partial ^jM:=\{ x\in M:  ind_M(x)=j \} $.  
\par For a real manifold with corners on the connecting mappings
$u_l\circ u_j^{-1}\in C^{\infty }$ 
of real charts is imposed only Condition $2.1.2.3(i)$.
\par {\bf 2.1.2.5. Definition of a submanifold with corners.} 
A subset $Y\subset M$ is called a 
complex or quaternion submanifold with corners of $M$ if for each
$y\in Y$ there exists a chart $(U,u,Q)$ of $M$ 
centered at $y$ (that is $u(y)=0$) and there exists a
quadrant $Q'\subset {\bf C^k}$ or in $l_2({\bf C})$
or $Q'\subset {\bf H^k}$ or $Q'\subset l_2({\bf H})$ respectively such that
$Q'\subset Q$ and $u(Y\cap U)=u(U)\cap Q'$.
A submanifold with corners $Y$ of
$M$ is called neat, if an index in $Y$ of each $y\in Y$ coincides with its
index in $M$. 
\par  Analogous definitions are for real manifolds with corners
for $\bf R^k$ and $\bf R^n$ or $l_2({\bf R})$ instead of $\bf C^k$
and $\bf C^n$ or $l_2({\bf C})$.
\par {\bf 2.1.2.6. Terminology.} 
Henceforth, the term a complex manifold or a quaternion manifold
$N$ modelled on $X=\bf C^n$ or $X=l_2({\bf C})$ or on
$X=\bf H^n$ or $X=l_2({\bf H})$ means a metric separable
space supplied with an atlas $\{ (U_j,\phi _j): j\in \Lambda _N \} $
such that:
\par $(i)$ $U_j$ is an open subset of $N$ for each $j\in \Lambda _N$
and $\bigcup_{j\in \Lambda _N}U_j=N$, where $\Lambda _N\subset \bf N$;
\par $(ii)$ $\phi _j: U_j\to \phi _j(U_j)\subset X$ is a
$C^{\infty }$-diffeomorphism for each $j$, where $\phi _j(U_j)$ is a 
$C^{\infty }$-domain in $X$;
\par $(iii)$ $\phi _j\circ \phi _m^{-1}$ is a complex biholomorphic
or quaternion biholomorphic mapping respectively
from $\phi _m(U_m\cap U_j)$ onto
$\phi _j(U_m\cap U_j)$ while $U_m\cap U_j\ne \emptyset $. 
When $X=l_2({\bf C})$ or $X=l_2({\bf H})$ it is supposed, that
$\phi _j\circ \phi _m^{-1}$ is Frech\'et (strongly)
$C^{\infty }$-differentiable for each $j$ and $m$
and certainly either condition $(ii)$ or $(iv)$ of \S 2.1.2.3
is satisfied.
\par {\bf 2.1.3.1. Remark.}  Let $X$ be either the standard separable
Hilbert space $l_2=l_2({\bf C})$ over the field $\bf C$ of complex numbers
or $X=\bf C^n$ or $l_2=l_2({\bf H})$ over the skew field of quaternions
or $X=\bf H^n$ correspondingly.  Let $t\in \bf
N_o$ $:={\bf N}\cup \{ 0\}$, ${\bf N}:=\{ 1,2,3,... \}$ 
and $W$ be a domain with a continuous piecewise
$C^{\infty }$-boundary $\partial W$ in $\bf R^{2m}$ or in $\bf R^{4m}$
respectively, $m\in \bf N$, that is, $W$ is a $C^{\infty }$-manifold
with corners and it is a canonical closed
subset of $\bf C^m$ or of $\bf H^m$, $cl (Int (W))=W$, where $cl (V)$
denotes the closure of $V$, $Int (V)$ denotes the interior of $V$
in the corresponding topological space.  As
usually $H^t (W,X)$ denotes the Sobolev space of functions 
$f:  W\to X$ for which there exists a finite norm
\par $\| f\|_{H^t(W,X)} := (\sum_{|\alpha |\le t}{\|
D^{\alpha }f\|^2}_{L^2 (W,X)})^{1/2}<\infty $, \\
where $f(x)=(f^j(x):  j\in {\bf
N})$, $f(x)\in l_2$, $f^j(x)\in \bf C$ or $f^j(x)\in \bf H$ correspondingly,
$x\in W$,
\par ${\| f\|^2}_{L^2
(W,X)}:=\int_W {\| f(x)\|^2}_X\lambda (dx)$, $\lambda $ is the Lebesgue
measure on $\bf R^{2m}$ or on $\bf R^{4m}$ respectively,
$\| z\|_{l_2}:= (\sum_{j=1}^{\infty }|z^j|^2)^{1/2}$, $z=(z^j:
j\in {\bf N})\in l_2$, $z^j\in \bf C$ or $z^j\in \bf H$.
Then $H^{\infty }(W,X):=
\bigcap_{t\in \bf N}H^t(W,X)$ is a uniform space with a uniformity
given by the family of norms $\{ \| f \|_{H^t(W,X)}: t\in {\bf N} \}$.
\par {\bf 2.1.3.2. Sobolev spaces for manifolds.}  Let now $M$ be a
compact Riemannian or complex or quaternion
$C^{\infty }$-manifold with corners with 
a finite atlas $At(M):=\{ (U_i, \phi _i,
Q_i); i\in \Lambda _M \} $, where $U_i$ is open in $M$ for each $i$, 
$\phi _i:  U_i\to \phi
_i(U_i)\subset Q_i\subset \bf R^m$ (or it is a subset in $\bf C^m$ or in
$\bf H^m$ correspondingly)
is a diffeomorphism (in addition holomorphic respectively as
in \S 2.1.2.3), $(U_i, \phi _i)$ is a chart, $i\in \Lambda _M \subset \bf N$.
\par Let also $N$ be a separable real or complex or quaternion
metrizable manifold with corners
modelled either on $X=\bf R^n$ or $X=l_2({\bf R})$ or
on $X=\bf C^n$ or $X=l_2({\bf C})$ or on $X=\bf H^n$
or $l_2({\bf H})$ respectively. Let $(V_i, \psi _i, S_i)$ be
charts of an atlas $At(N):=\{ (V_i, \psi _i, S_i):  i\in \Lambda _N \} $
such that $\Lambda _N\subset \bf N$ and
$\psi _i:  V_i\to \psi _i(V_i)\subset
S_i\subset X$ is a diffeomorphism for each $i$,
$V_i$ is open in $N$, $\bigcup_{i\in \Lambda _N}V_i=N$.
We denote by $H^t(M,N)$ the Sobolev
space of functions $f:  M\to N$ for which $f_{i,j}\in H^t(W_{i,j},X)$ for each
$j\in \Lambda _M$ and $i\in \Lambda _N$ for a domain $W_{i,j}\ne \emptyset $
of $f_{i,j}$, where $f_{i,j}:=\psi _i\circ f\circ {\phi _j}^{-1}$, and
$W_{i,j}=\phi _j(U_j\cap f^{-1}(V_i))$ is a canonical closed subset
of $\bf R^m$ (or $\bf C^m$ or $\bf H^m$ respectively). The
uniformity in $H^t(M,N)$ is given by the base $\{ (f,g)\in
(H^t(M,N))^2:  \sum_{i\in \Lambda _N, j\in \Lambda _M} {\| f_{i,j}-
g_{i,j}\|^2}_{H^t(W_{i,j}, X)}<\epsilon \}$, where $\epsilon >0$, $W_{i,j}$
is a domain of $(f_{i,j}- g_{i,j}).$ For $t=\infty $ as usually
$H^{\infty }(M,N):=\bigcap_{t\in \bf N} H^t(M,N)$.
\par {\bf 2.1.3.3. Proposition.} {\it Let $M$ be a quaternion manifold.
Then there exists a tangent bundle $TM$ which has the structure
of the quaternion manifold such that each fibre $T_xM$ is the vector
space over $\bf H$.}
\par {\bf Proof.} With the help of complex $2\times 2$ matrices
present each quaternion $z$, $z={{t\quad u} \choose { - {\bar u}
\quad {\bar t} }}$, where $t, u\in \bf C$, hence each $z\in \bf H$
is also a real $4\times 4$ matrix. Therefore, $M$ has also the
structure of real manifold. Since each quaternion holomorphic mapping
is infinite differentiable (see Theorems 2.15 and 3.10 \cite{luoyst}),
then there exists its tangent bundle $TM$ which is $C^{\infty }$-manifold
such that each fibre $T_xM$ is a tangent space, where $x\in M$,
$T$ is the tangent functor.
If $At (M)= \{ (U_j,\phi _j): j \} $, then $At (TM)= \{
(TU_j, T\phi _j): j \} $, $TU_j=U_j\times X$, where $X$ is the quaternion
vector space on which $M$ is modelled,
$T(\phi _j\circ \phi _k^{-1})=(\phi _j\circ \phi _k^{-1},
D(\phi _j\circ \phi _k^{-1}))$ for each $U_j\cap U_k\ne \emptyset $.
Each transition mapping $\phi _j\circ \phi _k^{-1}$ is quaternion
holomorphic on its domain, then its (strong) differential coincides
with the superdifferential $D(\phi _j\circ \phi _k^{-1})=
D_z(\phi _j\circ \phi _k^{-1})$, since ${\tilde \partial }
(\phi _j\circ \phi _k^{-1})=0$. Therefore, 
$D(\phi _j\circ \phi _k^{-1})$ is $\bf R$-linear and
$\bf H$-additive, hence is the automorphism of the
quaternion vector space $X$. Since $D_z(\phi _j\circ \phi _k^{-1})$
is also quaternion holomorphic, then $TM$ is the quaternion manifold.
\par {\bf 2.1.3.4. Proposition.} {\it Let $M$ be a quaternion
paracompact manifold on $X=\bf H^n$ or on $X=l_2({\bf H})$,
where $n\in \bf N$. Then $M$ can be supplied with the quaternion
Hermitian metric.}
\par {\bf Proof.} Mention that a continuous mapping $\pi $
of a Hausdorff space $A$ into another $B$ is called an $\bf H$-vector
bundle, if
\par $(i)$ $A_x:=\pi ^{-1}(x)$ is a vector space $X$ over $\bf H$ for each
$x\in B$ and
\par $(ii)$ there exists a neighbourhood $U$ of $x$ and a homeomorphism
$\psi : \pi ^{-1} (U)\to U\times X$ such that $\psi (A_x)=
\{ x \} \times X$ and $\psi _x: A_x\to X$ is a $\bf H$-vector isomorphism
for each $x\in B$, where $\psi ^x:=\psi \circ \pi _2$, $\pi _2:
U\times X\to X$ is the projection, $(U,\psi )$ is called a local
trivialization.
\par An $\bf H$-Hermitian inner product in a vector space
$X$ over $\bf H$ is a $\bf R$-bilinear $\bf H$-biadditive form
$<*,*> : X^2\to \bf H$ such that:
\par $(iii)$ $<x,x> \ge 0$ for each $x\in X$,
$<x,x>=0$ if and only if $x=0$;
\par $(iv)$ $<x,y>^* = <y,x>$ for each $x, y\in X$;
\par $(v)$ $<x,y+z>=<x,y>+<x,z>$, $<x+y,z>=<x,z>+<y,z>$
for each $x, y, z\in X$;
\par $(vi)$ $<ax,y>=a<x,y>$ and $<x,by>=<x,y>b^*$ for each
$x, y\in X$ and $a, b\in \bf H$.
\par If $\pi : A\to B$ is an $\bf H$-vector bundle, then
an $\bf H$-Hermitian metric $g$ on $A$ is an assignment of an
$\bf H$-Hermitian inner product $<*,*>_x$ to each fibre $A_x$
such that for each open subset $U$ in $B$ and $\xi ,\eta \in
C^{\infty }(U,A)$ the mapping $<\xi ,\eta >: U\to \bf H$ such that
$<\xi , \eta >(x)=<\xi (x), \eta (x)>_x$ is $C^{\infty }$.
An $\bf H$-vector bundle $A$ equipped with an $\bf H$-Hermitian metric
$g$ we call an $\bf H$-Hermitian vector bundle.
If $A$ is paracompact, then each its (open) covering contains a locally
finite refinement (see \cite{eng}). Take a locally finite covering
$\{ U_j: j \} $ of $B$. By the supposition of this proposition
$X=\bf H^n$ or $X=l_2({\bf H})$. Choose a subordinated real partition
of unity $\{ \alpha _j: j \} $ of class $C^{\infty }$:
$\quad \sum_j\alpha _j(x)=1$ and $\alpha _j(x)\ge 0$ for
each $j$ and each $x\in B$. Therefore, there exists a frame
$\{ e_k: k \in \Lambda \} $ at $x\in B$, where either
$\Lambda = \{ 1,2,...,n \} $ or $\Lambda =\bf N$, such that
$e_k\in C^{\infty }(U_j,A)$ for each $k$, $\{ e_k: k \in \Lambda \} $
are $\bf H$-linearly independent at each $y\in U_j$ relative to left
and right mulitplications on constants $a_k, b_k$ from $\bf H$,
that is, $\sum_k a_ke_kb_k=0$ if and only if $\sum_k |a_kb_k|=0$,
where $x\in U_j$.
Define $<\xi ,\eta >^j_x:=\sum_k \xi (e_k)(x) (\eta (e_k)(x))^*$
and $<\xi ,\eta >_x := \sum_j \alpha _j(x) <\xi ,\eta >^j_x$.
If $\xi ,\eta \in C^{\infty }(U,A)$, then the mapping
$x\mapsto <\xi (x), \eta (x)>_x$ is $C^{\infty }$ on $U$.
Since $|<\xi ,\eta >_x|\le \sum_j \alpha _j(x) (<\xi ,\xi >^j_x)^{1/2}
(<\eta ,\eta >^j_x)^{1/2}<\infty $ for each $x\in B$, hence
the $\bf H$-Hermitian metric is correctly defined.
For each $b\in \bf R$ we have $<\xi b, \eta b>_x=|b|^2<\xi ,\eta >_x$,
since $\bf R$ is the centre of the algebra $\bf H$ over $\bf R$.
Take in particular $\pi : TM\to M$ and this provides an $\bf H$-Hermitian
metric in $M$.
\par {\bf 2.1.3.5. Definitions.} A $C^1$-mapping $f: M\to N$ is called
an immersion, if $rang (df|_x: T_xM\to T_{f(x)}N) = m_M$ for each $x\in M$,
where $m_M := dim_{\bf R}M$.
An immersion $f: M\to N$ is called an embedding, if $f$ is bijective.
\par {\bf 2.1.3.6. Theorem.} {\it Let $M$ be a compact
quaternion manifold, $dim_{\bf H}M=m<\infty $. Then there exists a
quaternion holomorphic embedding $\tau : M\hookrightarrow {\bf H}^{2m+1}$
and a quaternion holomorphic immersion
$\theta : M\to {\bf H}^{2m}$ correspondingly. Each continuous mapping
$f: M\to {\bf H}^{2m+1}$ or $f: M\to {\bf H}^{2m}$ can be approximated
by $\tau $ or $\theta $ relative to the norm $\| * \|_{C^0}$.}
\par {\bf Proof.} Since $M$ is compact, then it is finite
dimensional over $\bf H$, $dim_{\bf H}M=m\in \bf N$, such that
$dim_{\bf R}M=4m$ is its real dimension. Choose an atlas $At' (M)$
refining initial atlas $At (M)$ of $M$ such that $({U'}_j, \phi _j)$
are charts of $M$, where each ${U'}_j$ is quaternion holomorphic
diffeomorphism to an interior of the unit ball $Int (B({\bf H^k},0,1))$,
where $B({\bf H^k},y,r) := \{ z\in {\bf H^k}: |z-y|\le r \} $.
In view of compactness of $M$ a covering $\{ {U'}_j, j \} $
has a finite subcovering, hence $At' (M)$ can be chosen finite.
Denote for convenience the latter atlas as $At (M)$.
Let $(U_j, \phi _j)$ be the chart of the atlas
$At (M)$, where $U_j$ is open in $M$, hence $M\setminus U_j$
is closed in $M$.
\par  Consider the space ${\bf H^k}\times {\bf R}$
as the $\bf R$-linear space $\bf R^{4k+1}$. The unit sphere
$S^{4k}:=S ({\bf R}^{4k+1},0,1) := \{ z\in {\bf R}^{4k+1}:$ $|z|=1 \} $
in ${\bf H^k}\times \bf R$ can be supplied with two charts
$(V_1, \phi _1)$ and $(V_2, \phi _2)$ such that $V_1:=
S^{4k}\setminus \{ 0,...,0, 1 \} $ and
$V_2:=S^{4k}\setminus \{ 0,...,0, - 1 \} $,
where $\phi _1$ and $\phi _2$ are stereographic projections
from poles $\{ 0,...,0, 1 \} $ and $ \{ 0,...,0, -1 \} $
of $V_1$ and $V_2$ respectively onto $\bf H^k$.
Since $z^* = - (z + JzJ+ KzK+ LzL)/2$ in $\bf H^k$,
then $\phi _1\circ \phi _2^{-1}$ (in the $z$-representation)
is quaternion holomorphic
diffeomorphism of ${\bf H^k}\setminus \{ 0 \} $, but certainly
with neither right nor left superlinear superdifferential
$D_z (\phi _1\circ \phi _2^{-1})$, where $J$, $K$, $L$ are complex Pauli
$2\times 2$-matrices.
Thus $S^{4k}$ can be supplied with the structure of the quaternion manifold.
\par Therefore, there exists a quaternion holomorphic mapping
$\psi _j$ (that is, locally $z$-analytic \cite{luoyst}) of $M$ into the unit
sphere $S^{4m}$ such that $\psi _j (M\setminus U_j) = \{ x_j \} $ is the
singleton and $\psi _j: U_j\to \psi _j(U_j)$ is the quaternion holomorphic
diffeomorphism onto the open subset $\psi _j(U_j)$ in
$S^{4m}$, $x_j\in S^{4m}\setminus \psi _j(U_j)$. There is evident embedding
of ${\bf H^m}\times \bf R$ into $\bf H^{m+1}$.
Then the mapping $\psi (z):=(\psi _1(z),...,\psi _n(z))$
is the embedding into $(S^{4m})^n$ and hence into
${\bf H}^{n{m+1}}$, since the rank
$rank [d_z\psi (z)]=4m$ at each point $z\in M$, because
$rank [d_z\psi _j(z)]=4m$ for each $z\in U_j$ and $dim_{\bf H}\psi (U_j)\le
dim_{\bf H}M=m$. Moreover, $\psi (z)\ne \psi (y)$ for each
$z\ne y\in U_j$, since $\psi _j(z)\ne \psi _j(y)$.
If $z\in U_j$ and $y\in M\setminus U_j$, then there exists
$l\ne j$ such that $y\in U_l\setminus U_j$,
$\psi _j(z)\ne \psi _j(y)=x_j$.
\par Let $M\hookrightarrow {\bf H}^N$ be the
quaternion holomorphic embedding as above.
There is also the quaternion holomorphic embedding of $M$
into $(S^{4m})^n$ as it is shown above, where $(S^{4m})^n$ is the quaternion
manifold as the product of quaternion manifolds.
Consider the bundle of all $\bf H$ straight
lines ${\bf H}x$ in ${\bf H}^N$, where $x\in {\bf H}^N$.
They compose the projective space
${\bf H}P^{N-1}$. Fix the standard orthonormal base
$ \{ e_1,...,e_N \} $ in ${\bf H}^N$ and projections
on $\bf H$-linear subspaces relative to this base
$P^L(x):=\sum_{e_j\in L}x_je_j$ for the $\bf H$-linear span
$L=span_{\bf H} \{ e_i:$ $i\in \Lambda _L \} $,
$\Lambda _L\subset \{ 1,...,N \} $, where
$x=\sum_{j=1}^Nx_je_j$, $x_j\in \bf H$ for each $j$,
$e_j=(0,...,0,1,0,...,0)$ with $1$ at $j$-th place.
In this base consider the $\bf H$-Hermitian scalar product
$<x,y> := \sum_{j=1}^Nx_jy^*_j$.
Let $l\in {\bf H}P^{N-1}$, take a $\bf H$-hyperplane
denoted by ${\bf H}^{N-1}_l$ and given by the condition:
$<x,[l]>=0$ for each $x\in {\bf H}^{N-1}_l$, where $0\ne [l]\in {\bf H}^N$
characterises $l$. Take $\| [l] \| =1$. Then the orthonormal base
$\{ q_1,...,q_{N-1} \} $ in ${\bf H}^{N-1}_l$ and together with $[l]=:q_N$
composes the orthonormal base $\{ q_1,...,q_N \} $ in ${\bf H}^N$.
This provides the quaternion holomorphic projection
$\pi _l: {\bf H}^N\to {\bf H}^{N-1}_l$
relative to the orthonormal base $ \{ q_1,...,q_N \} $.
The operator $\pi _l$ is $\bf H$ left and also right linear
(but certainly nonlinear relative to $\bf H$), hence $\pi _l$
is quaternion holomorphic.
\par To construct an immersion it is sufficient, that each
projection $\pi _l: T_xM\to {\bf H}^{N-1}_l$ has $ker [d(\pi _l(x))]=
\{ 0 \} $ for each $x\in M$. The set of all $x\in M$ for which
$ker [d(\pi _l(x))] \ne \{ 0 \} $
is called the set of forbidden directions of the first kind.
Forbidden are those and only those directions $l\in {\bf H}P^{N-1}$
for which there exists $x\in M$ such that $l'\subset T_xM$,
where $l'=[l]+z$, $z\in {\bf H}^N$. The set of all forbidden directions
of the first kind forms the quaternion manifold $Q$
of quaternion dimension $(2m-1)$ with points $(x,l)$, $x\in M$,
$l\in {\bf H}P^{N-1}$, $[l]\in T_xM$. Take $g: Q\to {\bf H}P^{N-1}$
given by $g(x,l):=l$. Then $g$ is quaternion holomorphic.
\par In each quaternion paracompact manifold $A$ modelled
on $\bf H^p$ there exists an $\bf H$-Hermitian metric
(see Proposition 2.1.3.4), hence it can be supplied with the Riemann
manifold structure also. Therefore, on $A$ there exists a Riemann
volume element. In view of the Morse theorem
$\mu (g(Q))=0$, if $N-1>2m-1$, that is, $2m<N$, where $\mu $ is the
Riemann volume element in ${\bf H}P^{N-1}$. In particular,
$g(Q)$ is not contained in ${\bf H}P^{N-1}$ and there exists
$l_0\notin g(Q)$, consequently, there exists $\pi _{l_0}: M\to
{\bf H}^{N-1}_{l_0}$.
This procedure can be prolonged, when $2m<N-k$, where $k$
is the number of the step of projection. Hence $M$ can be immersed
into ${\bf H}^{2m}$.
\par Consider now the forbidden directions of the second type:
$l\in {\bf H}P^{N-1}$, for which there exist $x\ne y\in M$
simultaneously belonging to $l$ after suitable parrallel translation
$[l]\mapsto [l]+z$, $z\in {\bf H}^N$. The set of the forbidden directions
of the second type forms the manifold $\Phi :=M^2\setminus \Delta $,
where $\Delta := \{ (x,x):$ $x\in M \} $. Consider $\psi :
\Phi \to {\bf H}P^{N-1}$, where $\psi (x,y)$ is the straight $\bf H$-line
with the direction vector $[x,y]$ in the orthonormal base.
Then $\mu (\psi (\Phi ))=0$
in ${\bf H}P^{N-1}$, if $2m+1<N$. Then the closure
$cl (\psi (\Phi ))$ coinsides with $\psi (\Phi )\cup g(Q)$ in
${\bf H}P^{N-1}$.
Hence there exists $l_0\notin cl (\psi (\Phi ))$. Then consider
$\pi _{l_0}: M\to {\bf H}_{l_0}^{N-1}$.
This procedure can be prolonged, when $2m+1<N-k$, where $k$
is the number of the step of projection. Hence $M$ can be embedded
into ${\bf H}^{2m+1}$.
\par {\bf 2.1.3.7. Theorem.} {\it Let $M$ be a quaternion
locally compact manifold with and atlas $At (M)= \{ (U_j,\phi _j) : j \} $
such that transition mappings $\phi _j\circ \phi _k^{-1}$ of
charts with $U_j\cap U_k\ne \emptyset $ 
are either right or left superlinearly superdifferentiable. Then
$M$ is orientable.}
\par {\bf Proof.} Since $M$ is locally compact, then $M$ is finite
dimensional over $\bf H$ such that $dim_{\bf H}X=m<\infty $,
where $X=T_xM$ for $x\in M$. Consider open subsets $U$ and $V$
in $X$ and a function $\phi : U\to V$ which is right superdifferentiable.
Write $\phi $ in the form $\phi = (\mbox{ }_1\phi ,..., \mbox{ }_m\phi )$,
where $\mbox{ }_v\phi \in \bf H$ for each $v=1,...,m$.
Then $\mbox{ }_v\phi =\mbox{ }_v\alpha + \mbox{ }_v\beta j$, where
$\mbox{ }_v\alpha $ and $\mbox{ }_v\beta \in \bf C$, $i, j, k$
are generators of $\bf H$ such that $i^2=j^2=k^2=-1$, $ij=-ji=k$,
we write also $\mbox{ }_vz=\mbox{ }_vx+\mbox{ }_vyj$ with
$\mbox{ }_vx$ and $\mbox{ }_vy\in \bf C$, where $z\in \bf H^m$.
From Proposition 2.2 \cite{luoyst} it follows, that
there exist right superlinear
$\partial a /\partial b$ with $\partial a/\partial {\bar b}=0$
for each $a\in \{ \mbox{ }_v\alpha , \mbox{ }_v\beta : v=1,...,m \} $
and each $b\in \{ \mbox{ }_vx , \mbox{ }_vy : v=1,...,m \} $
for each $z\in U$, where $(\partial \phi /\partial \mbox{ }_vx).h=
(\partial \phi /\partial \mbox{ }_vz).h$ for each $h\in \bf C$,
$(\partial \phi /\partial \mbox{ }_vy).h=
(\partial \phi /\partial \mbox{ }_vz).hj$ for each $h\in \bf C$.
\par Consider a right superlinear operator $A$
on $X$, then $A(h) = \sum_{l=0}^3 A(i_l)h_l$, where
$h_l\in \bf R^m$, $i_0:=1$, $i_1:=i$, $i_2:=j$, $i_3:=k$,
$h=h_0 + i_1h_1 + i_2h_2 + i_3h_3$, $h\in X=\bf H^m$.
From \S 3.28 \cite{luoyst} it follows, that $A(h)=
A(1)(h - i_1hi_1 - i_2hi_2 - i_3hi_3)/4 +
A(i_1)(h - i_1hi_1 + i_2hi_2 + i_3hi_3)/4 $
$+ A(i_2)(h + i_1hi_1 - i_2hi_2 + i_3hi_3)/4 +
A(i_3)(h + i_1hi_1 + i_2hi_2 - i_3hi_3)/4 $ for each $h\in X$.
There are identities $i_1hi_1 = -h$, $i_2hi_2 = - \bar h$ and
$i_3hi_3 = - \bar h$ for each $h\in \bf C^m$, consequently,
$A(h)=(A(1)+A(i_1))h/2 + (A(1)-A(i_1)){\bar h}/2$ for each $h\in \bf C^m$
and in view of right superlinearity $A(1)=A(i_1)$, hence
\par $(i)$ $A(h)=A(1)h$ for each
$h\in {\bf C^m}:={\bf R^m}\oplus {\bf R^m}i_1$. \\
Then $A(hi_2)=(A(i_2)+A(i_3))hi_2/2 + (A(i_2) - A(i_3)){\bar h}i_2/2$
for each $h\in \bf C^m$ and in view of right superlinearity of $A$
we get $A(i_2)=A(i_3)$. Restriction on ${\bf C^m}i_3$ of operator $A$ shows,
that $A(1)=A(i_3)$, consequently,
\par $(ii)$ $A(hi_2)=A(1)hi_2$ for each $h\in \bf C^m$. \\
Since $(\partial a /\partial b).h$ is complex valued
for each $h\in \bf C$, for each
$a\in \{ \mbox{ }_v\alpha , \mbox{ }_v\beta : v=1,...,m \} $
and each $b\in \{ \mbox{ }_vx , \mbox{ }_vy : v=1,...,m \} $
for each $z\in U$, then Formulas $(i,ii)$ imply that
$(\partial a /\partial b).h$ is complex linear in $h$.
\par Consider now real local coordinates in $M$: $ \{ \mbox{ }_{v,l}\phi :
v= 1,...,m; l= 0,1,2,3 \} $ such that $\mbox{ }_v\phi =
\sum_{l=0}^3 \mbox{ }_{v,l} \phi i_l $, then the determinant
of change of local coordinates between charts is positive:
$det ( \{ \partial \mbox{ }_{v,l}\phi /
\partial \mbox{ }_{w,p}z \}_{v,l;w,p} ) = | det ( \{ 
\partial \mbox{ }_{v,l} \eta /\partial \mbox{ }_{w,p} \xi
\} _{v,l;w,p} )|^2 $, where $\mbox{ }_{v,l}\phi \in \bf R$
for each $v=1,...,m$ and $l=0,1,2,3$; $\mbox{ }_{v,0}\eta
+ \mbox{ }_{v,1} \eta i_2=\mbox{ }_v\phi $ and
$\mbox{ }_{v,0}\xi + \mbox{ }_{v,1} \xi i_2 =\mbox{ }_vz$
and $\mbox{ }_{v,l}\xi , \mbox{ }_{v,l}\eta \in \bf C$
for each $v=1,...,m$ and each $l=0, 1$.
\par {\bf 2.1.3.8. Remark.} Theorem $2.1.3.6$
is the quaternion analog of the Witney theorem. For the proof
of it was important that superdifferentials of quaternion
holomorphic functions may be nonlinear by $\bf H$.
Without condition of right or left superlinearity of superdifferentials
of transition mappings Theorem $2.1.3.7$ is untrue, since
there are square $4m\times 4m$-matrices $A$ with entries in $\bf R$ such that
$det (A)< 0$ to which correspond $\bf R$-homogeneous $\bf H$-additive
operators. Certainly Riemannian manifold need not be a quaternion
manifold even if its real dimension is $4m$ for $m\in \bf N$,
since transition mappings $\phi _i\circ \phi _l^{-1}$ need not be
quaternion holomorphic for charts with $U_i\cap U_j\ne \emptyset $.
\par {\bf 2.1.3.9. Theorem.} {\it Let $M$ be a
quaternion manifold, then
there exists an open neighborhood ${\tilde T}M$ of $M$ in $TM$
and an exponential quaternion holomorphic mapping
$\exp : {\tilde T}M\to M$ of ${\tilde T}M$ on $M$.}
\par {\bf Proof.} It was shown in \S \S 2.1.3.4 and 2.1.3.7
each quaternion manifold has also the structure of the Riemann
manifold. Therefore, the geodesic equation
\par  $(i)$  $\nabla _{\dot c}{\dot c}=0$ with initial conditions
$c(0)=x_0$, ${\dot c}(0)=y_0$, $x_0\in M$, $y_0\in T_{x_0}M$ \\
has a unique $C^{\infty }$-solution, $c: (- \epsilon , \epsilon )
\to M$ for some $\epsilon >0$. For a chart $(U_j, \phi _j)$
containing $x$, put $\psi _j(\beta ) = \phi _j\circ c(\beta )$,
where $\beta \in ( - \epsilon , \epsilon )$.
Consider an $\bf H$-Hermitian metric $g$ in $M$, then $g$ is quaternion
holomorphic in local quaternion coordinates in $M$ (see \S 2.1.3.4),
$\partial g(z)/\partial {\tilde z}=0$, where $g(z)(*,*)=<*,*>_z$
is the $\bf H$-Hermitian inner product in $T_zM$, $z\in M$,
where $\phi _l(z)\in X$ is denoted for convenience also by $z$.
Consider real-valued inner product induced by $g$ of the form
$G(z)(x,y) := (g(z)(x,y) - \sum_{l=1}^3 i_l g(z)(x,y) i_l)/4$.
Then $G(z)$ is also quaternion holomorphic:
$\partial G(z)(x,y)/\partial {\tilde z}=0$ for each $x, y$.
In real local coordinates $G$ can be written as:
$G(z)(\eta ,\xi )=\sum_{l,p} G^{l,p}(z)\eta _l \xi _p$,
where $x = (\mbox{ }_1x,...,\mbox{ }_mx)$,
$\mbox{ }_vx=\sum_{l=4v-3}^{4v} \eta _l i_{l-4v+3}$
and $\mbox{ }_vy=\sum_{l=4v-3}^{4v} \xi _l i_{l-4v+3}$
for each $v=1,...,m$. Then the Christoffel symbols are
given by the equation $\Gamma ^a_{b,c} =
(\sum_{l=1}^{4m} G^{a,l} (\partial _bG_{c,l} + \partial _c
G_{l,b} - \partial _l G_{b,c})/2$ for each $a, b, c=1,...,4m$.
Using expression of the Christoffel symbol $\Gamma $ through
$g$, we get that $\Gamma (z)$ is quaternion holomorphic:
\par  $(ii)$  $\partial  \Gamma ^a_{b,c}(z)/\partial {\tilde z}=0$
for each $a, b, c=1,...,4m$. \\
Thus differential operator
corresponding to the geodesic equation has quaternion holomorphic form:
\par $(iii)$ $d^2 c(t)/dt^2 + \Gamma (c(t)) (c'(t), c'(t))=0$.
Therefore, the mapping ${\tilde T}V_1\times (-\epsilon , \epsilon )
\ni (z_0,y_0;\beta )\mapsto \psi _j(\beta ;x_0,y_0)$
is quaternion holomorphic by $(x_0,y_0)$, since components of $y_0$
can be expressed through $\bf R$-linear combinations of
$\{ i_ly_0i_l: l=0, 1, 2, 3 \} $, where
$0<\epsilon $, $z_0=\phi _j(x_0)\in V_1\subset V_2\subset
\phi _j(U_j)$, $V_1$ and $V_2$ are open, $\epsilon $ and ${\tilde T}V_1$
are sufficiently small, that to satisfy the inclusion
$\psi _j(\beta ;x_0,y_0)\in V_2$ for each $(z_0,y_0;\beta )\in
{\tilde T}V_1\times (-\epsilon , \epsilon )$.
\par Then there exists $\delta >0$ such that
$c_{aS}(t)=c_S(at)$ for each $a\in (- \epsilon , \epsilon )$
with $|aS(\phi _j(q))|<\delta $
since $dc_S(at)/dt=a(dc_S(z)/dz)|_{z=at}$.
The projection $\tau :=\tau _M: TM\to M$ is given by
$\tau _M(s)=x$ for each vector $s\in T_xM$, $\tau _M$ is the tangent bundle.
For each $x_0\in M$ there exists a chart $(U_j,\phi _j)$
and open neighbourhoods $V_1$ and $V_2$,
$\phi _j(x_0)\in V_1\subset V_2\subset
\phi _j(U_j)$ and $\delta >0$ such that from $S\in TM$ with
$\tau _MS=q\in \phi _j^{-1}(V_1)$ and $|S(\phi _j(q))|<\delta $
it follows, that the geodesic $c_S$ with $c_S(0)=S$ is defined
for each $t\in (-\epsilon ,\epsilon )$ and $c_S(t)\in \phi _j^{-1}(V_2)$.
Due to paracompactness of $TM$ and $M$ this covering
can be chosen locally finite \cite{eng}.
\par This means that there exists an open neighbourhood
${\tilde T}M$ of $M$ in $TM$ such that a geodesic $c_S(t)$
is defined for each $S\in {\tilde T}M$ and each $t\in (-\epsilon ,
\epsilon )$.
Therefore, define the exponential mapping $\exp : {\tilde T}M\to M$
by $S\mapsto c_S(1)$, denote by $\exp _x:=\exp |_{{\tilde T}M\cap T_xM}$
a restriction to a fibre. Then $\exp $ has a local representation
$(x_0,y_0)\in {\tilde T}V_1\mapsto
\psi _j(1;x_0,y_0)\in V_2\subset \phi _j(U_j).$
From Equations $(i - iii)$ it follows that
$\exp $ is quaternion holomorphic from ${\tilde T}M$ onto $M$.
\par {\bf 2.1.3.10. A uniform space of piecewise holomorphic mappings.}
For manifolds $M$ and $N$ with corners both either complex or quaternion
let either ${\sf O}_{\Upsilon }(M,N)$ or ${\sf H}_{\Upsilon }(M,N)$
denotes a space of continuous mappings $f: M\to N$ such that
for each $f$ there exists a partition $Z_f$ of $M$ via
a real $C^{\infty }$-submanifold ${M'}_f$, 
which may be with corners, such that its codimension 
over $\bf R$ in $M$ is $codim_{\bf R}{M'}_f=1$
and $M\setminus {M'}_f$ is a disjoint union of open either
complex submanifolds
or quaternion submanifolds $M_{j,f}$ respectively possibly with corners 
with $j=1,2,...$ such that each restriction $f|_{M_{j,f}}$
is either complex holomorphic or quaternion holomorphic correspondingly
with all its derivatives bounded on $M_{j,f}$.
For a given partition $Z$ (instead of $Z_f$) and the corresponding
$M'$ the latter subspace of continuous piecewise either complex
or quaternion holomorphic mappings
$f: M\to N$ is denoted either by ${\sf O}(M,N;Z)$ or ${\sf H}(M,N;Z)$
respectively. The family $\{ Z \} $ of all such partitions is denoted
$\Upsilon $.
That is ${\sf O}_{\Upsilon } (M,N)=str-ind_{\Upsilon }{\sf O}(M,N;Z)$
and ${\sf H}_{\Upsilon } (M,N)=str-ind_{\Upsilon }{\sf H}(M,N;Z)$.
Let also ${\sf O}(M,N)$ and ${\sf H}(M,N)$ denote the uniform spaces
of all either complex holomorphic or quaternion holomorphic
respectively mappings $f: M\to N$, $Diff^{\infty }(M)$ denotes a group of 
$C^{\infty }$-diffeomorphisms of $M$ and
$Diff^{\sf O}_{\Upsilon }(M) := Hom(M)\cap {\sf O}_{\Upsilon }(M,M)$,
$Diff^{\sf H}_{\Upsilon }(M):=
Hom(M)\cap {\sf H}_{\Upsilon }(M,M)$, where $Hom (M)$ is a group
of all homeomorphisms of $M$.
\par Let $A$ and $B$ be two manifolds either both either real or complex
or quaternion with corners such that $B$ is a submanifold of $A$.
Then $B$ we call a strong $C^r([0,1]\times A,A)$-retract 
(or $C^r([0,1],{\sf O}_{\Upsilon }(A,A))$-retract
or $C^r([0,1],{\sf H}_{\Upsilon }(A,A))$-retract respectively)
of  $A$ if there exists a mapping $F: [0,1]\times A\to A$
such that $F(0,z)=z$ for each
$z\in A$ and $F(1,A)=B$ and $F(x,A)\supset B$ for each $x\in
[0,1]:=\{ y: 0\le y \le 1, y\in {\bf R} \} $,
$F(x,z)=z$ for each $z\in B$ and $x\in [0,1]$,
where $F\in C^r([0,1]\times A,A)$ or $F\in C^r([0,1],{\sf O}_{\Upsilon }
(A,A))$ or $F\in C^r([0,1],{\sf H}_{\Upsilon }(A,A))$
respectively, $r\in [0,\infty )$, $F=F(x,z)$, $x\in [0,1]$,
$z\in A$. Such $F$ is called the retraction. In the case of $B=\{ a_0 \} $,
$a_0\in A$ we say that $A$ is $C^r([0,1]\times A,A)$-contractible
(or $C^r([0,1], {\sf O}_{\Upsilon }(A,A))$-contractible
or $C^r([0,1], {\sf H}_{\Upsilon }(A,A))$-contractible
correspondingly). Two maps $f: A\to E$ and $h: A\to E$ we call 
$C^r([0,1]\times A,E)$-homotopic (or 
$C^r([0,1],{\sf O}_{\Upsilon }(A,E))$-homotopic
or $C^r([0,1],{\sf H}_{\Upsilon }(A,E))$-homotopic) if there
exists $F\in C^r([0,1]\times A,E)$ (or $F\in C^r([0,1],{\sf O}_{\Upsilon }
(A,E))$ or $F\in C^r([0,1],{\sf H}_{\Upsilon }(A,E))$
respectively) such that $F(0,z)=f(z)$ and $F(1,z)=h(z)$
for each $z\in A$, where $E$ is also either a real or complex 
or quaternion manifold. Such $F$ is called the homotopy.
\par Let $M$ be either a real or complex manifold or a quaternion manifold
with corners satisfying the following conditions:
\par $(i)$ it is compact;
\par $(ii)$ $M$ is a union of two closed submanifolds both either
real or complex or quaternion $A_1$ and $A_2$
with corners, which are canonical closed subsets in $M$
with $A_1\cap A_2=\partial A_1\cap \partial A_2=:A_3$ 
and a codimension over $\bf R$ of $A_3$ in $M$ is $codim_{\bf R}A_3=1$;
\par $(iii)$ a marked point $s_0$ is in $A_3$;
\par $(iv)$ $A_1$  and $A_2$ are  either
$C^0([0,1]\times A,A)$-contractible or
$C^0([0,1],{\sf O}_{\Upsilon }(A_j,A_j))$-contractible
or $C^0([0,1],{\sf H}_{\Upsilon }(A_j,A_j))$-contractible
respectively into a marked point
$s_0\in A_3$ by mappings $F_j(x,z),$ where either $j=1$ or $j=2$.
In the complex and quaternion cases more general condition of
$C^0([0,1],{\sf O}_{\Upsilon }(A_j,A_j))$-contractibility
or $C^0([0,1],{\sf H}_{\Upsilon }(A_j,A_j))$-contractibility
of $A_j$ on $X_0\cap A_j$ can be considered,
where $X_0$ is a closed subset in $M$, $j=1$ or $j=2$, $s_0\in X_0$,
$dim X_0< dim_{\bf R}M$, $dim X_0$ is a covering dimension
of $X_0$ (see its defintion in \cite{eng}).
\par  We consider all finite partitions $Z:=\{ M_k:  k\in \Xi _Z\}$ of $M$
such that $M_k$ are submanifolds either all complex or quaternion
respectively (of $M$), which may be with corners and
$\bigcup_{k=1}^sM_k=M$, $\Xi _Z=\{ 1,2,...,s \}$, $s\in \bf N$ 
depends on $Z$, $M_k$ is a canonical closed subset of $M$ for each $k$.
We denote by $\tilde diam(Z):=\sup_k(diam
(M_k))$ a diameter of a partition $Z$, where $diam (A)=\sup_{x, y\in A}
|x-y|_X$ is a diameter of a subset $A$ in a normed space $X$, since
each finite dimensional manifold $M$ can be embedded into $\bf C^n$
or in $\bf H^n$ with a corresponding $n\in \bf N$. We suppose
also that $M_i\cap M_j\subset M'$ and $\partial M_j\subset M'$
for each $i\ne
j$, where $M'$ is a closed $C^{\infty }$-submanifold (which may be
with corners) in $M$ with
the codimension $codim_{\bf R}(M')=1$ of $M'$ in $M$,
$M'=\bigcup_{j\in \Gamma _Z}{M'}_j$, ${M'}_j$ are
$C^{\infty }$-submanifolds of $M$, $\Gamma _Z$ is a finite
subset of $\bf N$. 
\par  We denote by $H^t(M,N;Z)$ a space of continuous functions
$f: M\to N$ such that $f|_{(M\setminus M')}\in H^t(M\setminus M',N)$ and
$f|_{[Int(M_i)\cup (M_i\cap {M'}_j)]}\in
H^t(Int(M_i)\cup (M_i\cap {M'}_j),N)$,
when $\partial M_i\cap {M'}_j\ne \emptyset $, $h^Z_{Z'}:  H^t(M,N;Z)\to
H^t(M,N;Z')$ is an embedding for each $Z\le Z'$ in $\Upsilon $.
\par An ordering $Z\le Z'$ means by our definition,
that each submanifold $M_i^{Z'}$
from a partition $Z'$ either belongs to the family
$(M_j:  j=1,...,k)=(M_j^Z: j=1,...,k)$ for $Z$ or 
there exists $j$ such that $M_i^{Z'}\subset M_j^Z$ and
$M_j^Z$ is a finite union of $M_l^{Z'}$ for which $M_l^{Z'}\subset
M_j^Z$. Moreover, $M_l^{Z'}$ is a submanifold (may be
with corners) in $M_j^Z$ for each $l$ and a corresponding $j$. 
\par Then we consider the following uniform space
$H^t_p(M,N)$ that is the strict inductive limit $str-ind \{ H^t(M,N;Z);
h^{Z'}_{Z}; \Upsilon \} $ (the index $p$ reminds about the procedure of
partitions), where $\Upsilon $ is the directed family of all such $Z$, 
for which $\lim_{\Upsilon }\tilde diam(Z)=0$.
\par {\bf 2.1.4. Notes and definitions of loop monoids.}
Let now $s_0$ be a marked
point in $M$ such that $s_0\in A_3$ (see \S 2.1.3.10) 
and $y_0$ be a marked point in a manifold $N$.
\par $(i).$ Suppose that $M$ and $N$ are connected.
\par Let
$H^t_p(M,s_0;N,y_0):=\{ f\in H^t(M,N)| f(s_0)=y_0 \} $ denotes the closed
subspace of $H^t(M,N)$ and $\omega _0$ be its element such that $\omega
_0(M)= \{ y_0 \} $, where $\infty \ge t\ge m+1$, $2m=dim_{\bf R}M$
such that $H^t\subset C^0$ due to the Sobolev embedding theorem.  
The following uniform subspace
$ \{ f: f\in H^{\infty }_p(M,s_0;N,y_0),\quad
{\bar \partial }f=0 \} $ is uniformly isomorphic with
${\sf O}_{\Upsilon }(M,s_0;N,y_0)$ for complex manifolds,
while the uniform subspace $ \{ f: f\in H^{\infty }_p(M,s_0;N,y_0),$ 
${\tilde \partial }f=0 \} $ is uniformly isomorphic with
${\sf H}_{\Upsilon }(M,s_0;N,y_0)$ for quaternion manifolds $M$ and $N$,
since $f|_{(M\setminus M')} \in H^{\infty }(M\setminus M',N)=
C^{\infty }(M\setminus M',N)$ and ${\bar \partial }f=0$
or ${\tilde \partial }f=0$ respectively.
\par Let as usually $A\vee B:=A\times \{
b_0\}\cup \{ a_0\}\times B\subset A\times B$ 
be the wedge sum of pointed spaces $(A,a_0)$ and $(B,b_0)$, where $A$
and $B$ are topological spaces with marked points $a_0\in A$ and $b_0\in B$.
Then the wedge combination $g\vee f$ 
of two elements $f, g\in H^t_p(M,s_0;N,y_0)$ is
defined on the domain $M\vee M$.
\par The uniform spaces
${\sf O}_{\Upsilon }(J,A_3;N,y_0):=\{ f\in {\sf O}_{\Upsilon }
(J,N): f(A_3)=\{ y_0 \} \} $ 
have the manifold structure and have embeddings into
${\sf O}_{\Upsilon }(M,s_0;N,y_0)$ for complex manifolds,
while the uniform spaces
${\sf H}_{\Upsilon }(J,A_3;N,y_0):=\{ f\in {\sf H}_{\Upsilon }
(J,N): f(A_3)=\{ y_0 \} \} $ have the manifold structure
and have embeddings into ${\sf H}_{\Upsilon }(M,s_0;N,y_0)$ due to Condition 
$2.1.3.10(ii)$ for quaternion manifolds $M$ and $N$,
where either $J=A_1$ or $J=A_2$.
This induces the following embedding
$\chi ^*: {\sf O}_{\Upsilon }(M\vee M,s_0\times s_0;N,y_0)\hookrightarrow
{\sf O}_{\Upsilon }(M,s_0;N,y_0)$ for complex manifolds, also
$\chi ^*: {\sf H}_{\Upsilon }(M\vee M,s_0\times s_0;N,y_0)\hookrightarrow
{\sf H}_{\Upsilon }(M,s_0;N,y_0)$ for quaternion manifolds $M$ and $N$.
Considering $H^t_p(M,X_0;N,y_0)=\{ f\in H^t(M,N): f(X_0)=\{ y_0 \} \} $
and ${\sf O}_{\Upsilon }(J,A_3\cup X_0;N,y_0)$ we get the embedding 
$\chi ^*: {\sf O}_{\Upsilon }(M\vee M,X_0\times X_0;N,y_0)\hookrightarrow
{\sf O}_{\Upsilon }(M,X_0;N,y_0)$ for complex manifolds, also
$\chi ^*: {\sf H}_{\Upsilon }(M\vee M,X_0\times X_0;N,y_0)\hookrightarrow
{\sf H}_{\Upsilon }(M,X_0;N,y_0)$ for quaternion manifolds $M$ and $N$.
Therefore $g\circ f:=\chi ^*(f\vee g)$ is the composition in
${\sf O}_{\Upsilon }(M,s_0;N,y_0)$, also in
${\sf H}_{\Upsilon }(M,s_0;N,y_0)$.
\par The space $C^{\infty }(M,N)$ is dense in $C^0(M,N)$
and there is the inclusion ${\sf O}_{\Upsilon }(M,N)\subset
H^{\infty }_p(M,N)$ for complex manifolds, also
${\sf H}_{\Upsilon }(M,N)\subset
H^{\infty }_p(M,N)$ for quaternion manifolds.
Let $M_{\bf R}$ be a Riemannian manifold generated by
a manifold $M$ considered over $\bf R$.
Then $Diff^{\infty }_{s_0}(M_{\bf R})$
is a group of $C^{\infty }$-diffeomorphisms $\eta $ of $M_{\bf R}$ 
preserving the marked point $s_0$, that is $\eta (s_0)=s_0$.
There exists the following equivalence relation
$R_{\sf O}$ in ${\sf O}_{\Upsilon }(M,s_0;N,y_0)$:
$fR_{\sf O}h$  (also $R_{\sf H}$ in ${\sf H}_{\Upsilon }(M,s_0;N,y_0)$:
$fR_{\sf H}h$) if and only if there exist nets
$\eta _n\in Diff^{\infty }_{s_0}(M_{\bf R})$, 
also $f_n$ and $h_n\in H^{\infty }_p(M,s_0;N,y_0)$ with $\lim_{n}f_n=f$ 
and $\lim_{n}h_n=h$ such that $f_n(x)=h_n(\eta _n(x))$ for each
$x\in M$ and $n\in \omega $, where $\omega $ is a directed set,
$f, h \in {\sf O}_{\Upsilon }(M,s_0;N,y_0)$
(or $f, h \in {\sf H}_{\Upsilon }(M,s_0;N,y_0)$  respectively)
and converegence is considered in $H^{\infty }_p(M,s_0;N,y_0)$. 
In general, considering $Diff^{\infty }_{X_0}(M_{\bf R}):=
\{ f\in Diff^{\infty }(M_{\bf R}): f(X_0)=X_0 \} $ and elements
$f, h$ in ${\sf O}_{\Upsilon }(M,X_0;N,y_0)$ and a convergence
in $H^{\infty }(M,X_0;N,y_0)$ we get the equivalence relation
$R_{\sf O}$ in ${\sf O}_{\Upsilon }(M,X_0;N,y_0)$ for complex manifolds,
also with $f, h$ in ${\sf H}_{\Upsilon }(M,X_0;N,y_0)$ and a convergence
in $H^{\infty }(M,X_0;N,y_0)$ this leads to the equivalence relation
$R_{\sf H}$ in ${\sf H}_{\Upsilon }(M,X_0;N,y_0)$.
\par The quotient uniform space ${\sf O}_{\Upsilon }(M,X_0;N,y_0)/R_{\sf O}=:
(S^MN)_{\sf O}$ for complex manifolds also
${\sf H}_{\Upsilon }(M,X_0;N,y_0)/R_{\sf H}=:
(S^MN)_{\sf H}$ for complex manifolds we call the loop monoid.
For the spaces $H^t_p(M,s_0;N,y_0)$ the corresponding equivalence 
relations are denoted $R_{t,H}$, for them the loop monoids are denoted 
by $(S^M_{\bf R}N)_{t,H}$. When real  manifolds are considered 
we omit the index $\bf R$.
\par {\bf 2.1.4.1. Theorem.} {\it  The uniform spaces
$(S^MN)_{t,H}$ for real manifolds, $(S^MN)_{\sf O}$
for complex manifolds, $(S^MN)_{\sf H}$ for quaternion manifolds
$M$ and $N$ have sructures of complete topological Abelian monoids with a
unit and with a cancellation property, where $t>m+5$,
$m := dim_{\bf R}M$. It is nondiscrete for $dim_{\bf R}N>1$.
From these monoids it is possible to construct
uniform spaces $(L^mN)_{t,H}$, $(L^MN)_{\sf O}$ and $(L^MN)_{\sf H}$
respectively which have structures of
the complete separable Abelian topological groups, where $t>m+5$.}
\par {\bf Proof.} Consider the case of quaternion manifolds $M$ and $N$.
In view of the results from \cite{omo}
$Hom(M)\cap H^t(M,M)$ is the topological group of diffeomorphisms
for $t>m+5$, where $Hom(M)$ is the group of homeomorphisms of
$M$. Then $C^t(M,{\bf R})\subset H^{t+[m/2]+1}(M,{\bf R})$
for each $m\in \bf N$ and $t\ge 0$ due to the Sobolev embedding 
theorem \cite{miha}. 
Hence for each $f\in H^t_p(M,s_0;N,y_0)$ the range $f(M)$ is
compact and connected in $N$. In view of \cite{seel,touger,tri}
each $f\in H^t(M_k,N)$ has an extension of the same class
of smoothness onto an open neighbourhood $U$ of $M_k$ in $M$.
Therefore, in view of Lemmas 6.8 and 6.9 \cite{swit}
for each partition $Z$ there exists $\delta >0$
such that for each partition $Z"$ with $\sup_i \inf_j dist(M_i,{M"}_j)
<\delta $ and $f\in  H^t(M,N;Z)$, $f(s_0)=y_0$, there exists 
$f_1\in H^t(M,N;Z")$, $f_1(s_0)=y_0$, such that $fR_{t,H}f_1$,
hence we can choose $f_1$ with $fR_{\sf H}f_1$ for
${\tilde \partial }f=0$ and
${\tilde \partial }f_1=0$ on corresponding quaternion submanifolds
of $M$ prescribed by the partitions $Z$ and $Z"$,
where $dist(A,B)=\max (\sup_{a\in A}D(a,B); \sup_{b\in B}D(b,A)),$
$D(a,B):=\inf_{b\in B}d(a,b)$, $A$ and $B$ are subspaces of  
the metric space $M$ with the metric $d(a,b)=|a-b|_{\bf H^{2m+1}}$
(see Theorem 2.1.3.6).
\par Hence there exists a countable subfamily
$\{ Z_j: j \in {\bf N} \} $ in $\Upsilon $ such that $Z_j\subset Z_{j+1}$
for each $j$ and $\lim_j \tilde diam Z_j=0$. Then 
\par $(i)$ $str-ind \{ {\sf H}
(M,s_0;N,y_0;Z_j); h^{Z_i}_{Z_j}; {\bf N} \} /R_{\sf H}=(S^MN)_{\sf H}$
is separable, since
each space ${\sf H}(M,s_0;N,y_0;Z_j)$ is separable. The continuity of the 
composition and the inversion follows from notes in \S 1.4.
The space $str-ind \{ {\sf H}
(M,s_0;N,y_0;Z_j); h^{Z_i}_{Z_j}; {\bf N} \}$ is complete due to 
theorem 12.1.4 \cite{nari}, each class of $R_{\sf H}$-equivalent elements 
is closed in it, where ${\sf H}(M,s_0;N,y_0;Z):=
{\sf H}_{\Upsilon }(M,s_0;N,y_0) \cap H^t(M,N;Z)$.
Then to each Cauchy sequence in $(S^MN)_{\sf H}$
there corresponds a Cauchy sequence in 
$str-ind \{ {\sf H}(M\times [0,1],s_0\times 0;N,y_0;Z_j\times Y_j); h^{Z_i
\times Y_i}_{Z_j\times Y_j}; {\bf N} \}$ due to theorems about 
extensions of functions, where $Y_j$ are partitions of $[0,1]$
with $\lim_j\tilde diam(Y_j)=0$, $Z_j\times Y_j$ are the corresponding 
partitions of $M\times [0,1]$. Hence $(S^MN)_{\sf H}$ is complete.
\par  In view of Lemma 2.27 \cite{swit} it can be shown that there
exists  a $C^{\infty }([0,1],{\sf H})$-homotopy $h: (M_1 \vee M_2)
\times [0,1] \to M_2\vee M_1$, where $M=M_j$ for $j=1$ or $j=2$. Therefore,
$(g\circ f)(h(s,z))=(g\circ f)(s)$ for $z=0$ and $(g\circ 
f)(h(s,z))=f\circ g(s)$ for $z=1$ for each $s\in M$, consequently,
$(f\circ g)R_{\sf H}(g\circ f)$ for each $f, g\in {\sf H}_{\Upsilon }
(M,s_0;N,y_0)$, since $M_2\vee M_1\setminus \{ s_{0,1}\times s_{0,2} \} $
and $M_1\vee M_2\setminus \{ s_{0,2}\times s_{0,1} \}$
are $C^{\infty }$-diffeomorphic. Hence $(S^MN)_{\sf H}$ is Abelian.
\par For a commutative monoid $(S^MN)_{\sf H}$ with a unit and
with the cancellation property
there exists a commutative group $(L^MN)_{\sf H}$.
Algebraically this group is the quotient group
$F/\sf B$, where $F$ is a free Abelian group generated by 
$(S^MN)_{\sf H}$ and $\sf B$ is a subgroup of $F$ generated by
elements $[f+g]-[f]-[g]$, $f$ and $g\in (S^MN)_{\sf H}$,
$[f]$ denotes an element of $F$ corresponding to $f$. The natural mapping 
$\gamma : (S^MN)_{\sf H}\to (L^MN)_{\sf H}$ is injective.
We supply $F$ with a topology inherited from
the Tychonoff product topology of $(S^MN)_{\sf H}^{\bf Z}$,
where each element $z$ of $F$ is $z=\sum_fn_{f,z}[f]$, 
$n_{f,z}\in \bf Z$ for each $f\in (S^MN)_{\sf H}$,
$\sum_f|n_{f,z}|<\infty $. In particular $[nf]-n[f]\in \sf B$,
hence $(L^MN)_{\sf H}$ is the complete topological group and
$\gamma $ is the topological embedding
such that $\gamma (f+g)=\gamma (f)+ \gamma (g)$
for each $f, g \in (S^MN)_{\sf H}$, $\gamma (e)=e$,
since $(z+B)\in \gamma (S^MN)_{\sf H}$, when $n_{f,z}\ge 0$
for each $f$, so in general $z=z^{+}-z^{-}$, where
$(z^{+}+B)$ and $(z^{-}+B)\in \gamma (S^MN)_{\sf H}$.
Cases of complex and real manifolds are proved analogously.
\par {\bf 2.1.5.1. Notes and Definitions.} In view of \S I.5 \cite{kob}
a complex manifold $M$ or a quaternion manifold (see also \S 2.1.3.4) 
considered over $\bf R$
admits a (positive definite) Riemannian metric $g$, since $M$ is paracompact
(see \S \S 1.3 and 1.5 \cite{kob}). Due to Theorem IV.2.2
\cite{kob} there exists the Levi-Civit\`a connection (with vanishing torsion)
of $M_{\bf R}$. Each complex manifold and each quaternion manifold
with right (or left) superlinearly superdifferentiable transition mappings
of charts is orientable (see Theorem 2.1.3.7).
For the orientable manifold $M$ we suppose that $\nu $ is a measure on
$M$ corresponding to the Riemannian volume element $w$ ($m$-form)
$\nu (dx)=w(dx)/w(M)$. The Riemannian volume element $w$ is non-degenerate
and non-negative, since $M$ is orientable.
For a nonorientable $M$ consider its double covering orientable manifold
$\tilde M$ and the quotient mapping $\theta _M: \tilde M\to M$,
then the Riemannian volume element $w$ on $\tilde M$
produces the following measure $\nu (S):=w(\theta _M^{-1}(S))/
(2w(\tilde M))$ for each Borel subset $S$ in $M$.
\par The Christoffel
symbols $\Gamma ^k_{i,j}$ of the Levi-Civit\`a derivation
(see \S 1.8.12 \cite{kling}) are of class $C^{\infty }$ for $M$.
Then the
equivalent uniformity in $H^t(M,N)$ for $0\le t<\infty $
is given by the following base $\{
(f,g)\in (H^t(M,N))^2:  {\| (\psi _j\circ f -\psi _j\circ g)\|"}_{
H^t(M,X)}<\epsilon $, 
where $D^{\alpha }=\partial ^{|\alpha |}/\partial (x
^1)^{\alpha ^1} ...\partial (x ^{km})^{\alpha ^{km}}$, $\epsilon >0$, 
${{\| (\psi
_j\circ f -\psi _j\circ g)\|"}^2}_{ H^t(M,X)}:= \sum_{|\alpha |\le t}
\int_{M}|D^{\alpha } (\psi _j\circ f(x ) 
-\psi _j\circ g(x ))|^2\nu
(dx) \} $, $j\in \Lambda _N$, $X$ is the Hilbert space  over $\bf C$
either $\bf C^n$ or $l_2({\bf C})$, 
or over $\bf H$ either $\bf H^n$ or $l_2({\bf H})$ respectively,
$x$ are local normal coordinates in $M_{\bf R}$, $k=2$ or $k=4$ respectively.
We consider submanifolds $M_{i,k}$ and ${M'}_{j,k}$
for each partition $Z_k$ as in \S 2.1.3.10
(with $Z_k$ instead of $Z$), $i\in \Xi _{Z_k}$, $j\in \Gamma _{Z_k}$,
where $\Xi _{Z_k}$ and $\Gamma _{Z_k}$ are finite subsets of $\bf N$. 
We supply $H^{\gamma }(M,X;Z_k)$ with the following metric 
\par $\rho _{k, \gamma }(y) := [\sum_{i\in
\Xi }$ $ {{\| y|_{M_{i,k}}\|"}^2}_{\gamma ,i,k}]^{1/2}$
for $y\in H^{\gamma }(M,X;Z_k)$ \\
and $\rho _{k, \gamma }(y)=+\infty $ in the
contrary case, where $\Xi =\Xi _{Z_k}$, $\infty >t\ge \gamma \in \bf N$,
$\gamma \ge m+1$, ${\|
y|_{M_{i,k}}\|"}_{\gamma ,i,k}$ is given analogously to ${\| y\|"}_{H^{\gamma
}(M,X)}$, but with $\int_{M_{i,k}}$ instead of $\int_{M}$, where
$m:=dim_{\bf R}M$.
\par Let $Z^{\gamma }(M,X)$ be the completion of
$str-ind \{ H^{\gamma }(M,X;Z_j); h^{Z_i}_{Z_j};
{\bf N} \}=:Q$ relative to the following norm 
${\| y\| '}_{\gamma }:= \inf_k\rho _{k, \gamma }(y)$, as usually
$Z^{\infty }(M,X)=\bigcap_{\gamma \in \bf N}Z^{\gamma }(M,X)$. Let
\par ${\bar Y}^{\infty }(M,X):= \{ f: f\in Z^{\infty }(M,X),$
${\bar \partial }f_j|_{M_{j,k}}=0$ $\mbox{for each }k $ \\
$\mbox{for}$ $M$ $\mbox{and}$ $X$ $\mbox{over}$ $\bf C$,
$\mbox{or}$ ${\tilde \partial }f_j|_{M_{j,k}}=0$ $\mbox{for}$ $M$
$\mbox{and}$ $N$ $\mbox{over}$ ${\bf H} \} $, \\
where $f\in Z^{\infty }(M,X)$ imples $f=\sum_jf_j$ with
$f_j\in H^{\infty }(M,X;Z_j)$ for each $j\in \bf N$.
\par For a domain $W$ in $\bf C^n$ or in $\bf H^n$,
which is a complex or quaternion respectively manifold with corners,
let $Y^{\Upsilon ,a}(W,X)$ (and $Z^{\Upsilon ,a}(W,X)$)
be a subspace of those
$f\in {\bar Y}^{\infty }(W,X)$ (or $f\in Z^{\infty }(W,X)$ respectively)
for which 
$$\| f\|_{\Upsilon ,a}:=(\sum_{j=0}^{\infty }
({\| f\|^{*}}_j)^2/[(j!)^{a_1}j^{a_2}])^{1/2}<\infty ,$$ 
where ${({\| f\|^{*}}_j)^2}:=
({\| f\| '}_j)^2-
({\| f\| '}_{j-1})^2$ for $j\ge 1$ and ${{\| f\|^*}_0}={\| f\|'}_0$, 
$a=(a_1,a_2)$, $a_1$ and $a_2\in \bf R$, $a<a'$ if either
$a_1<{a'}_1$ or $a_1={a'}_1$ and $a_2<{a'}_2$.
\par Using the atlases $At(M)$ and $At(N)$ 
for $M$ and $N$ of class of smoothness $Y^{\Upsilon ,b}\cap C^{\infty }$ 
with $a\ge b$
we get the uniform space $Y^{\Upsilon ,a}(M,X_0;N,y_0)$ (and also
$Z^{\Upsilon ,a}(M,X_0;N,y_0)$)
of mappings $f:  M\to N$ with $f(X_0)=y_0$ such
that $\psi _j\circ f\in Y^{\Upsilon ,a}(M,X)$ 
(or $\psi _j\circ f\in Z^{\Upsilon ,a}(M,X)$ respectively) for each $j$,
where $\sum_{p\in \Lambda _M, j\in \Lambda _N}\| f_{p,j}-(w_0)_{p,j}
\|^2_{Y^{\Upsilon ,a}(W_{p,j},X)}<\infty $ for each 
$f\in Y^{\Upsilon ,a}(M,X_0;N,y_0)$ is satisfied with
$w_0(M)=\{ y_0 \} $, since $M$ is compact.
Substituting $w_0$ on a fixed mapping
$\theta : M\to N$ we get the uniform space
$Y^{\Upsilon ,a,\theta }(M,N)$. To each equivalence class
either $\{ g:  gR_{\sf O}f \}=:<f>_{\sf O}$
or $\{ g:  gR_{\sf H}f \}=:<f>_{\sf H}$
there corresponds an equivalence class either
$<f>_{\Upsilon ,a}:= cl(<f>_{\sf O}\cap Y^{\Upsilon ,a}(M,X_0;N,y_0))$
in the complex case
or $<f>_{\Upsilon ,a}:= cl(<f>_{\sf H}\cap Y^{\Upsilon ,a}(M,X_0;N,y_0))$
in the quaternion case
(or $<f>^{\bf R}_{\Upsilon ,a}:= cl(<f>_{\infty ,H}\cap Z^{\Upsilon ,a}
(M,X_0;N,y_0))$ in the real case), where the closure is taken
in $Y^{\Upsilon ,a}(M,X_0;N,y_0)$ 
(or $Z^{\Upsilon ,a}(M,X_0;N,y_0)$ respectively). 
This generates equivalence relations either $R_{\Upsilon ,a}$ in the complex
or quaternion cases or $R^{\bf R}_{\Upsilon ,a}$ in the real case
respectively. We use here the same notation in the complex and quaternion
cases, because it can not cause a confusion, when manifolds are described
as either complex or quaternion correspondingly.
We denote the quotient uniform spaces
$Y^{\Upsilon ,a}(M,X_0;N,y_0)/R_{\Upsilon ,a}$ 
and $Z^{\Upsilon ,a}(M,X_0;N,y_0)/R^{\bf R}_{\Upsilon ,a}$ 
by $(S^MN)_{\Upsilon ,a}$ and $(S^M_{\bf R}N)_{\Upsilon ,a}$
correspondingly.
\par {\bf 2.1.5.2. Gevrey-Sobolev classes of smoothness
of loops. Notes and definitions.} 
Let $M$ be an infinite dimensional complex or quaternion
$Y^{\xi '}$-manifold with corners modelled on $l_2({\bf C})$
or on $l_2({\bf H})$ such that
\par $(i)$ there is the sequence of the canonically embedded
complex or quaternion respectively submanifolds $\eta _m^{m+1}:
M_m\hookrightarrow M_{m+1}$ for each $m\in \bf N$
and to $s_{0,m}$ in $M_m$ it corresponds $s_{0,m+1}=
\eta _m^{m+1}(s_{0,m})$ in $M_{m+1}$, $dim_{\bf C}M_m=n(m)$
or $dim_{\bf H} M_m = n(m)\in \bf N$,
$0<n(m)<n(m+1)$ for each $m\in \bf N$, $\bigcup_m M_m$ is dense in $M$;
\par $(ii)$ $M$ and $At(M)$ are foliated, that is, 
\par $(\alpha )$ $u_i\circ u_j^{-1}|_{u_j(U_i\cap U_j)}\to l_2$
are of the form: $u_i\circ u_j^{-1}((z^l:$ $l\in {\bf N}))=
(\alpha _{i,j,m}(z^1,...,z^{n(m)}), \gamma _{i,j,m}(z^l:$ $l>n(m)))$
for each $m$, when $M$ is without a boundary. 
If $M$ is with a boundary, $\partial M\ne \emptyset $, then 
\par $(\beta )$ for each boundary component $M_0$ of $M$ 
and $U_j\cap M_0\ne \emptyset $ we have $\phi _j: U_j\cap M_0\to
H_{l,Q}$, moreover, $\partial M_m\subset \partial M$ for each $m$,
where $H_{l,Q} := \{ z\in Q_j:$ $x^{2l-1}\ge 0 \} $ in the complex case
or $H_{l,Q}:= \{ z\in Q_j:$ $x^{4l-3}\ge 0 \} $ in the quaternion case,
$Q_j$ is a quadrant in $l_2$ such that $Int_{l_2}H_{l,Q}\ne \emptyset $
(the interior of $H_{l,Q}$ in $l_2$), $z^l=x^{2l-1}+x^{2l}i$ in
the complex case or $z^l := \sum_{s=1}^4 x^{4l+s-4} i_{s-1}$ in the
quaternion case, $x^j\in \bf R$, $z^l\in \bf C$ or $z^l\in \bf H$
respectively (see also \S 2.1.2.4);
\par $(iii)$ $M$ is embedded into $l_2$ as a bounded closed subset;
\par $(iv)$ each $M_m$ satisfies conditions $2.1.3.10 (i-iv)$
with $X_{0,m}:=X_0\cap M_m$, where $X_0$ is a closed subset in $M$.
\par Let $W$ be a bounded canonical closed subset in $l_2({\bf K})$,
${\bf K}=\bf C$ or ${\bf K}=\bf H$ with a continuous piecewise
$C^{\infty }$-boundary and $H_m$ be an increasing sequence of
finite dimensional subspaces over $\bf K$,
$H_m\subset H_{m+1}$ and $dim_{\bf K} H_m = n(m)$ for
each $m\in \bf N$. Then there are spaces 
$P^{\infty }_{\Upsilon ,a}(W,X):=str-ind_mY^{\Upsilon ,a}(W_m,X)$,
where $W_m=W\cap H_m$ and $X$ is a separable Hilbert space over $\bf K$.
\par Let $Y^{\xi }(W,X)$
be the completion of $P^{\infty }_{\Upsilon ,a}(W,X)$
relative to the following norm 
$$\| f\|_{\xi }:=
[\sum_{m=1}^{\infty }{{\| f|_{W_m} \| " }^2}_{Y^{\Upsilon ,a}(W_m,X)}/[
(n(m)!)^{1+c_1}n(m)^{c_2}]]^{1/2},$$
where ${{\| f|_{W_m} \| " }^2}_{Y^{\Upsilon ,a}(W_m,X)}:=
{\| f|_{W_m} \| ^2}_{Y^{\Upsilon ,a}(W_m,X)}
-{\| f|_{W_{m-1}} \| ^2}_{Y^{\Upsilon ,a}(W_{m-1},X)}$
for each $m>1$ and
${\| f|_{W_1} \| "}_{Y^{\Upsilon ,a}(W_1,X)}:=
\| f|_{W_1} \| _{Y^{\Upsilon ,a}(W_1,X)}$;
$c=(c_1,c_2)$, $c_1$ and $c_2\in \bf R$, $c<c'$ if either 
$c_1<{c'}_1$ or $c_1={c'}_1$ and $c_2<{c'}_2$; $\xi =(\Upsilon ,a,c)$.
Let $M$ and $N$ be the $Y^{\Upsilon ,a',c'}$-manifolds with 
$a'<a$ and $c'<c$.
\par If $N$ is a finite dimensional complex or quaternion
$Y^{\Upsilon ,a'}$-manifold, then it is also a
$Y^{\Upsilon ,a',c'}$-manifold. There exists a strict inductive limit
$str-ind_mQ^m=:Q^{\infty }_{\Upsilon ,a}(N,y_0)$,
\par $str-ind_mY^{\Upsilon ,a}(M_m;N)=:Q^{\infty }_{\Upsilon ,a}(N)$,
where $Q^m:=Y^{\Upsilon ,a}(M_m,X_{0,m};N,y_0)$.
Then with the help of charts of $At(M)$ and $At(N)$ the space
$Y^{\xi }(W,X)$
induces the uniformity $\tau $ in $Q^{\infty }_{\Upsilon ,a}(N,y_0)$
and the completion of it relative to $\tau $ we denote by
$Y^{\xi }(M,X_0;N,y_0)$, where $\xi =(\Upsilon ,a,c)$ and
$\sum_{p\in \Lambda _M, j\in \Lambda _N}\| f_{p,j}-(w_0)_{p,j}
\|^2_{Y^{\xi }(W_{p,j},X)}<\infty $ for each 
$f\in Y^{\xi }(M,X_0;N,y_0)$ is supposed to be satisfied with
$w_0(M)=\{ y_0 \} $, since each $M_m$ is compact.
Substituting $w_0$ on a fixed mapping $\theta : M\to N$ we get a 
uniform space $Y^{\xi ,\theta }(M,N)$. Therefore,
using classes of equivalent elements from $Q^{\infty }_{\Upsilon ,a}(N,y_0)$ 
and their closures in $Y^{\xi }(M,X_0;N,y_0)$ as in \S 2.1.5.1 we get
the corresponding loop monoids which are denoted $(S^MN)_{\xi }$.
Substituting spaces $Y^{\Upsilon ,a}$ over $\bf C$
onto $Z^{\Upsilon ,a}$ over $\bf R$ with the respective modifications we get
spaces $Z^{\Upsilon ,a,c}(M,N)$ over $\bf R$ and
loop monoids $(S^M_{\bf R}N)_{\xi }$ 
for the multi-index $\xi =(\Upsilon ,a,c)$.
\par A relation between manifolds with corners and usual 
manifolds is given by the following lemma.
\par {\bf 2.1.6. Lemma.} {\it If $M$ is a real or complex or quaternion
manifold modelled on $X=\bf K^m$ or $X=l_2({\bf K})$, where
$\bf K$ is either $\bf R$ or $\bf C$ or $\bf H$
with an atlas $At (M) = \{ (V_j,\phi _j): j \}$,
then there exists an atlas $At'(M)= \{ (U_k,u_k,Q_k): k \} $
which refines $At(M)$, where $(V_j,\phi _j)$ is an usual chart
for each $j$ with a diffeomorphism $\phi _j: V_j\to \phi _j(V_j)$ such that 
$\phi _j(V_j)$ is a $C^{\infty }$-domain in $\bf K^n$ or $l_2({\bf K})$, 
$(U_k,u_k,Q_k)$ is a chart corresponding to quadrants $Q_k$ in $\bf K^n$
or $l_2({\bf K})$.}
\par {\bf Proof.} The covering $\{ V_j: j \} $ of $M$ has a refinement
$\{ W_l: l \} $ such that for each $j$ there exists $l=l(j)$
with $W_l\subset V_j$ so that each $\phi _j(W_l)$ is a simply connected 
region in $\bf K^m$ or $l_2({\bf K})$ which is not the whole space.
We choose $W_l$ such that 
\par $(i)$ either $W_l\cap \partial M=\emptyset $
or $W_l\cap \partial M$ is open in $\partial M$;
\par $(ii)$ $ \{ \pi _k(z)=z^k:$ $z\in \phi _j(W_l) \} =:E_{l,k}$,
$z\in X$, $\pi _k:$ $X\to \bf K$ are canonical projections
associated with the standard orthonormal base $\{ e_j:$ $j \} $ in  $X$,
$E_{l,k}$ are $C^{\infty }$-regions in $\bf K$, $\phi _j(W_l)=
\prod_{k=1}^mE_{l,k}$ for $X=\bf K^m$, or $\pi _J(\phi _j(W_l))=
\prod_{k\in J}E_{l,k}$ for each finite subset $J$ in $\bf N$
and the corresponding projection $\pi _J: X\to sp_{\bf K}
\{ e_k:$ $k\in J \} $. In the real case $E_{l,k}$ are open intervals
in $\bf R$. In the complex case in view of the Riemann Mapping Theorem 
for each $E_{l,k}$ there exists holomorphic 
diffeomorphism $q_{l,k}$ either onto $B^{-}_r:=\{ z\in {\bf C} :$ $|z|<r \} $
or $F_r:=\{ z\in {\bf C}:$ $|z|<r, x^1\ge 0 \} $, where
$z=x^1+x^2i$, $x^1, x^2\in \bf R$, $z\in \bf C$
(see \S 2.12 in \cite{heins}). The latter case appears
while treatment of $\pi _k(\phi _j(W_l\cap \partial M))\ne \emptyset $
(see \S 10.5.2 \cite{sanger} and \S 12 \cite{heins}).
\par In the quaternion case consider $\bf H$ as the linear
space ${\bf C}\oplus j {\bf C}$ over $\bf C$, where $j=i_2$.
Then right superlinear superdifferentiable mapping $q_{l,k}(z)$ of
a quaternion variable $z=a+jb$, $a, b\in \bf C$, $z\in \bf H$,
is characterized by the conditions $\partial q_{l,k}(a,b)/\partial
{\bar a}=0$ and $\partial q_{l,k}(a,b)/\partial {\bar b}=0$
while $q_{l,k}$ is written in $a, b$ variables (see Formulas
$2.1.3.7 (i,ii)$). Moreover,
$q_{l,k}$ is right linearly differentiable by $a$ and $b$,
that is, $s_{l,k}$ and $p_{l,k}$ are holomorphic in variables $a, b$,
where $q_{l,k}(z) = s_{l,k} + j p_{l,k} $. Then each substitution
$q_{l,k}\mapsto i_vq_{l,k}$ for each $i_v \in \{ 1, i_1, i_2, i_3 \} $
preserves the property of right superlinear superdifferentiability
and gives others dependent quaternion Cauchy-Riemann conditions
in such notation (see also Propositions 2.2 and 3.13 \cite{luoyst}).
Applying the Riemann mapping
theorem to this situation (over $\bf C^2$) we get, that
there exists a right superlinearly superdifferentiable mapping
$q_{l,k}(z)$ such that $E_{l,k}$ is quaternion holomorphically
diffeomorphic with $B^{-}_r := \{ z\in {\bf H} :$ $|z|<r \} $
or $F_r := \{ z\in {\bf H}:$ $|z|<r, x^1\ge 0 \} $, where
$z = x^1 + x^2i_1 + x^3i_2 + x^4i_3$, $x^1,..., x^4\in \bf R$, $z\in \bf H$.
The latter statement can be proved also analogously to the complex case
with the help of Theorem 3.15 and Remark 3.16 \cite{luoyst}
in the class of quaternion right superlinearly superdifferentaible
functions.
Since this is true in the class of right superlinearly superdifferentiable
functions, then this is true in the more abundant class
of superdifferentiable functions (quaternion holomorphic).
Slightly shrinking
covering if necessary we can choose $\{ W_l : l \} $ such
that  each $q_{l,k}$ and its derivatives are bounded on $E_{l,k}$.
In view of Central Theorem from \S 6.3 \cite{heins}
$q_{l,k}$ are boundary preserving maps. In view of Chapter 13 \cite{moret}
$B^{-}_r$ and $F_r$ have finite atlases with charts corresponding to
quadrants.
\par {\bf 2.1.6.1.} {\bf Note.} Vice versa there are complex
and quaternion manifolds with corners, which are not usual complex
or quaternion manifolds correspondingly, for example,
canonical closed domains $F$ in $\bf C^m$ or in $\bf H^m$ with piecewise 
$C^{\infty }$-boundaries, which are not of class $C^1$.
Since each complex or quaternion manifold $G$ has a boundary
$\partial G$ of class $C^{\infty }$ by Definition 2.1.2.6.
\par {\bf 2.1.6.2. Lemma.} {\it The uniform spaces
${\sf O}_{\Upsilon ,a,c}(M,N)$ and ${\sf H}_{\Upsilon ,a,c}(M,N)$
from \S \S 2.1.3 and 2.1.5
are infinite dimensional complex and quaternion manifolds respectively
dense in $C^0(M,N)$. Moreover, there exist their tangent bundles
$T{\sf O}_{\Upsilon ,a,c}(M,N)={\sf O}_{\Upsilon ,a,c}(M,TN)$
and $T{\sf H}_{\Upsilon ,a,c}(M,N)={\sf H}_{\Upsilon ,a,c}(M,TN).$
If $N=\bf K^n$ or $N=l_2({\bf K})$, where ${\bf K}=\bf C$ or $\bf H$,
then ${\sf O}_{\Upsilon ,a,c}(M,N)$ and ${\sf H}_{\Upsilon ,a,c}(M,N)$
are infinite dimensional topological vector spaces over $\bf C$
and $\bf H$ respectively.}
\par {\bf Proof.} The connecting mappings $\phi _j\circ \phi _k^{-1}$
of charts $(U_j,\phi _j)$ and $(U_k,\phi _k)$ with $U_j\cap U_k\ne 
\emptyset $ are complex or quaternion holomorphic on the corresponding
domains $\phi _k(U_j\cap U_k)$ for ${\bf K}=\bf C$ or $\bf H$ respectively. 
Consider the quaternion case. For each submanifold $M_j$ in $M$
we have $T{\sf H}(M_j,N)={\sf H}(M_j,TN)$ (see \cite{eliass}
and Proposition 2.1.3.3).
For each $f\in {\sf H}_{\Upsilon }(M,N)$ there exists a partition
$Z_f$ of $M$ such that $f|_{M_{j,f}}\in {\sf H}(M_{j,f},N)$ for each 
submanifold $M_{j,f}$ with corners  defined by $Z_f$.
In accordance with Proposition 2.1.3.3 and \S 2.1.4
$\phi _j^{-1}\circ \phi _k$ induce connecting mappings
$(\phi _j^{-1}\circ \phi _k)^*$
of the corresponding charts in ${\sf H}_{\Upsilon }(M,N)$ such that 
$(\phi _j^{-1}\circ \phi _k)^*(f(z)):=f\circ (\phi _j^{-1}\circ \phi _k)(z)$
for each $z\in U_j\cap U_k$ such that its Frech\'et derivatives are
the following $[\partial (\phi _j^{-1}\circ \phi _k)^*(f)/\partial f].h=
(\phi _j^{-1}\circ \phi _k)^*(h)$ and 
$[\partial (\phi _j^{-1}\circ \phi _k)^*(f)/\partial \bar f].h=0$,
where $h$ are vectors in $T_f{\sf H}_{\Upsilon }(U_j\cap U_k,N)$.
Thus transition mappings are quaternion holomorphic. Therefore,
$T{\sf H}_{\Upsilon }(M,N)={\sf H}_{\Upsilon }(M,TN)$, since
${\sf H}_{\Upsilon }(M,N)$ is the quaternion manifold (certainly of class 
$C^{\infty }$). 
\par In particular ${\sf H}_{\Upsilon }(M,Y)$ is a topological vector space
over $\bf H$ for $Y=\bf H^n$ or $Y=l_2({\bf H})$.
It remains to prove that the manifold ${\sf H}_{\Upsilon }(M,N)$ is
infinite dimensional and dense in $C^0(M,N)$. This follows from 
Corollary 3.2.3, Exer. 1.28 and Conjecture in Exer. 3.2 \cite{henle} and
\cite{lufsqv}, since for each quadrant $Q$ and a given function
$s$ on $\partial Q$, which is a restriction $q|_{\partial Q}$ of a
holomorphic function $q$ on a neighbourhood of $\partial Q$ in $\bf H^m$,
$m=dim_{\bf H}M$, there exists a space of functions $u: W\to \bf H^n$
such that $u|_{Int(Q)}$ and $u|_{W\setminus Q}$ 
are holomorphic and bounded together with each
partial derivative, where $W$ is
an open ball in $\bf H^n$ containing $Q$ and such that
$u^+(z)-u^-(z)=s(z)$ for each $z\in \partial Q$.
Then using Cauchy integration along $C^{\infty }$-curves
we construct a space of continuous functions 
$f: W\to \bf H^n$ holomorphic on $U$ and $W\setminus Q$
with prescribed  $(\partial f)^-(z)-(\partial f)^{+}(z)$ 
for each $z\in \partial Q$ and analogously for $f\in C^l$
with jump conditions for higher order derivatives.
In the case of $n>1$ there also can be used holomorphic extension
of holomorphic functions from proper quaternion submanifolds $K$
of $\partial Q$ (see Theorems 2(b) and 3(b) in \cite{adch}
and Theorem 4.1.11 \cite{henle} and \cite{lufsqv},
since a space of rational functions $f: K\to T_yN$ such that
$f|_K$ are holomorphic is infinite-dimensional, see also
Corollaries 3.4 and 3.5 in \cite{huya}).
Using charts of the atlas of $M$ we get that ${\sf H}_{\Upsilon }(M,N)$
is infinite dimensional.
\par In view of the Stone-Weierstrass theorem \cite{fell} the space
${\sf H}_{\Upsilon }(M,Y)$ is dense in $C^0(M,Y)$, hence
${\sf H}_{\Upsilon }(M,N)$ is dense in $C^0(M,N)$.
In the general case ${\sf H}_{\Upsilon ,a,c}(M,N)$
we finish the proof using the standard procedure of an increasing sequence
of $C^{\infty }$-domains $W_n$ in a quadrant
$Q$ with $dim_{\bf H}Q < \infty $ such that $cl (\bigcup_nW_n)=Q$.
The proof in the complex case is analogous due to Theorem VI.9 \cite{chau}.
\par {\bf 2.1.7. Theorems. (1).} {\it A unform space
$(S^MN)_{\xi }$ is a complete monoid and there exists a generated by it
topological group $(L^MN)_{\xi }=:G$
for $\xi =(\Upsilon ,a)$ or $\xi =(\Upsilon ,a,c)$ from \S 2.1.5
which is complete separable Abelian. 
Moreover, in $G$ there is a dense subgroup  
$(L^MN)_{\sf O}$ for complex manifolds or $(L^MN)_{\sf H}$
for quaternion manifolds with $\xi =(\Upsilon ,a)$;
$G$ is non-discrete non-locally compact and locally connected
for $dim_{\bf R} N>1$.}
\par {\bf (2).} {\it  A uniform space $X^{\xi }(M,N):=T_e(L^MN)_{\xi }$  
is a Hilbert space for each $1\le m=dim_{\bf K}M\le \infty $
and $dim_{\bf R} N >1$,
where ${\bf K}=\bf R$ or $\bf C$ or $\bf H$ respectively.}
\par {\bf (3.)} {\it Let $N$ be a Hilbert 
$Y^{\xi '}$-manifold over $\bf K$ with $a>a'$ and $c>c'$ 
for $\xi '=(\Upsilon ,a')$ or $\xi '=(\Upsilon ,a',c')$ respectively,
then there exists a mapping 
$\tilde E: {\tilde T}(L^MN)_{\xi }\to (L^MN)_{\xi }$ 
defined by $\tilde E_{\eta
}(v)=\exp_{\eta (s)}\circ v_{\eta }$ on a neighbourhood $V_{\eta }$ 
of the zero section in $T_{\eta }(L^MN)_{\xi }$ and it is 
a $C^{\infty }$-mapping for $Y^{\xi '}$-manifold $N$ by $v$ onto
a neighbourhood $W_{\eta }=W_e\circ \eta $ of $\eta \in (L^MN)_{\xi }$;
$\tilde E$ is the uniform isomorphism of 
uniform spaces $V_{\eta }$ and $W_{\eta
}$, where $s\in M$, $e$ is a unit element in $G$, $v\in V_{\eta },$
$1\le m\le \infty $.}
\par {\bf (4).} {\it For complex or quaternion manifolds
a group $(L^MN)_{\xi }$ is a closed proper subgroup in
$(L^M_{\bf R}N_{\bf R})_{\xi }$ of infinite codimension of
$T_e(L^MN)_{\xi }$ in $T_e(L^M_{\bf R}N_{\bf R})_{\xi }$.}
\par {\bf Proof.} Consider the orientable quaternion manifolds.
For $\xi =(\Upsilon ,a)$ or $\xi =(\Upsilon ,a,c)$ classes $<f>_{\xi }$
are closed due to Lemma 2.1.6.2 for the considered class 
of smoothness, since the uniform spaces $Y^{\xi }(M,s_0;N,y_0)$ are complete.
Via the construction of an Abelian group from an
Abelian monoid with unit and cancellation property we get
the loop groups $(L^MN)_{\xi }$ and $(L^M_{\bf R}N)_{\xi }$
respectively for finite dimensional $M$, where $\xi = (\Upsilon ,
a)$ (see \S 2.1.4.1). There exists a strict inductive limit of loop groups
$(L^{M_m}N)_{\Upsilon ,a}=:L^m$, since
there are natural embeddings $L^m\hookrightarrow L^{m+1}$,
such that each element $f\in Y^{\Upsilon ,a}(M_m,X_{0,m};N,y_0)$ is
considered in $Y^{\Upsilon ,a}(M_{m+1},X_{0,m+1};N,y_0)$
as independent from $(z^{n(m)+1},...,z^{n(m+1)-1})$ in the local normal
coordinates $(z^1,...,z^{n(m+1)})$ of $M_{m+1}$.
We denote it $str-ind_mL^m=:(L^MN)_{\Upsilon ,a}$ (see \S 2.1.5.2).
Via the construction of an Abelian group from
an Abelian monoid with unity and cancellation property we get
loop groups $(L^MN)_{\xi }$ for $\xi = (\Upsilon ,a,c)$ also.
The space $T_e(L^{M_m}N)_{\sf H}$ is linear over $\bf H$,
where $e$ is the unit element of the groop $(L^{M_m}N)_{\sf H}$.
Then in particular $X^{\xi }(M_m,N)$ is
the Banach space with $\| f\|_{X^{\xi }(M_m,N)}=\inf_{y\in f}{\|
y\|'}_{\xi }$, where $f=<y>_{\xi }$, $y\in Y^{\xi }(M,s_0;X,0)$.  On
the other hand, $X^{\xi }(M_m,N)$ is isomorphic with the completion of
$T_e(L^{M_m}N)_{\sf H}$ by the norm
$\| f\|_{X^{\xi }(M_m,N)}$.  Then $(\rho
_{k, \gamma }(y^1+y^2))^2+(\rho _{k, \gamma }(y^1-y^2))^2= 2[(\rho _{k,
\gamma }(y^1))^2+(\rho _{k, \gamma }(y^2))^2]$ for
each $y^1, y^2\in H^{\gamma }(M,X;Z_k)$ due to the 
choices of $\nu $ and $\rho
_{k, \gamma }$. If $y\in H^{\gamma }(M,X;Z_k)$ then $\rho _{k, \gamma }
(y)=\rho _{l, \gamma }(y)$ for each
$l>k$, since $\nu ({M'}_l)=0$ and $y\in H^{\gamma }(M,X;Z_l)$.  
For each $y^1, y^2 \in {H^{\gamma }}_p (M,s_0;X,0)$ there exists 
$Z\in \Upsilon $ such that
$y^1, y^2 \in H^{\gamma }(M,X;Z)$.  Therefore, from
\S 2.1.5 it follows that ${\| f_1+f_2\|^2}_{X^{\xi }(M,N)}+ 
{\| f_1-f_2\|^2}_{X^{\xi
}(M,N)}=2[{\| f_1\|^2}_{X^{\xi }(M,N)} +{\| f_2\|^2}_{X^{\xi
}(M,N)}]$ for each $f_1, f_2 \in X^{\xi }(M,N)$. 
Then ${\| f_1+f_2\|^*}_k^2+ {\| f_1-f_2\|^*}_k^2=
2[{\| f_1\|^*}_k^2+{\| f_2\|^*}_k^2]$
for each $0\le k\in \bf Z$ and each $f_1$ and $f_2\in 
Q^{\infty }_{\Upsilon ,a}(X)$,
consequently, $\| f_1+f_2\|_{X^{\xi }(M,N)}^2+
\| f_1-f_2\|_{X^{\xi }(M,N)}^2=2[
\| f_1\|_{X^{\xi }(M,N)}^2+\| f_2\|_{X^{\xi }(M,N)}^2]$. Hence the formula
$4(f_1,f_2):= { \| f_1+f_2\|^2}_{X^{\xi }(M,N)} -
{\| f_1 - f_2\|^2}_{X^{\xi }(M,N)} + [ \sum_{l=1}^3
i_l { \| f_1 + i_l f_2 \|^2}_{X^{\xi }(M,N)}
- i_l { \| f_1 - i_l f_2\|^2}_{X^{\xi }(M,N)} ]/3 $
gives the scalar product $(f_1,f_2)$ in 
$X^{\xi }(M,N)$ and this is the Hilbert space over $\bf H$,
where $i_l\in \{ i, j, k \} $, $l=1, 2, 3$.
\par The spaces $Y^{\xi }(M,s_0;N,y_0)$ and 
$X^{\xi }(M,N)$ are complete, consequently, $G$ is complete. 
The space $X^{\xi }(M,N)$ is separable and $\Lambda _N \subset
\bf N$, consequently, $G$ is separable.
The composition and the inversion in
$(L^{M_m}N)_{\sf H}$ induces these operations in $G$, that are continuous 
due to Theorem 2.1.4.1 and using completions we get, that
$G$ is the Abelian topological group.
\par Let $\beta :  M_m\to
N$ be a $Y^{\xi }$-mapping such that $\beta (s_0)=y_0$. 
If $C_0$
is the connected component of $y_0$ in $N$ then $\beta (M_m) \subset C_0$.
On the other hand, $N$ was supposed to be connected.
In view of Theorems about extensions of functions of different classes of
smoothness \cite{seel,touger} 
and using completions in the described above
spaces  there exists a neighbourhood $W$ of $w_0$ such that
for each $f:  M_m\times \{ 0, 1 \}\to N$ of class $Y^{\xi }$ with
$f(M_m,0)=\{ y_0 \} $ and $f(s,1)=\beta (s)$ for each $s\in M$ there exists
its $Y^{\xi }$-extension $f:  M_m\times [0,1]\to N$, where $\{ 0, 1\} :=\{
0\} \cup \{ 1\} $, $\beta \in W$,
since there exists a neighbourhood $V_0$ of $y_0$ in $N$
such that it is $C^0([0,1]\times V_0,N)$-contractible into a point.
Hence for each class $<\beta >_{\xi }$ 
in a sufficiently small (open) neighbourhood of $e$ there exists a
continuous curve $h:  [0,1]\to G$ such that 
$h(0)=e$ and $h(1)=<\beta >_{\xi }$.
\par In view of Theorem 2.1.4.1 and Lemma 2.1.6.2
the tangent space $T_eG$ is infinite dimensional over $\bf H$, 
consequently, $G$ is not locally compact, where
$e$ is the unit element in $G$.  
\par Let $\nabla $ be a covariant
differentiation in $N$ corresponding to the Levi-Civit\'a 
connection in $N$ due to Theorem 5.1 \cite{flas}.  
This is possible, since $N$ is the Hilbert manifold $C^{\infty }$-manifold
over $\bf H$ and hence over $\bf R$, consequently, it
has the partition of unity \cite{lan}.  
Therefore, there exists the
exponential mapping $\exp:  {\tilde T}N\to N$ such that for each 
$z\in N$ there are a ball
$B(T_zN,0,r):=\{ y\in T_zN:  \| y\|_{T_zN}\le r \} $ 
and a neighbourhood $S_z$
of $z$ in $N$ for which $\exp_z:  B(T_zN,0,r)\to S_z$ 
is the homeomorphism, since
$\phi $ is in the class of smoothness $Y^{\xi }$ due to Theorem IV.2.5
\cite{lan}, where $\exp_z w = \phi (1)$, $\phi (q)$ is a geodesic, $\phi :
[0,1] \to N$, $d\phi (q)/dq|_{q=0}=w$, $\phi (0)=z$, 
$w\in B(T_zN,0,r)$ \cite{kling}. 
\par In view of Theorems 5.1 and 5.2 \cite{eliass} a mapping
\par $(i)$ $E:  {\tilde T}Y^{\xi }(M_m,s_0;N,y_0)\to
Y^{\xi }(M_m,s_0;N,y_0)$ is a local isomorphism, since
$a> b$, so $\sum_{k=1}^{\infty }k^2(k!)^{b-a} <\infty $, where 
\par $(ii)$ $E_g(h):=\exp_{g(s)}\circ h_{g}$, $s\in M_m$, $h\in TY^{\xi
}(M_m,s_0;N,y_0)$, $h_g\in T_gY^{\xi }(M_m,s_0;N,y_0)$, $g\in Y^{\xi
}(M_m,s_0;N,y_0) .$ 
In view of \cite{ebi,eliass} the tangent bundle $TY^{\xi }
(M_m,s_0;N,y_0)$ is isomorphic with $Y^{\xi }(M_m,s_0;TN, y_0\times \{ 0\}
)\times T_{y_0}N.$ 
Then $E$ induces $\tilde E$ with the help of factorisation by 
$R_{\xi }$ and the subsequent construction of the loop group from
the loop monoid. This mapping $\tilde E$ is of class of smoothness 
$C^{\infty }$ as follows
from equations for geodesics (see \S IV.3 \cite{lan}), since $TTN$ is the
$Y^{\xi "}$-manifold with $a>a">0$ and $c>c">0$. 
Indeed, this construction at first may be applied for
$(L^{M_m}N)_{\sf H}$ and then using the completion to $(L^MN)_{\xi }$.
\par The uniform space of quaternion holomorphic mappings from $M$ into
$N$ is proper in the set of infinite Fr\'echet differentiable mappings
from $M$ into $N$. In view of \S 2.1.5 uniformities in
groups $(L^MN)_{\xi }$ and $(L^M_{\bf R}N_{\bf R})_{\xi }$ are equivalent.
Therefore, $(L^MN)_{\xi }$ has an embedding into
$(L^M_{\bf R}N_{\bf R})_{\xi }$ as a closed proper subgroup
of infinite codimension of $T_e(L^MN)_{\xi }$ in
$T_e(L^M_{\bf R}N_{\bf R})_{\xi }$, since $dim_{\bf R} N>1$.
The proofs in the complex and orientable real cases are analogous.
\par It remains the case of nonorientable real or quaternion manifolds
which can be deduced also from Theorems 2.1.8 below and the case
of orientable manifolds.
\par {\bf 2.1.8. Theorems.} {\it Suppose that both either real
or quaternion manifolds $M$ and $N$ together with their covering manifolds 
$\tilde M$ and $\tilde N$ satisfy the conditions imposed in \S 2.1.5.
\par $\bf (1)$. Let $N$ be a nonorientable manifold and let
$\theta _N: \tilde N\to N$ be a quotient mapping of its
double covering manifold $\tilde N$. 
Then there exists a quotient group homomorphism 
$\tilde \theta _N : (L^M{\tilde N})_{\xi }\to (L^MN)_{\xi }$.
\par $\bf (2)$. Let $M$ be a nonorientable manifold,
then a quotient mapping $\theta _M: \tilde M\to M$ induces a 
group embedding $\tilde \theta _M:
(L^MN)_{\xi }\hookrightarrow (L^{\tilde M}N)_{\xi }.$}
\par {\bf Proof.} If $M$ is a nonorientable manifold,
then there exists the homomorphism $h$ of the fundamental group
$\pi _1(M,s_0)$ onto the two-element group $\bf Z_2$.
For connected $M$ the group $\pi _1(M,s_0)$ does not depend 
on a marked point $s_0$ and it is denoted by $\pi _1(M)$.
If $M$ is a connected manifold, 
then it has a universal covering manifold $M^*$ which is 
linearly connected and it has a fiber bundle with the group
$\pi _1(M)$ and a projection $p: M^*\to M$.
Using the homomorphism $h$ one gets the orientable double covering
$\tilde M$ of $M$ such that $\tilde M$ is connected, if
$M$ is connected (see Proposition 5.9 \cite{kob} and Theorem 78
\cite{pont}). Moreover, for each $x\in M$ there exists
a neighborhood $U$ of $x$ such that 
$\theta _M^{-1}(U)$ is the disjoint union of two diffeomorphic
open subsets $V_1$ and $V_2$ in $\tilde M$, where $\theta _M: \tilde M\to M$
is the quotient mapping, $g: V_1\to V_2$ is a diffeomorphism. 
\par $(1)$. For each $\tilde f \in Z^{\Upsilon ,a,c}(M,\tilde N)$
there exists $f=\theta _N\circ \tilde f$ in $Z^{\Upsilon ,a,c}(M,N)$.
This induces the quotient mapping $\bar \theta _N :
Z^{\Upsilon ,a,c}(M,\tilde N)\to Z^{\Upsilon ,a,c}(M,N)$,
hence it induces the quotient mapping $\bar \theta _N :
Z^{\Upsilon ,a,c}(M\vee M,\tilde N)\to Z^{\Upsilon ,a,c}(M\vee M,N)$
such that $\bar \theta _N (f\vee h) =
\bar \theta _N(f)\vee \bar \theta _N(g)$.
Considering the equivalence relation in $Y^{\xi }(M,s_0;\tilde N,y_0)$
and then loop monoids we get the quotient homomorphism
$(S^M\tilde N)_{\xi }\to (S^MN)_{\xi }$. With the help of construction
of the loop group from the loop monoid it induces the loop groups
quotient homomorphism.
\par $(2)$. On the other hand, let $M$ be a nonorientable manifold
then the quotient mapping $\theta _M: \tilde M\to M$ induces a 
locally finite open covering $\{ U_x: x\in M_0 \} $ of $M$,
where $M_0$ is a subset of $M$, such that each $\theta _M^{-1}(U_x)$
is the disjoint union of two open subsets $V_{x,1}$ and $V_{x,2}$ 
in $\tilde M$ and there exists a diffeomorphism $g_x$ of $V_{x,1}$
on $V_{x,2}$ of the same class of smoothness as $M$.
This produces the closed subspace of all $\tilde f\in  Z^{\Upsilon ,a,c}(
\tilde M,N)$ for which 
\par $(i)$ $\tilde f|_{V_{x,1}}(g_x^{-1}(y))=\tilde f|_{V_{x,2}}(y)$
for each $y\in V_{x,2}$ and for each $x\in M_0$, 
where $M_0=M_0^f$ and $\{ U_x=U_x^f: x\in M_0 \} $ may be dependent on $f$. 
If $s_0$ is a marked point in $M$, then let one of the points $\tilde s_0$
of $\theta _M^{-1}(s_0)$ be the marked point in $\tilde M$.
Then $\theta _M$ induces the quotient mapping 
$\theta _M: \tilde M\vee \tilde M \to M\vee M$.
If both $M$ and $\tilde M$ 
satisfy the imposed above conditions on manifolds, then
this induces the embedding
$\bar \theta _M: Z^{\Upsilon ,a,c}(M,N)\hookrightarrow
Z^{\Upsilon ,a,c}(\tilde M,N)$. 
The identity mapping $id(x)=x$ for each $x\in \tilde M$
evidently satisfy Condition $(i)$.
If $\tilde f$ is the diffeomorphism of $\tilde M$
satisfying Condition $(i)$, then
applying $\tilde f^{-2}$ to both sides of the equality
we see, that it is satisfied for $\tilde f^{-1}$
with the same covering $\{ U_x: M_0 \} $.
If $\tilde f$ and $\tilde h$ are two diffeomorphisms of $\tilde M$,
then there exists $\{ U_x^h: x\in M_0^h \} $ such that
$\{ \tilde f^{-1} (\theta _M^{-1}(U_x^h)): x\in M_0^h \} $
is the locally finite covering of $\tilde M$. 
Two manifolds $M$ and $\tilde M$ are metrizable, consequently, paracompact
(see Theorem 5.1.3 \cite{eng}).
Due to paracompactness of $M$ and $\tilde M$
there exists a locally finite covering
$\{ U_x^{h\circ f}: x\in M_0^{h\circ f} \} $ 
for which Condition $(i)$ is satisfied, since 
$\{ U_x^f\cap \theta _M\circ \tilde f^{-1} 
(\theta _M^{-1}(U_z^h)): x\in M_0^f, z\in M_0^h \} $
has a locally finite refinement.
This means, that $\bar \theta _M: Diff^{\infty }(M_{\bf R})\hookrightarrow
Diff^{\infty }(\tilde M_{\bf R})$ is the group embedding. 
In a complete uniform space $(X,{\sf U})$ for its subset $Z$ 
a uniform space $(Z,{\sf U}_Z)$ is complete if and only if
$Z$ is closed in $X$ relative to the topology induced by $\sf U$
(see Theorem 8.3.6 \cite{eng}).
Since both groups are complete and the uniformity of 
$Diff^{\infty }(\tilde M_{\bf R})$ induces the uniformity in 
$Diff^{\infty }(M_{\bf R})$ equivalent to its own, then
$\bar \theta _M( Diff^{\infty }(M_{\bf R}))$ is closed in
$Diff^{\infty }(\tilde M_{\bf R})$.
Considering the equivalence relation in
$Z^{\Upsilon ,a,c}(\tilde M,N)$ we get the embedding of loop monoids
$\tilde \theta _M :(S^MN)_{\xi }\hookrightarrow
(S^{\tilde M}N)_{\xi }$ (see \S \S 2.1.4 and 2.1.5).
This produces via the construction of an Abelian group from an
Abelian monoid with unity and cancellation property
the loop groups embedding (respecting their topological group structures)
$\tilde \theta _M:
(L^MN)_{\xi }\hookrightarrow (L^{\tilde M}N)_{\xi }.$
\par {\bf 2.2. Note.} For a diffeomorphism group we also consider
a compact complex manifold $M$.
For noncompact complex $M$, satisfying conditions of \S 2.1.1
and $(N1)$ a diffeomorphism group is considered
as consisting of diffeomorphisms $f$ of class $Y^{\Upsilon ,a,c}$
(see \S 2.1.5), that is, 
$(f_{i,j}-id_{i,j})\in Y^{\Upsilon ,a,c}(U_{i,j},\phi _i(U_i))$
for each $i, j$, $U_{i,j}$ is a domain of definition of 
$(f_{i,j}-id_{i,j})$ and then analogously to the real case
the diffeomorphism group $Diff^{\xi }(M)$
is defined, where $\xi =(\Upsilon ,a,c)$, 
$a=(a_1,a_2)$, $c=(c_1,c_2)$, $a_1 \le -1 $ and $c_1\le -1.$
This means that $Diff^{\xi }(M):=Y^{\xi ,id}(M,M)\cap Hom(M)$.
\par {\bf 2.2.1. Theorem.} {\it Let $M$ and $N$ be two complex manifolds,
then there exist quaternion manifolds $P$ and $Q$ such that
$M$ and $N$ have complex holomorphic emebeddings into $P$ and $Q$
respectively as closed submanifolds. Moreover, $(L^MN)_{\xi }$ is the proper
closed subgroup in $(L^PQ)_{\xi }$ and $Diff^{\xi }(M)$ is the closed
proper subgroup in $Diff^{\xi }(P)$ such that the codimensions
over $\bf R$ of $T_e(L^MN)_{\xi }$ in $T_e(L^PQ)_{\xi }$
and $T_eDiff^{\xi }(M)$ in $T_eDiff^{\xi }(P)$ are infinite.}
\par {\bf Proof.} At first prove that if $N$ is a complex manifold,
then there exists a quaternion manifold $Q$ and a
complex  holomorphic embedding $\theta : N\hookrightarrow Q$.
Suppose $At (N)=\{ (V_a,\psi _a): a \in \Lambda \} $
is any holomorphic atlas of $N$, where $V_a$ is open in $N$,
$\bigcup_aV_a=N$, $\psi _a: V_a\to \psi _a(V_a)\subset X=\bf C^n$
or $X=l_2({\bf C})$ is a homeomorphism for each $a$,
$n=dim_{\bf C}M\in \bf N$ or $n=\infty $ respectively,
$\{ V_a: a \in \Lambda \} $ is a locally finite
covering of $N$, $\psi _b\circ \psi _a^{-1}$ is a holomorphic function
on $\psi _a(V_a\cap V_b)$ for each $a, b\in \Lambda $
such that $V_a\cap V_b\ne \emptyset $.  For each complex holomorphic
function $f$ on an open subset $V$ in $X$ there exists a quaternion
holomorphic function $F$ on an open subset $U$ in $Y=\bf H^n$
or $Y=l_2({\bf H})$ respectively 
such that $\pi _{1,1}(U)=V$ and $F_{1,1}|_{Ve}=f|_V$,
where $\pi _{1,1}: Z e\oplus Z i\oplus Z j \oplus Z k\to Z e\oplus Z i=X$
is the projection with $Z=\bf R^n$ or $Z=l_2({\bf R})$ correspondingly,
$F_{1,1}:=\pi _{1,1}\circ F$
(see Proposition $3.13$ \cite{luoyst} and use a locally finite
covering of $V$ by balls). Therefore, for each two charts
$(V_a,\psi _a)$ and $(V_b,\psi _b)$ with $V_{a,b}:=
V_a\cap V_b\ne \emptyset $ there exists $U_{a,b}$ open in
$Y$ and a quaternion holomorphic function $\Psi _{b,a}$
such that $\Psi _{b,a}|_{\psi _a(V_{a,b})e}=\psi _{b,a}|_{\psi _a(V_{a,b})}$,
where $\psi _{b,a}:=\psi _b\circ \psi _a^{-1}$,
$\pi _{1,1}(U_{a,b})=\psi _a(V_{a,b})$.
Consider $F :=\bigoplus_aF_a$, where $F_a$ is open in $Y$,
$\pi _{1,1}(F_a)=\psi _a(V_a)$ for each $a\in \Lambda $.
The equivalence relation $\cal C$ in the topological space
$\bigoplus_a\psi _a(V_a)$ generated by functions $\psi _{b,a}$
has an extension to the equivalence relation $\cal H$
in $F$. Then $Q:=F/\cal H$ is the desired quaternion manifold
with $At (Q)=\{ (\Psi _a,U_a): a\in \Lambda  \} $
such that $\Psi _b\circ \Psi _a^{-1}=\Psi _{b,a}$ for each
$U_a\cap U_b\ne \emptyset $, $\Psi _a^{-1}|_{\psi _a(V_a)e}=
\psi _a^{-1}|_{\psi _a(V_a)}$ for each $a$, $\Psi _a^{-1}: F_a\to U_a$
is the quaternion homeomorphism.
Moreover, each homeomorphism $\psi _a: V_a\to \psi _a(V_a)\subset
X$ has the quaternion extension up to the homeomorphism
$\Psi _a: U_a\to \Psi _a(U_a)\subset Y$. The family of embeddings
$\eta _a: \psi _a(V_a)\hookrightarrow F_a$ such that
$\pi _{1,1}\circ \eta _a=id$ together with $At(M)$ induces the
complex holomorphic embedding $\theta : N\hookrightarrow Q$.
If $M$ is foliated then choose also $P$ foliated such that
$dim_{\bf C}M_m=dim_{\bf H}P_m=n(m)$ for each $m\in \bf N$ (see
\S 2.1.5).
\par In view of Proposition $3.13$ \cite{luoyst} there exists an embedding
of $Y^{\xi }(M,N)$ into $Y^{\xi }(P,Q)$ as the closed submanifold
such that a codimension over $\bf R$ of $T_fY^{\xi }(M,N)$
in $T_gY^{\xi }(P,Q)$ is infinite, where $g$ is a quaternion
extension of $g$. This is true even in the case when for these manifolds
underlying real manifolds $M_{\bf R}$ and $P_{\bf R}$,
$N_{\bf R}$ and $Q_{\bf R}$ are $C^{\infty }$-diffeomorphic,
since quaternion holomorphic functions need
not be a right or left superlinearly superdifferentiable and hence
need not in general satisfy complex Cauchy-Riemann conditions.
Uniformities in $Y^{\xi }(M,N)$ and in $Y^{\xi }(P,Q)$
are consistent, hence embeddings of corresponding groups of loops
and groups of diffeomorphisms are closed.
\par{\bf 2.2.2. Remarks and definitions.}
For investigations of representations of diffeomorphism groups
with the help of quasi-invariant transition measures induced by
stochastic processes at first there are given below necessary
definitions and statements on special kinds of diffeomorphism groups
having Hilbert manifold structures.
\par Let $M$ and $N$ be real manifolds on $\bf  R^n$ or $l_2$ 
and satisfying Conditions 2.2.(i-vi) \cite{lurim1}
or they may be also canonical closed submanifolds of that of
in \cite{lurim1}.
For a field ${\bf K}=\bf R$, $\bf C$ or $\bf H$
let $l_{2,\delta }({\bf K})$
be a Hilbert space of vectors $x=(x^j: x^j\in {\bf K}, j\in {\bf N} )$
such that $\| x \| _{l_{2,\delta }}:=\{ \sum_{j=1}^{\infty }
|x^j|^2j^{2\delta } \} ^{1/2}<\infty $.
For $\delta =0$ we omit it as the index.
Let also $U$ be an open subset in $\bf R^m$ and $V$ be an open subset in
$\bf R^n$ or $l_2$ over $\bf R$, 
where $0\in U$ and $0\in V$ with $m$ and $n\in \bf N.$
By $H^{l, \theta }_{\beta ,\delta }(U,V)$
is denoted the following completion relative to the metric
$q^l_{\beta ,\delta }(f,g)$ of the family of all strongly
infinite differentiable functions 
$f,g: U \to V$ with $q^l_{\beta ,\delta }(f,\theta )<\infty $, where
$\theta \in C^{\infty }(U,V)$, $0\le l\in \bf Z$,
$\beta \in \bf R$, $\infty >\delta \ge 0$, $q^l_{\beta ,\delta }(f,g)
:=(\sum_{0 \le |\alpha |\le l} \| \bar m^{\alpha \delta }
<x> ^{\beta +|\alpha |}D_x^{\alpha }
(f(x)-g(x)) \|^2_{L^2})^{1/2}$, $L^2:=L^2(U,F)$
(for $F:=\bf R^n$ or $F=l_{2, \delta }=l_{2,\delta }({\bf R})$)
is the standard Hilbert space of all classes of equivalent 
measurable functions $h: U \to F$ for which there exists
a finite norm $\| h\|_{L^2}:=(\int_U|h(x)|_F^2
\mu _m(dx)))^{1/2} <\infty $, $\mu _m$ denotes the Lebesgue measure
on $\bf R^m$. 
\par Let also $M$ and $N$ have finite atlases
such that $M$ be on
$X_M:=\bf R^m$ and $N$ on $X_N:=\bf R^n$ or $X_N:=l_2$, 
$\theta : M\hookrightarrow N$ be a $C^{\infty }$-mapping,
for example, embedding.
Then $H^{l,\theta }_{\beta ,\delta }(M,N)$ denotes the completion
of a family of all $C^{\infty }$-functions $g, f:
M \to N$ with $\kappa ^l_{\beta ,\delta }(f,\theta )< \infty $,
where the metric is given by the following formula
$\kappa ^l_{\beta ,\delta }(f,g)=(\sum_{i,j}[q^l_{
\beta ,\delta }(f_{i,j},
g_{i,j})]^2)^{1/2}$, where $f_{i,j}:=\psi _i\circ f\circ \phi _j^{-1}$ 
with domains $\phi _j(U_j)\cap \phi _j(f^{-1}(V_i))$,
$At(M):= \{ (U_i,\phi _i): i \} $ and $At(N):= \{ (V_j,\psi _j): j \} $
are atlases of $M$ and $N$, $U_i$ is an open subset in $M$ for each $i$
and $V_j$ is an open subset in $N$ for each $j$, $\phi _i: U_i\to X_M$ and 
$\psi _j: V_j\to X_N$ are homeomorphisms of $U_i$ on $\phi _i(U_i)$ and
$V_j$ on $\psi _j(V_j)$, respectively. 
Hilbert spaces $H^{l,\theta }_{\beta ,\delta }(U,F)$
and $H^l_{\beta ,0}(TM)$ are called weighted Sobolev spaces, 
where $H^l_{\beta ,\delta }(TM):=\{  f:M\to TM:$ $f \in H^l_{\beta ,\delta
}(M,TM),$ $\pi \circ
f(x)=x\mbox{ for each } x \in M \} $ with $\theta (x)=(x,0)\in T_xM$
for each $x\in M$. From the latter definition it follows, that
for such $f$ and $g$ there exists $\lim_{R\to \infty }q^l_{\beta ,\delta }
(f|_{U_R^c}, g|_{U_R^c})=0$, when $(U,\phi )$ is a chart Hilbertian 
at infinity, $U_R^c$ is an exterior of a ball of radius $R$ 
in $U$ with center in the fixed point $x_0$
relative to the distance function $d_M$ in $M$
induced by a Riemannian metric $g$ (see \S 2.2(v) \cite{lurim1}). 
For $\beta =0$ or $\gamma =0$ we omit $\beta $
or $\gamma $ respectively in the notation 
$Dif^t_{\beta ,\gamma }(M):=H^{t,id}_{\beta ,\gamma }(M,M)\cap Hom(M)$ 
and $H^{l,\theta }_{\beta ,\gamma }$.
\par The uniform space $Dif^t_{\beta ,\gamma }(M)$ has the group structure
relative to the composition of diffeomorphisms and is called 
the diffeomorphism group, where $Hom(M)$ is the group of all homeomorphisms 
of $M$. 
\par Each topologically adjoint space
$(H^l_{\beta }(TM))'=:H^{-l}_{-\beta }(TM)$ also is the Hilbert space
with the standard norm in
$H'$ such that $\| \zeta \|_{H'}=\sup_{\| f\| _H=1}| \zeta (f)|$.
\par  {\bf 2.3. Diffeomorphism groups of Gevrey-Sobolev classes of 
smoothness. Notes and definitions.} 
Let $U$ and $V$ be open subsets in the Euclidean space
$\bf R^k$ with $k\in \bf N$ or in the standard separable Hilbert 
space $l_2$ over $\bf R$,
$\theta : U\to V$ be a $C^{\infty }$-function (infinitely strongly 
differentiable), $\infty >\delta \ge 0$ be a parameter.
There exists the following metric space 
$H^{ \{ l \} ,\theta }_{ \{ \gamma \} ,\delta }(U,V)$
as the completion of a space of all functions
$Q:= \{ f:$ $f \in E^{\infty ,\theta }
_{\infty ,\delta }(U,V),$ $\mbox{ there exists }n \in {\bf N}
\mbox{ such that }$
$supp (f) \subset U\cap {\bf R^n}, d_{ \{ l\}, \{ \gamma \},\delta }
(f,\theta )<\infty \} $
relative to the given below metric
$d_{ \{ l\}, \{ \gamma \},\delta }:$ 
$$(i)\mbox{ }d_{ \{ l\}, \{ \gamma \},\delta }(f,g):=
\sup_{x \in U}(\sum_{n=1}^{\infty }(\bar \rho ^l_{\gamma ,n,\delta }
(f,g))^2)^{1/2}<\infty $$ and 
\par $\lim_{R \to \infty }d_{\{ l\},\{ \gamma \},\delta }(f|_{U_R^c},
g|_{U_R^c})=0$, when $U$ is a chart Euclidean or Hilbertian 
correspondingly at infinity,
$f$ as an argument in $\bar \rho ^l_{\gamma ,n,\delta }$ 
is taken with the restriction on
$U\cap \bf R^n$, that is,  $f|_{U\cap \bf R^n}: U\cap {\bf R^n}\to
f(U)\subset V$ (see also \S \S 2.1-2.5 \cite{lurim1} and \cite{lurim2}
about $E^{t,\theta }_{\beta ,\gamma }$), 
$\bar \rho ^l_{\gamma ,n,\delta }(f,id)^2$ 
$:=\omega _n^2(\kappa ^l_{\gamma ,\delta }(f|_{(U\cap {\bf R^n})},
id|_{(U \cap {\bf R^n})})^2 -$
$\kappa ^{l(n-1)}_{\gamma (n-1),\delta }(f|_{(U\cap {\bf R^{n-1}})},
id|_{(U \cap {\bf R^{n-1}})})^2)$ 
for each $n>1$ and \\
$\bar \rho ^l_{\gamma ,1,\delta }(f,id)
:=\omega _1(\kappa ^l_{\gamma ,\delta }(f|_{(U\cap {\bf R^1})},
id|_{(U \cap {\bf R^1})})$ 
with $q^l_{\gamma ,\delta }$
and the corresponding terms $\kappa ^l_{\gamma ,\delta }$ 
from \S 2.2.2,
$l=l(n)>n+5$, $\gamma =\gamma (n)$ and $l(n+1)\ge l(n)$
for each $n$, $l(n)\ge [t]+sign \{ t \}+[n/2]+3$, $\gamma
(n)\ge \beta +sign \{ t \} +[n/2]+7/2$, $\omega _{n+1}\ge n\omega _n\ge 1$. 
Moreover, $\bar \rho ^l_{\gamma ,n,\delta }
(f,id)(x^{n+1},x^{n+2},...)\ge 0$ is the metric by variables
$x^1,...,x^n$ in $H^l_{\gamma ,\delta }(U\cap
{\bf R^n}, V)$ for $f$ as a function by 
$(x^1,...,x^n)$ such that $\bar \rho ^l_{\gamma ,n,\delta }$
depends on parameters $(x^j: j>n)$. The index $\theta $ 
is omitted when $\theta =0$.
The series in $(i)$ terminates $n\le k$, when $k\in \bf N$.
\par Let for $M$ connecting mappings of charts be such that
$(\phi _j\circ \phi _i^{-1}-id_{i,j})
\in H^{ \{ l' \} }_{ \{ \gamma ' \} ,\chi }(U_{i,j},l_2)$ for each $U_i
\cap U_j\ne \emptyset $ and the Riemannian metric $g$ be of class
of smoothness $H^{
\{ l' \} }_{ \{ \gamma ' \} ,\chi }$, where subsets $U_{i,j}$ are open in
$\bf R^k$ or in $l_2$ correspondingly
domains of $\phi _j \circ \phi _i^{-1}$, $l'(n)\ge l(n)+2$,
$\gamma '(n)\ge \gamma (n)$ for each $n$, $\infty >\chi \ge \delta $, 
submanifolds $\{ M_k: k=k(n), n \in {\bf N} \}$
are the same as in Lemma 3.2 \cite{lurim1}. 
Let $N$ be some manifold satisfying analogous conditions as $M$.
Then there exists the following uniform space
$H^{ \{ l \} ,\theta }_{ \{ \gamma \} ,\delta ,\eta }(M,N):=
\{ f\in E^{\infty ,\theta }_{\infty ,\delta }(M,N)|
(f_{i,j}-\theta _{i,j})\in H^{ \{ l \} ,\theta }_{ \{ \gamma \} ,\delta }
(U_{i,j},l_2)$ $\mbox{for each charts }$
$\{ U_i,\phi _i \}$ and $\{ U_j,\phi _j \}$ $\mbox{ with }U_i\cap U_j
\ne \emptyset , \chi _{\{ l\},\{ \gamma \},\delta ,\eta }
(f,\theta )<\infty $ $\mbox{and }\lim_{R\to \infty }
\chi _{\{ l\},\{ \gamma \},\delta ,\eta }(f|_{M^c_R},\theta 
|_{M^c_R})=0 \} $ and there exists the corresponding diffeomorphism group
$Di^{ \{ l \} }_{ \{ \gamma \} , \delta ,\eta }(M):=$ $ \{ f : f\in Hom(M),$
$f^{-1}$ $\mbox{ and }$ $f\in
H^{ \{ l \} ,id}_{ \{ \gamma \} ,\delta ,\eta }(M,M) \} $ 
with its topology given by the following left-invariant metric
$\chi _{\{ l\}, \{ \gamma \},\delta ,\eta }(f,g):=\chi _{ \{ l\},
\{ \gamma \}, \delta ,\eta }(g^{-1}f,id)$,
$$(ii)\mbox{ }\chi _{\{ l\}, \{ \gamma \}, \delta ,\eta }
(f,g):=(\sum_{i,j}(d_{ \{ l\}, \{ \gamma \}, \delta }(f_{i,j},
g_{i,j})i^{\eta }j^{\eta })^2)^{1/2}<\infty ,$$  
$g_{i,j}(x)\in l_2$ and $f_{i,j}(x) \in l_2$,
$\phi _i(U_i)\subset l_2$, $U_{i,j}=U_{i,j}(x^{n+1},x^{n+2},...)
\subset l_2$ is a domain of $f_{i,j}$ by variables 
$x^1,...,x^n$ for chosen variables $(x^j: j>n)$
due to the given foliations in $M$,
$U_{i,j} \subset {\bf R^n}\hookrightarrow l_2$, when $(x^j: j>n)$
are fixed and $U_{i,j}$ is a domain in $\bf R^n$ by variables
$(x^1,...,x^n)$, where $\infty >\eta \ge 0$.
\par In particular, for the finite dimensional manifold
$M_n$ the group $Di^{ \{ l\} }_{ \{ \gamma \} ,\delta ,\eta }(M_n)$ 
is isomorphic to the diffeomorphism group 
$Dif^l_{\gamma ,\delta }(M_n)$
of the weighted Sobolev class of smoothness
$H^l_{\gamma ,\delta }$ with $l=l(n)$, $\gamma =\gamma (n)$,
where $n=dim_{\bf R}(M_n)<\infty $.
\par If manifolds $M$ and $N$ are both supplied (or considered with)
either complex or quaternion structures with corresponding
foliations in the infinite dimensional case, then in the definitions
of $H^{ \{ l \} }_{ \{ \gamma \} , \delta , \eta } (M,N)$,
$Dif^t_{\beta , \gamma }(M)$, $Di ^{ \{ l \} }_{ \{ \gamma \} , \delta ,
\eta }(M)$ impose conditions that these uniform spaces are completions
relative to their uniformities of complex or quaternion holomorphic
mappings respectively, such that their elements $f$ satisfy in the
generalized sense condition ${\bar \partial }f=0$ in the complex case
and ${\tilde \partial }f=0$ in the quaternion case. This generalized
sense of complex or quaternion holomorphicity is induced by such conditions
in the sense of distributions, since each Sobolev space $H$ on
finite dimensional manifold $M_n$ for functions with values in
a $\bf K$-Hilbert space $X$ has the topologically conjugated space $H^*$
relative to the complex or quaternion Hermitian inner product
correspondingly. Elements of the above mentioned uniform spaces
are therefore generalized functions on the subspaces
of complex or quaternion holomorphic functions in $H^*$.
In the case of a manifold $N$ modelled on $X$ over $\bf K$
use an atlas of $N$ and transition mappings of charts in the standard
way to get holomorphicity conditions in the sense of distributions.
\par {\bf 2.4. Remarks.} Let two sequences be given
$\{ l \}:=\{ l(n): n \in {\bf N} \}\subset
\bf Z$ and $\{ \gamma \} :=\{ \gamma (n): n \in {\bf N} \}\subset
\bf R$, where manifolds $M$ and $\{ M_k: k=k(n), n \in {\bf N} \} $ are
the same as in \S \S 2.2 and 2.3. Then there exists the following space 
$H^{ \{ l \} ,\theta }_{ \{ \gamma \} ,\delta ,\eta }(M,TN)$.
By $H^{ \{ l \} ,\theta }_{ \{ \gamma \} ,\delta ,\eta }(M|TN)$
it is denoted its subspace of all functions $f: M\to TN$
such that $\pi _N(f(x))=\theta (x)$ for each $x\in M$,
where $\pi _N: TN\to N$ is the natural projection,
that is, each such 
$f$ is a vector field along $\theta $, $\theta : M\to N$ is a
fixed $C^{\infty }$-mapping. For $M=N$ and $\theta =id$ the metric space
$H^{ \{ l \} ,\theta }_{ \{ \gamma \} ,\delta ,\eta }(M|TM)$
is denoted by
$H^{ \{ l \} }_{ \{ \gamma \} ,\delta ,\eta }(TM)$. Spaces
$H^{ \{ l \} ,id}_{ \{ \gamma \} ,\delta ,\eta }(M_k|TN)$
and $H^{ \{ l \} }_{ \{ \gamma \} ,\delta ,\eta }(TM)$
are Banach spaces with the norms 
$\| f\|_{ \{ l\}, \{ \gamma \} ,\delta ,\eta }
:=\chi _{\{ l\}, \{ \gamma \}, \delta ,\eta }(f,f_0)$
denoted by the same symbol, where $f_0(x)=(x,0)$ and $pr_2f_0(x)=0$
for each $x\in M$. This definition can be spread on the case
$l=l(n)<0$, if take $\sup_{\| \tau \| =1}| <x>_m^{| \alpha |- \gamma (m)}$ 
$(D_x^{\alpha }\tau _{i,j},
[\zeta _{i,j}-\xi _{i,j}])_{L^2(U_{i,j,m},l_{2,\delta })}|$
instead of $\| <x>_m^{\gamma (m)+|
\alpha |}$ $D_x^{\alpha }(\zeta _{i,j}-\xi _{i,j})(x) \| _{L^2(
U_{i,j,m}, l_{2,\delta })}$, where
$\tau \in H_{-\gamma } ^{-l}(M_k|TN)$, $<x>_m=(1+\sum_{i=1}^m(x^i)^2)^{1/2}$,
$U_{i,j,m}=U_{i,j,m}(x^{m+1},x^{m+2},...)$ denotes the 
domain of the function $\zeta _{i,j}$
by $x^1,...,x^m$ for chosen $(x^j:$ $j>m)$, 
$\| \zeta \| _k(x)$ are functions by variables $(x^i: i>k)$.
Further the traditional notation is used:  
$sign( \epsilon )=1$ for $\epsilon >0$, $sign(\epsilon )=-1$
for $\epsilon <0$, $sign (0)=0$, $\{ t \}=t-[t]\ge 0$.
\par  {\bf 2.5. Lemma.} {\it Let a manifold $M$
and spaces $E^t_{\beta ,\delta }(TM)$
and $H^{ \{ l \} }_{\{ \gamma \} ,\delta ,\eta }(TM)$ 
be the same as in \S \S 2.2-2.4 with $l(k)\ge
[t]+[k/2]+3+sign \{ t \}$, $\gamma (k)\ge \beta +[k/2]+7/2+sign\{ t \}$.
Then there exist constants $C>0$ and $C_n>1$ for each $n$
such that $\| \zeta \| _{E^t_{\beta ,\delta }(TM)} \le 
C \| \zeta \| _{ \{ l \}, \{ \gamma \} ,\delta ,0}$ for each
$\zeta \in H^{ \{ l \} }_{ \{ \gamma \} ,\delta ,0}(TM)$, 
moreover, there can be chosen $\omega _n\ge C_n$, $C_{n+1}\ge k(n+1)
(k(n+1)-1)...(k(n)+1)C_n$
for each $n$ such that the following inequality be valid:
$\| \xi \| _{C^{l'(k)}_{\gamma '(k)}(TM_k)}$
$ \le C_n \| \xi \| _{H^{l(k)}
_{\gamma (k)}(TM_k)}$ for each $k=k(n)$, $l'(k)=l(k)-[k/2]-1$,
$\gamma '(k)=
\gamma (k)-[k/2]-1$ for each $\xi \in H^{l(k)}_{\gamma (k)}(TM_k)$.}
\par  {\bf Proof.} In view of theorems from ~\cite{tri} 
and the inequality $\int_{\bf R^m}
<x>_m^{-m-1}dx< \infty $ (for $<x>_m$ taken in $\bf R^m$ with
$x \in \bf R^m$) there exists the embedding $H^{l(n)}_{\gamma (n)}(TM_n)$
$\hookrightarrow C^{l'(n)}_{\gamma '(n)}(TM_n)$ for each $n$, 
since $2([n/2]+1)
\ge n+1$.  Moreover, due to results of \S III.6 \cite{miha}
there exists a constant $C_n>0$ for each $k=k(n)$, 
$n \in \bf N$ such that $\| \xi \| _{C^{l'(k)}_{\gamma '(k)}(TM_k)}$
$ \le C_n \| \xi \| _{H^{l(k)}
_{\gamma (k)}(TM_k)}$ for each $\xi \in H^{l(k)}_{\gamma (k)}(TM_k)$.
Then $D^{\alpha }f(x^1,...,x^n,...)-D^{\alpha }f(y^1,...,y^n,...)=$
$$\sum_{n=0}^{\infty }(D^{\alpha }f(y^1,...,y^{n-1},x^n,...)-D^{\alpha }
f(y^1,...,y^n,x^{n+1},...))$$ 
for each $f \in H^{ \{ l \} }_{ \{ \gamma \} ,\delta }(TM)$
in local coordinates, where $f(y^1,...,y^{n-1},x^n,...)=f(x^1,x^2,...,
x^n,...)$, if $n=0$; $\alpha =(\alpha ^1,...,\alpha ^m)$, $m \in \bf N$,
$\alpha ^i \in {\bf N_o}:=\{ 0,1,2,... \} $. 
Therefore, for each $x^n< y^n$ the following inequlity is satisfied: \\
$|D^{\alpha }f(y^1,...,y^{n-1},
x^n,x^{n+1},...)-D^{\alpha }f(y^1,...,y^n,x^{n+1},...)| _{l_{2,\delta }}
\bar m^{\alpha \delta }\le $
$$[\int_{\phi _j(U_j\cap M_k)\ni z:=
(y^1,...,y^{n-1},z^n), x^n\le z^n\le y^n}
|D^{\alpha }\partial f(y^1,...,y^{n-1},z^n,x^{n+1},...)/\partial z^n
|_{l_{2,\delta }}dz^n]\bar m^{\alpha \delta }\le $$
$$C^1\int \int_{\phi _j(U_j\cap M_{k(n+1)}) 
\ni z:=(y^1,...,y^{n-1},z^n, z^{n+1}), x^n\le z^n\le y^n}
sup_{x \in M}(\| f \| _{H^{l(k(n+8))}_{ \gamma (k(n+8))}(M_{k(n+8)}|TM)}$$
$<z>_{n+1}^{-5/2})dz^ndz^{n+1}(n+1)^{-2}
\le C' \| f \| _{ \{ l \}, \{ \gamma \} ,\delta ,0} \times (n+1)^{-2}$, \\
when $| \alpha |=\alpha ^1+...+\alpha ^m \le l(k)$,
$k=k(n)\ge n$, $m \le n$,
where $C^1=const >0$ and $C'=const >0$ are constants not depending on
$n$ and $k$; $x,$ $y$ and $(y^1,..,y^n,x^{n+1},x^{n+2},...)\in \phi _j(U_j)$
for each $n\in \bf N$. This is possible due to local convexity of the subset 
$\phi _j(U_j)\subset l_2$. Therefore, $H^{ \{ l \} }_{ \{
\gamma \} ,\delta ,0}(TM) \subset E^t_{\beta ,\delta }(TM)$ and $\| f \|
_{E^t_{\beta ,\delta }(TM)}\le C \| f \| _{ \{ l \}, \{ \gamma \} ,
\delta ,0}$ for each
$f \in H^{ \{ l \} }_{ \{ \gamma \} ,\delta ,0}(TM)$, moreover, $C=C'
\sum_{n=1}^{\infty }
n^{-2} <\infty $, since $\sup _{x \in M}\sum_{j=1}^{\infty }g_j(x)
\le \sum_{j=1}^{\infty }\sup_{x \in M}g_j(x)$ for each function 
$g: M \to [0, \infty )$ and
$\lim_{R\to \infty }\| f|_{M_R^c} \| _{E^t_{\beta ,\delta }(TM)}$
$\le C\times \lim_{R\to
\infty }\| $ $f|_{M_R^c}\|_{\{ l \}, \{ \gamma \} ,\delta ,0}=0$.
\par The uniform space 
$E^t_{\beta ,\delta }(TM)\cap H^{ \{ l \} }_{ \{ \gamma \} ,\delta,0}
(TM)$ contains the corresponding cylindrical functions
$\zeta $, in particular with
$supp(\zeta ) \subset U_j\cap M_n$ for some $j \in \bf N$ and
$k=k(n)$, $n \in \bf N$. The linear span of the family $\sf K$
over the field $\bf R$ of all such functions $\zeta $ is dense in 
$E^t_{\beta ,\delta }(TM)$ and in $H^{ \{ l \} }_{ \{ \gamma \} ,\delta ,0}
(TM)$ due to the Stone-Weierstrass theorem,
consequently, $H^{ \{ l \} }_{ \{ \gamma \},\delta ,0}(TM)$ 
is dense in $E^t_{\beta ,\delta }(TM)$,
since $\partial f/\partial x^{n+1}=0$
for cylindrical functions $f$ independent from $x^{n+j}$ for each $j>0$.
\par {\bf 2.5.1. Lemma.} {\it Let $M$ and $N$ be both either
complex or quaternion manifolds and $M_{\bf R}$ and $N_{\bf R}$
be underlying their Riemann manifolds over $\bf R$.
Then $H:=H^{ \{ l \} }_{ \{ \gamma \} , \delta , \eta }(M,N)$
has a holomorphic embedding into $H^{\bf R}:=
H^{ \{ l \} }_{ \{ \gamma \} , \delta , \eta }(M_{\bf R},N_{\bf R})$
as the closed uniform subspace and the codimension of $H$ in $H^{\bf R}$
over $\bf R$ is infinite.}
\par {\bf Proof.} Elements of $H^{\bf R}$ are Fr\'echet differentiable,
consequently, operators $f\mapsto Df$ are continuous from
$H^{\bf R}$ into $D(H^{\bf R})$. Since $D=\partial + {\bar \partial }$
in the complex case and $D=\partial + {\tilde \partial }$
in the quaternion case, then the generalized holomorphicity condition
${\bar \partial }f=0$
in the complex case or ${\tilde \partial f}=0$ in the quaternion case
defines the closed uniform subspace $H$ in $H^{\bf R}$.
This is correct, since transition mappings of charts are either
complex or quaternion holomorphic respectively.
The uniform subspace of mappings $Df\ne \partial f$ is infinite dimensional,
since subspaces of locally analytic functions having nonvanishing
series in variables $\bar z$ or $\tilde z$ (in local coordinates) are
infinite dimensional. Thus $H$ has an infinite codimension in $H^{\bf R}$
over $\bf R$.
\par {\bf 2.6.1. Note.}
For a diffeomorphism group $Diff^t_{\beta ,\gamma }(\tilde M)$
of a Banach manifold $\tilde M$ let $M$ be a dense Hilbert submanifold
in $\tilde M$ as in \cite{lurim1,lurim2}.
\par  {\bf 2.6.2. Lemma.} {\it Let 
$Di^{ \{ l \} }_{ \{ \gamma \} , \delta , \eta }(M)$ and $M$ 
be the same as in \S 2.3 with values of parameters $C_n$ from Lemma 2.5 
for given $l(k)$, $\gamma (k)$ and $k=k(n)$ with $\omega _n=l(k(n))!C_n$, 
then $Di^{ \{ l \} }_{ \{ \gamma \} , \delta ,\eta }(M)$ is the 
separable metrizable topological group dense in
$Diff^t_{\beta ,\delta }(\tilde M)$.
\par In the case of a complex or quaternion manifold $M$ 
the groups $Diff^{\xi }(M)$, $Dif^t_{\beta , \gamma }(M)$
and $Di^{ \{ l \} }_{ \{ \gamma \} , \delta , \eta }(M)$
(see \S \S 2.1.5, 2.2, 2.3) are the separable metrizable topological
groups. They have embeddings as closed subgroups into
the groups $Diff^{\xi }(M_{\bf R})$, $Dif^t_{\beta , \gamma }(M_{\bf R})$
and $Di^{ \{ l \} }_{ \{ \gamma \} , \delta , \eta }(M_{\bf R})$
respectively of infinite codimensions.}
\par  {\bf Proof.} Consider at first the real case.
From the results of the paper \cite{omo}
it follows that the uniform space
$Di^{ \{ l \} }_{ \{ \gamma \} , \delta ,\eta }(M_k)$ 
is the topological group for each finite dimensional submanifold
$M_k$, since $l(k)>k+5$ and $dim_{\bf R}M_k=k$. The minimal algebraic group
$G_0:=gr(Q)$ generated by the family $Q:=\{ f: f\in E^{ \{ l\} ,id}_{
\{ \gamma \} ,\delta }(U,V)$ $\mbox{for all possible pairs of charts }
U_i$ $\mbox{and}$ $U_j$ $\mbox{with}$ $U=\phi _i(U_i)$ $\mbox{and}$
$V=\phi _j(U_j),$ $supp(f)\subset U\cap \bf R^n,$ 
$f\in Hom(M),$ $dim_{\bf R}M\ge n \in {\bf N} \}$
is dense in $Di^{ \{ l \} }_{ \{ \gamma \} , \delta ,\eta }(M)$ 
and in $Diff^t_{\beta ,\delta }(M)$ due to the Stone-Weierstrass theorem, 
since the union $\bigcup_kM_k$ is dense in $M$, where
$supp(f):=cl\{ x\in M: f(x)\ne x \}$, $cl(B)$ denotes the closure
of a subset $B$ in $M$. 
Therefore, $Di^{ \{ l \} }_{ \{ \gamma \} , \delta ,\eta }(M)$
and $Diff^t_{\beta ,\delta }(\tilde M)$ are the separable topological
spaces. It remains to verify that
$Di^{ \{ l \} }_{ \{ \gamma \} , \delta ,\eta }(M)$
is the topological group. For it can be used 
Lemma 2.5. For $a>0$ and $k\ge 1$ using integration by parts 
formula we get the following equality $\int^{\infty }_{-\infty }(
a^2+x^2)^{-(k+2)/2}dx=((k-1)/(ka^2))\int^{\infty }_{-\infty }
(a^2+x^2)^{k/2}dx$, which takes into account the weight multipliers.
Let $f, g\in Di^{ \{ l\} }_{ \{ \gamma \}, \delta ,\eta }(V)$
for an open subset $V=\phi _j(U_j)\subset l_2$ and
$\chi _{ \{ l\}, \{ \gamma \},\delta ,\eta }(f,id)<1/2$
and $\chi _{ \{ l\}, \{ \gamma \},\delta ,\eta }(g,id)<\infty $, then
$\bar \rho ^l_{\gamma ,n,\delta }(g^{-1}\circ f,
id) $ $\le C_{l,n,\gamma ,\delta }(\bar \rho ^{4l}_{\gamma ,n,\delta }
(f,id)+\bar \rho ^{4l}_{\gamma ,n,\delta }(g,id))$,
where $0<C_{l,n,\gamma ,\delta }\le 1$ is a constant dependent on 
$l, n, \gamma $ and independent from $f$ and $g$.
For the Bell polynomials $Y_n$ there is the following inequality
$Y_n(1,...,1) \le n!e^n$ for each $n$ and $Y_n(F/2,...,F/(n+1))
\le (2n)!e^n$ for $F^p:=F_p=(n+p)_p:=(n+p)...(n+2)(n+1)$ (see
Chapter 5 in ~\cite{rio} and Theorem 2.5 in ~\cite{ave}). 
The Bell polynomials are given by the following formula
$Y_n(fg_1,...,fg_n):=\sum_{\pi (n)}(n!f_k/(k_1!...k_n!))(g_1/1!)^{k_1}
...(g_n/n!)^{k_n}$, where the sum is by all partitions
$\pi (n)$ of the number $n$, this partition is denoted by $1^{k_1}2^{k_2}...
n^{k_n}$ such that $k_1+2k_2+...+nk_n=n$ and $k_i$ is a number of terms
equal to $i$, the total number of terms in the partition is equal to
$k=k(\pi )=k_1+....+k_n$, $f^k:=f_k$ in the Blissar calculus notation.
For each $n\in \bf N$, $l=l(n)$ and $\gamma =\gamma (n)$
the following inequality is satisfied:
$\bar \rho ^l_{\gamma ,n,\delta }(f\circ g,id)
\le Y_l(\bar f\bar g_1,...,\bar f\bar g_m)$, 
$\bar \rho ^{l+1}_{\gamma ,n,\delta }(f^{-1},id)\le (3/2)Y_l(Fp_1/2,
...,Fp_m/(m+1))$, where $\bar f^m:=\bar f_m=\bar \rho ^m_{\gamma ,n,
\delta }(f,id)$ and $F^k:=F_k=(n+k)_k$, $(n)_j:=n(n-1)...(n-j+1)$,
$p_k=-\bar f_{k+1}(3/2)^{k+1}$.
Then $\sum_{n=1}^{\infty }(2l(k(n)))!e^{l(k(n))}b^{4l(k(n))}
(l(k(n))!)[(4l(k(n)))!]^{-1}<\infty $ for each $0<b<\infty $.
Hence due the Cauchy-Schwarz-Bunyakovskii inequality 
and the condition $C_n^2>C_n$
for each $n$ we get: $f\circ g$ and $f^{-1}\in Di^{ \{ l\} }_{ \{ \gamma \},
\delta ,\eta }(M)$ for each $f$ and $g\in 
Di^{ \{ l\} }_{ \{ \gamma \},\delta ,\eta }(M)$, moreover, the operations 
of composition an inversion are continuous.
\par The base of neighborhoods of $id$ in
$Di^{ \{ l \} }_{ \{ \gamma \} , \delta ,\eta }(M)$
is countable, hence this group is metrizable, moreover, 
a metric can be chosen left-invariant due to Theorem 8.3 \cite{hew}.
The case $Diff^{\xi }(M)$ for a complex or quaternion
manifold $M$ is analogous to that of considered above.
The latter statement follows from Lemma 2.5.1.
\par {\bf 2.7. Lemma.} {\it Let 
$G':=Di^{ \{ l" \} }_{ \{ \gamma " \} ,\delta ",\eta "}(M)$ 
be a subgroup of $G:=Di^{ \{ l \} }_{ \{ \gamma \} ,\delta ,\eta }(M)$ 
such that $m(n)>n/2$, $l"(n)=l(n)+m(n)n$,
$\gamma "(n)=\max (\gamma (n)-m(n)n,0)$
for each $n$, $inf-\lim_{n\to \infty }m(n)/n=c>1/2$, 
$\delta ">\delta +1/2$, $\infty >\eta ">\eta +1/2,$ $\eta \ge 0$ 
(see \S 2.3).  Let also $G':=Diff^{\xi '}(M')$ be a subgroup of
$G=Diff^{\xi }(M)$ with either ${a'}_1<a_1$ 
and ${c'}_1<c_1$ or ${a'}_1=a_1$ and ${a'}_2<a_2-1$ 
and ${c'}_1=c_1$ and ${c'}_2<c_2-1$ for the complex or quaternion
manifold $M$ (see \S \S 2.1.5 and 2.2).
Then there exists a Hilbert-Schmidt operator of
embedding $J : Y'\hookrightarrow Y$, where
$Y:=T_eG$ and $Y':=T_eG'$ are tangent Hilbert spaces
over $\bf R$, $\bf C$ or $\bf H$ respectively.}
\par  {\bf Proof.} Consider at first the real case.
The natural embedding 
$\theta _k$ of the Hilbert spaces $H^{l(k)-m(k)k,b(k)}_{
\gamma (k)+m(k)k,\delta }(M_k,{\bf R})$
into $H^{l(k),b(k)}_{\gamma (k)k,\delta }(M_k,l_{2,\delta +1+\epsilon })$
is the Hilbert-Schmidt operator for each $k=k(n)$, $n\in \bf N$
(see their definition for general Banach spaces in \cite{pietsch}).
For each chart $(U_j,\phi _j)$ there are 
linearly independent functions
$x^me_l<x>^{\zeta }_n/m!=:f_{m,l,n}(x)$, where $\{ e_l: l \in {\bf N} \}
\subset l_2$ is the standard orthonormal basis in $l_2$, 
$x^m:=x_1^{m_1}...x_n^{m_n}$,
$m!=m_1!...m_n!$, $<x>_n=$ $(1+\sum_{i=1}^n(x^i)^2)^{1/2}$, $n \in \bf N$,
$\zeta (n)=\zeta \in \bf R$.
The linear span over $\bf R$ of the family of all such functions
$f(x)$ is dense in $Y$.
Moreover, $D^{\alpha }f(x)=e_l\sum {\alpha \choose \beta }(D^{\beta }x^m/m!)$
$(D^{\alpha -\beta }<x>^{\zeta }_n)$, where 
$D^{\alpha }=\partial _1^{\alpha ^1}...\partial _n^{\alpha ^n},$
$\partial _i=\partial /\partial x_i$,
$\alpha =(\alpha ^1,...,
\alpha ^n)$, ${\alpha \choose \beta }$
$=\prod _{i=1}^n{{\alpha ^i} \choose {\beta ^i}}$,
$0\le \alpha ^i\in \bf Z$, $\lim_{n \to \infty }
q^n/n!=0$ for each $\infty >q>0$, $\sum_{j,l,n=1}^{\infty }
\sum_{|m|
\ge m(n), m}[jln^nm_1...m_n]^{-(1+2\epsilon )}<\infty $
for each $0<\epsilon <\min (c-1/2, \eta "- \eta -1/2, \delta "-\delta -1/2)$, 
where $m=(m_1,...,m_n)$, $|m|:=m_1+...+m_n$, $0\le m_i \in \bf Z$.
Hence due to \S \S 2.3 and 2.4 the embedding $J$ is the Hilbert-Schmidt 
operator.
\par In the complex case and quaternion cases it is possible
to use the convergence of the series \\
$\sum_{j=1}^{\infty }\sum_{n=1}^{\infty}(j!)^{{a'}_1-a_1}
(n!)^{{c'}_1-c_1}<\infty $ and 
$\sum_{j=1}^{\infty }\sum_{n=1}^{\infty}j^{{a'}_2-a_2}
n^{{c'}_2-c_2}<\infty $.
\par  {\bf 2.8. Theorems.} {\it Let diffeomorphism groups
$G:=Di^{ \{ l \} }_{ \{ \gamma \} , \delta ,\eta }(M)$
and $G:=Diff^{\xi }(M)$ be the same as in  \S \S 2.2, 2.3. Then
$$(i) \mbox{for each }H^{ \{l \} ,id }_{ \{ \gamma \} ,\delta ,\eta }
(M,TM)\mbox{-vector field } V \mbox{ its flow }\eta _t$$ 
$\mbox{ is a one-parameter subgroup of }
Di^{ \{ l \} }_{ \{ \gamma \} , \delta ,\eta }(M)$,
$\mbox{ the curve }t\mapsto \eta _t\mbox{ is of class }C^1,$ 
$\mbox{ the mapping}$
$\tilde Exp: \tilde T_e Di^{ \{ l \} }_{ \{ \gamma \} , \delta ,\eta }(M)\to
Di^{ \{ l \} }_{ \{ \gamma \} , \delta ,\eta }(M)$,
$\mbox{ is continuous }$ and defined on the neighbourhood 
$\tilde T_e Di^{ \{ l \} }_{ \{ \gamma \} , \delta ,\eta }(M)$
of the zero section in
$T_eDi^{ \{ l \} }_{ \{ \gamma \} , \delta ,\eta }(M)$,
$V\mapsto \eta _1$;
$$(ii)\mbox{ }T_fDi^{ \{ l \} }_{ \{ \gamma \} , \delta ,\eta }(M),
=\{ V \in H^{ \{l \} ,id }_{ \{ \gamma \} ,\delta ,\eta }
(M,TM)| \pi \circ V=f \};$$
$$(iii)\mbox{ }(V,W)=\int_Mg_{f(x)}(V_x,W_x)\mu (dx)$$
is a weak Riemannian structure on a Hilbert 
manifold $Di^{ \{ l \} }_{ \{ \gamma \} , \delta ,\eta }(M)$, where
$\mu $ is a measure induced on $M$ by $\phi _j$ and a Gaussian measure
with zero mean value on $l_2$ produced by an injective self-adjoint
operator $Q: l_2 \to l_2$ of trace class, $0< \mu (M) <\infty $;
$$(iv)\mbox{ the Levi-Civita connection }\nabla \mbox{ on } M
\mbox{ induces the Levi-Civita connection }$$ $\hat \nabla $ on 
$Di^{ \{ l \} }_{ \{ \gamma \} , \delta ,\eta }(M)$;
$$(v)\mbox{ } \tilde E: T Di^{ \{ l \} }_{ \{ \gamma \} , \delta ,\eta }(M)
\to Di^{ \{ l \} }_{ \{ \gamma \} , \delta ,\eta }(M)
\mbox{ is
defined by }$$ $\tilde E_{\eta }(V)=exp_{\eta (x)} \circ V_{\eta }$
on a neighbourhood $\bar V$ of the zero section
in $T_{\eta }Di^{ \{ l \} }_{ \{ \gamma \} , \delta ,\eta }(M)$
and is a $H^{ \{l \} ,id }_{ \{ \gamma \} ,\delta ,\eta }$-mapping 
by $V$ onto a neighbourhood $W_{\eta }=W_{id}\circ \eta $ of $\eta
\in Di^{ \{ l \} }_{ \{ \gamma \} , \delta ,\eta }(M)$; 
$\tilde E$ is the uniform isomorphism of uniform spaces $\bar V$ and $W$.
Analogous statements are true for $Diff^{\xi }(M)$
with the class of smoothness $Y^{\xi ,id}$
instead of $H^{ \{l \} ,id }_{ \{ \gamma \} ,\delta ,\eta }$.}
\par   {\bf Proof.} We have that $T_f
H^{ \{l \} ,\theta }_{ \{ \gamma \} ,\delta ,\eta }(M,N)
= \{ g \in H^{ \{l \} ,\theta }_{ \{ \gamma \} ,\delta ,\eta }(M,TN):$
$ \pi _N\circ g=f \} $, where $\pi _N: TN\to N$ is the canonical projection.
Therefore, the tangent space is given by the formula
$TH^{ \{l \} ,\theta }_{ \{ \gamma \} ,\delta ,\eta }(M,N)
=H^{ \{l \} ,\theta }_{ \{ \gamma \} ,\delta ,\eta }(M,TN)
=\bigcup_fT_fH^{ \{l \} ,\theta }_{ \{ \gamma \} ,\delta ,\eta }(M,N)$ 
and the following mapping $w_{exp}:T_f
H^{ \{l \} ,\theta }_{ \{ \gamma \} ,\delta ,\eta }(M,N) \to
H^{ \{l \} ,\theta }_{ \{ \gamma \} ,\delta ,\eta }(M,N)$, 
$w_{exp}(g)=exp \circ g$ gives
charts for $H^{ \{l \} ,\theta }_{ \{ \gamma \} ,\delta ,\eta }(M,N)$, 
since  $TN$ has an atlas of class 
$H^{ \{ l'(n)-1: n \} }_{ \{ \gamma ' (n)+1 : n \} , \chi }.$
Apply Theorem 5 about differential equations
on Banach manifolds in \S 4.2 \cite{lan} to the case considered here.
Then a vector field $V$ of class $H^{ \{l \} ,\theta }_{ \{ \gamma \} ,
\delta ,\eta }$ on $M$ defines a flow $\eta _t$
of such class, that is $d \eta _t/dt=V \circ \eta _t$ and
$\eta _0=e$. From the proofs of Theorem 3.1 and Lemmas 3.2,
3.3 in \cite{ebi} we get that $\eta _t$ is a one-parameter subgroup of
$Di^{ \{ l \} }_{ \{ \gamma \} , \delta ,\eta }(M)$, 
the curve $t \mapsto \eta _t$ is of class $C^1$, 
the map $\tilde Exp: T_eDi^{ \{ l \} }_{ \{ \gamma \} , \delta ,\eta }(M)
\to Di^{ \{ l \} }_{ \{ \gamma \} , \delta ,\eta }(M)$ defined by $V \mapsto
\eta _1$ is continuous.
\par Each curve of the form $t \mapsto \tilde E(tV)$ is a geodesic for $V
\in T_{\eta }Di^{ \{ l \} }_{ \{ \gamma \} , \delta ,\eta }(M)$ 
such that $d \tilde E(tV)/dt$ is the map
$m \mapsto d(exp(tV(m))/dt=\gamma '_m(t)$ for each $m\in M$, 
where $\gamma _m(t)$ is a
geodesic on $M$, $\gamma _m(0)=\eta (m)$, $\gamma _m'(0)=V(m)$.
Indeed, this follows from an existence of a solution of a corresponding
differential equation in the Hilbert space 
$H^{ \{l \} ,\eta }_{ \{ \gamma \} ,\delta ,\eta }(M|TM)$,
then we proceed as in the proof of Theorem 9.1 \cite{ebi}.
\par From the definition of $\mu $ it follows that for each $x \in M$
there exists an open neighbourhood $Y \ni x$ such that $\mu (Y)>0$
\cite{sko}. 
Since $t\ge 1$, the scalar product $(iii)$ gives a weaker topology
than the initial $H^{ \{l \} }_{ \{ \gamma \} ,\delta ,\eta }$.
\par Then the right multiplication $\alpha _h(f)=f \circ h$, 
$f \to f \circ h$ is of class $C^{\infty }$ on 
$Di^{ \{ l \} }_{ \{ \gamma \} , \delta ,\eta }(M)$ 
for each $h \in Di^{ \{ l \} }_{ \{ \gamma \} , \delta ,\eta }(M)$. 
Moreover, $Di^{ \{ l \} }_{ \{ \gamma \} , \delta ,\eta }(M)$
acts on itself freely from the right, hence we have
the following principal vector bundle $\tilde \pi :
T Di^{ \{ l \} }_{ \{ \gamma \} , \delta ,\eta }(M) \to 
Di^{ \{ l \} }_{ \{ \gamma \} , \delta ,\eta }(M)$
with the canonical projection $\tilde \pi $.
\par Analogously to \cite{ebi,lurim1} we get the connection $\hat \nabla
=\nabla \circ h$ over $\bf R$ on 
$Di^{ \{ l \} }_{ \{ \gamma \} , \delta ,\eta }(M)$. 
If $\nabla $ is torsion-free then $\hat \nabla $
is also torsion-free. From this it follows that the existence of $\tilde E$
and $Di^{ \{ l \} }_{ \{ \gamma \} , \delta ,\eta }(M)$
is the Hilbert manifold of class
$H^{ \{ l'(n)-1: n \} }_{ \{ \gamma '(n)+1: n \} ,\chi ,\eta }$, 
since $exp$ for $M$ is of class
$H^{ \{ l'(n)-1: n \} }_{ \{ \gamma '(n)+1: n \} ,\delta }$,
$f\to f\circ h$ is a $C^{\infty }$ mapping with the derivative $\alpha _h: 
H^{ \{l \} ,\eta }_{ \{ \gamma \} ,\delta ,\eta }(M',TN)
\to H^{ \{l \} ,\eta }_{ \{ \gamma \} ,\delta ,\eta }(M,TN)$ 
whilst $h \in H^{ \{l \} ,\eta }_{ \{ \gamma \} ,\delta ,\eta }(M,M'),$ 
\par $(vi)$ $\tilde E_h(\hat V):=exp_{h(x)}(V(h(x)))$, where 
\par $(vii)$ $\hat V_h=V\circ h$, 
$V$ is a vector field in $M$, $\hat V$ is a vector field 
in $Di^{ \{ l \} }_{ \{ \gamma \} , \delta ,\eta }(M).$
\par {\bf 2.9. Theorem.} {\it Let $M$ and $N$ be manifolds
both either real or complex or quaternion and $(L^MN)_{\xi }$ and
$Diff^{\xi }_{y_0}(N)$ be a group of loops and a group of diffeomorphisms
preserving a marked point $y_0\in N$, where $\xi $ is a class of smoothness
(see \S \S 2.1, 2.2 and 2.3). Then there exists a topological locally
connected nonlocally compact group which is their semidirect product
$Diff^{\xi }_{y_0}(N) \otimes ^s(L^MN)_{\xi }$ such that
for complex or quaternion manifolds it has an embedding as
the closed subgroup into $Diff^{\xi }_{y_0}(N_{\bf R})
\otimes ^s(L^{M_{\bf R}}N_{\bf R})_{\xi }$ with infinite codimension.
If $M$ and $N$ are complex manifolds, then there exist quaternion
manifolds $P$ and $Q$ and complex holomorphic embeddings of
$M$ and $N$ into $P$ and $Q$ respectively and an embedding
of $Diff^{\xi }_{y_0}(N) \otimes ^s(L^MN)_{\xi }$ into
$Diff^{\xi }_{y_0}(Q) \otimes ^s(L^PQ)_{\xi }$ as a closed
subgroup of infinite codimension of
$T_e (Diff^{\xi }_{y_0}(N) \otimes ^s(L^MN)_{\xi })$ in
$T_e (Diff^{\xi }_{y_0}(Q) \otimes ^s(L^PQ)_{\xi })$.}
\par {\bf Proof.} Let $Diff^{\xi }_{y_0}(N)$ be the subgroup of the
group of diffeomorphisms preservind a marked point $y_0$ in $N$,
$\zeta \in Diff^{\xi }_{y_0}(N)$ be a marked element in this subgroup,
then it induces the internal automoprhism $Diff^{\xi }_{y_0}(N)\ni
\psi \mapsto \zeta \circ \psi \circ \zeta ^{-1}$.
For each function $f: M\to N$ with $f(s_0)=y_0$ and each diffeomorphism
$\psi \in Diff^{\xi }_{y_0}(N)$ there is defined a mapping
$\psi (f): M\to N$ such that $\psi (f)(s_0)=y_0$.
Consider the equivalence relation caused by the action of
$Diff^{\xi }_{s_0}(M)$ on the space of all such mappings
$f$ of the class of smootness corresponding to $\xi $.
From Theorem 2.1.4.1 it follows, that $\psi $
is an automorphism of the loop monoid $(S^MN)_{\xi }$,
since $\psi (\omega _0)=\omega _0$ and
$\psi (f_1\vee f_2)=\psi (f_1)\vee \psi (f_2)$
for each $f_1$ and $f_2: M\to N$ with $f_1(s_0)=y_0$ and $f_2(s_0)=y_0$,
where $\omega _0(M):=\{ y_0 \} $. In view of Theorem 2.1.7(1) the
diffeomorphism $\psi $ induces the automorphism
of the loop group, which it is convenient to denote by
$\psi : g\ni (L^MN)_{\xi }\mapsto \psi (g)$.
Thus there exists a semidirect product of these groups
for which products of elements $(\psi _1,g_1)$ and $(\psi _2,g_2)
\in Diff^{\xi }_{y_0}(N) \otimes ^s(L^MN)_{\xi }$ are given by the formula:
$(\psi _1,g_1)(\psi _2,g_2)=(\psi _1\circ \psi _2,g_1g_2^{\psi _1})$,
where $g^{\psi }:=\zeta \circ \psi \circ \zeta ^{-1} (g)$
for each $g\in (L^MN)_{\xi }$.
Since the diffeomorphism group $Diff^{\xi }_{y_0}(N)$ is nonlocally
compact and locally connected, then such is also
$Diff^{\xi }_{y_0}(N) \otimes ^s(L^MN)_{\xi }$.
The latter statements about embeddings follows from Theorems
2.1.7(4), 2.2.1 and Lemma 2.6.2.
\par {\bf 2.10. Theorem.} {\it Let
$G:=Diff^{\xi }_{y_0}(N) \otimes ^s(L^MN)_{\xi }$ be a semidirect
product of a group of diffeomorphisms and a group of loops
as in \S 2.9. Then $G$ is an infinite dimensional uniformly complete
Lie group which does not satisfy even locally the Campbell-Hausdorff formula
as well as its closed subgroups $Diff^{\xi }_{y_0}(N)$ and $(L^MN)_{\xi }$,
besides the degenerate case of $Diff^{\sf O}(N)$ for a compact
complex manifold $N$. If both manifolds $M$ and $N$ are either complex or
quaternion, then $G$ has a structure of a complex or a quaternion
manifold respectively.}
\par {\bf Proof.} If $dim_{\bf R}N>1$, then $(L^MN)_{\xi }$
is infinite dimensional and nonlocally compact manifold
(see Theorem 2.1.7). A group $Diff^{\xi }_{y_0}(N)$ may be locally
compact only for finite dimensional complex manifold $N$ and $\xi =\sf O$,
but then $dim_{\bf R}N>1$, hence $G$ is nonlocally compact in all cases.
Groups $Diff^{\xi }_{y_0}(N)$ and $(L^MN)_{\xi }$ have structures
of smooth manifolds, hence so is $G$ also, since $Diff^{\xi }_{y_0}(N)$
acts smoothly on $(L^MN)_{\xi }$. Group operations $(f,g)\mapsto
fg^{-1}$ are smooth in them, hence they are Lie groups.
For both either complex or quaternion manifolds $M$ and $N$
groups of loops and diffeomorphisms have structures of complex
and quaternion manifolds respectively, hence such is their
semidirect product also.
It is known, that the diffeomorphism group does not satisfy
the Campbell-Hausdorff formula besides the case of a compact complex
manifold $N$ for $Diff^{\sf O}(N)$ (see \cite{bomon,ebi,kobtg}). 
The uniformities in $Diff^{\xi }_{y_0}(N)$ and $(L^MN)_{\xi }$
correspond to the class of smoothness characterized by $\xi $.
Thus these subgroups are closed in their semidirect product.
\par The loop group $(L^MN)_{\xi }$ for $dim_{\bf R} N > 1$
is nondiscrete. To prove that it does not satisfy locally
the Campbell-Hausdorff formula it is sufficient to prove it for its
subgroup $(L^{M_m}N)_{\xi }$ with compact submanifold $M_m$.
Suppose that it has a nontrivial local one-parameter
subgroup $ \{ g^b:$ $b\in (-a,a) \} $ with $a>1$ for an element $g$
corresponding to an equivalence class of a mapping $f: M_m\to N$,
$f (s_0) = y_0$, when $f$ is such that
$\sup_{y\in N} [card (f^{-1}(y))] = k <\aleph _0$.
The condition $a>1$ is not restrictive, since it is possible
to consider arbitrary small neighbourhood of $e$.
The existence of a nontrivial local one-parameter subgroup
imply that for each integer $0\ne p\in {\bf Z}$ there exists $g^{1/p}$
in $(L^M_mN)_{\xi }$ such that $(g^{1/p})^p=g$.
This number $p$ may be arbitrary large. The mapping $f_p$ belonging to
the class of equivalence corresponding to $g^{1/p}$ has at least one
restriction $f_p (s_0) = y_0$, but then wedge product of $f_p$
with itself $|p|$ times and the corresponding equivalence class
would give an element $h_p$ in it having at least one point $y\in N$
with $|p|\le card (h_p^{-1}(y))$. On the other hand, $h_p(M)=f(M)$.
Then $g^{k/p}$ with relatively prime $k$ and $p$, $(k,p)=1$,
with $p$ arbitrary large would give in the symmetric neighbourhood $U=U^{-1}
\ni g \ne e$, $e\notin U$, elements $g^{k/p}$ in $U$ with mapping
representatives $h_{k/p}$ of these classes having arbitrary large amount
of distinct points in $M_m$ belonging to $h_{k/p}^{-1}(y)$ for each
$y\in f(M)\subset N$ while $p$ tends to the infinity. But $h_1=f$ has
$\sup_{y\in N} [card (h_1^{-1}(y))] = k <\aleph _0$ by the supposition above.
The diameter of $M_m$ as the metric space is positive together with
$M_m^{\vee p}$.
Since $g\ne e$, then $<h_{k/p}>_{\xi }=g^{k/p}$ can not converge
to $g$ while $k/p\in \bf Q$ tends to $1$. Thus for each neighbourhood
$V$ of $e$ in $(L^MN)_{\xi }$ there exists $e\ne g\in V$ such that
$g$ does not belong to any local one-parameter subgroup.
In $G:=Diff^{\xi }_{y_0}(N) \otimes ^s(L^MN)_{\xi }$ an element
$(e,g)$ does not belong to any local one-parameter subgroup.
Moreover, $G$ does not satisfy the Campbell-Hausorff formula.
\par Let now both manifolds $M$ and $N$ be either complex or quaternion,
then the tangent bundle $TC^1(M,N)$ is isomorphic with $C^1(M,TN)$,
since $M$ and $N$ are $C^{\infty }$-manifolds particularly.
Considering continuous piecewise complex or quaternion holomorphic
mappings $f$ from $M$ into $N$ and the completion of their family
by $Y^{\xi }$-uniformity, we get, that
$TY^{\xi }(M,N)$ is isomorphic with $Y^{\xi }(M,TN)$.
If $f: M\to N$ then $\psi _j\circ f\circ \psi _k^{-1}$ has domain
in the complex or quaternion vector space $T_yN$, where $(V_j, \psi _j)$
is a chart of $At (N)$, $y\in N$, $V_j\cap V_k \ne \emptyset $.
Consider charts $W_k (U) := \{ f\in Y^{\xi }(M,N):$  $f(U)\subset V_k \} $,
where $U$ is open in $M$, then the transition mapping between $W_k(U)$
and $W_j(U)$ is $\psi _k\circ \psi _j^{-1}$, since from $f\in W_j(U)$
it follows, that $(\psi _j\circ f) (U)\subset \psi _j(V_j)\subset T_yN$.
Since $\psi _k\circ \psi _j^{-1}$ is either complex or quaternion
holomorphic, then $Y^{\xi }(M,N)$ is complex or quaternion manifold
respectively. Since $T_yN$ is either complex or quaternion vector space,
then $T_fY^{\xi }(M,N)$ is either complex or quaternion vector space
respectively for each $f\in Y^{\xi }(M,N)$. Using constructions above
of $G$ from $Y^{\xi }(M,N)$ and Lemma 2.1.6.2,
we get that $G$ has the structure of
either complex or quaternion manifold correspondingly.
\par {\bf 2.11. Theorem.} {\it A loop group $G:=(L^MN)_{\xi }$ from
\S \S 2.1.5, 2.3 and diffeomorphism groups $G:=Diff^t_{\beta ,\gamma }(M)$
from \cite{lurim1} and $G:=Di^{ \{ l \} }_{ \{ \gamma \} ,\delta ,\eta }(M)$ 
from \S 2.3 and $G:=Diff^{\xi }(M)$ from \S 2.2 and their semidirect product
$G := Diff^{\xi }_{y_0}(N) \otimes ^s(L^MN)_{\xi }$ have uniform atlases.}
\par {\bf Proof.} In view of Theorems 3.1 and 3.3 \cite{lurim1}
and Theorems 2.1.7, 2.2.1, 2.8 above a diffeomorphism group 
$G$ and a loop group $G:=(L^MN)_{\xi }$ have uniform atlases (see \S 2.1)
consistent with their topology, where $M$ is a real
manifold $1\le t<\infty $, $0\le \beta <\infty $, $0\le \gamma \le \infty $
for a diffeomorphism group $Diff^t_{\beta ,\gamma }(M)$
(see \cite{lurim1}). Others parameters are specified in the cited paragraphs. 
They also include the particular cases of finite dimensional 
manifolds $M$ and $N$.
\par The case of complex compact $M$ for $G:=Diff^{\infty }(M)$ is trivial,
since $Diff^{\infty }(M)$ is the finite dimensional Lie group 
for such $M$ \cite{kobtg}. 
\par In view of Theorems 2.1.7, 2.8, 2.10 and Formulas 2.8.$(vi,vii)$
above and Theorem 3.3 \cite{lurim1} to satisfy conditions $(U1,U2)$
of \S 2.1.1 it is sufficient to find an atlas $At(G)$
of each such group $G$, for which $U_1$ is a neighbourhood of $e$,
$U_x^1$ and $U_x^2$ are for $x=e$ such that $\phi _1(U_1)$
contains a ball of radius $r>0$. Due to an existence of  
left-invariant metrics in each such topological group
and its paracompactness and separability 
we can take a locally finite covering
$\{ U_j: g_j^{-1}U_j\subset U_1 : j\in {\bf N} \} $, where
$\{ g_j: j\in {\bf N } \} $ is a countable subset of pairwise 
distinct elements of the group, $g_1=e$. Using
uniform continuity of $\tilde E$ we can satisfy $(U1,U2)$ 
with $r>0$, since the manifolds $M$ for diffeomorphism groups
and $N$ for loop groups also have uniform
atlases. Choosing $U_1$ in addition such that $\tilde E$ is bounded
on $U_1U_1$ and using left shifts $L_hg:=hg$, where $h$ and $g\in G$,
$AB:=\{ c: c=ab, a\in A, b\in B \} $ for $A\cup B\subset G$, 
and Condition $(U3)$ for $M$ and $N$ we get, that  
there exist sufficiently small neighbourhoods $U_1$, 
$U_e^1$ and $U_e^2$ with $U_e^2U_e^2\subset U_e^1$
and $U_x^1\subset xU_e^1$, $U_x^2\subset xU_e^2$ 
for each $x\in G$ such that Conditions $(U1-U3)$ are fulfilled,
since uniform atlases exist on the Banach or Hilbert 
tangent space $T_eG$.
\section{Differentiable transition probabilities on groups.}
\par {\bf 3.1. Definitions and Notes.}  Let $G$ be a
Hausdorff topological group, we denote by 
$\mu : Af(G,\mu )\to [0,\infty )\subset \bf R$ 
a $\sigma $-additive measure. Its left shifts 
$\mu _{\phi }(E):=\mu (\phi ^{-1}\circ E)$ are considered for each 
$E \in Af(G,\mu )$, where $Af(G,\mu )$ is
the completion of $Bf(G)$ by $\mu $-null sets, 
$Bf(G)$ is the Borel $\sigma $-field on $G$,
$\phi \circ E:=\{ \phi \circ h: h\in E \} $, $\phi \in G$.
Then $\mu $ is called quasi-invariant if there exists a dense subgroup
$G'$ such that $\mu _{\phi }$ is equivalent to $\mu $ for each $\phi \in G'$.
Henceforth, we assume that a
quasi-invariance factor $\rho _{\mu }(\phi ,g)=\mu _{\phi }(dg)/\mu (dg)$
is continuous by $(\phi ,g) \in G' \times G$,
$\mu (V)>0$ for some (open) neighbourhood $V\subset
G$ of a unit element $e \in G$ and $\mu (G)<\infty $. 
\par Let $({\sf M,F})$ be a space $\sf M$ of measures on $(G,Bf(G))$
with values in $\bf R$ and $G"$ be a dense subgroup 
in $G$ such that a topology $\sf F$ on
$\sf M$ is compatible with $G"$, that is, $\mu \mapsto \mu _h$
is the homomorphism
of $({\sf M,F})$ into itself for each $h \in G"$. Let $\sf F$ be the
topology of convergence for each $E \in Bf(G)$.
Suppose also that $G$ and $G"$ are real Banach manifolds such that
the tangent space $T_eG"$ is dense in $T_eG$, then $TG$ and $TG"$
are also Banach manifolds. Let $\Xi (G")$ denotes
a set of all differentiable vector fields $X$ on $G"$, that is, 
$X$ are sections of the tangent bundle $TG"$. We say that a measure
$\mu $ is continuously differentiable if there exists its tangent mapping
$T_{\phi }\mu _{\phi }(E)(X_{\phi })$ corresponding to the strong
differentiability relative to Banach structures of the manifolds 
$G"$ and $TG"$. Its differential we denote by $D_{\phi }\mu _{\phi }(E)$, 
hence $D_{\phi }\mu _{\phi }(E)(X_{\phi })$ is the $\sigma $-additive
real measure by subsets $E\in Af(G, \mu )$ for each $\phi \in G"$ 
and $X\in \Xi (G")$ such that $D\mu (E): TG"\to \bf R$ is continuous
for each $E\in Af(G, \mu )$, 
where $D_{\phi }\mu _{\phi }(E)=pr_2\circ (T\mu )_{\phi }(E)$,
$pr_2: p\times {\bf F}\to \bf F$ is the projection in $TN$, $p\in N$,
$T_pN=\bf F$, $N$ is another real Banach differentiable manifold modelled on 
a Banach space $\bf F$, for a differentiable mapping $V: G"\to N$ 
by $TV: TG"\to TN$ is denoted the corresponding tangent mapping,
$(T\mu )_{\phi }(E):=T_{\phi }\mu _{\phi }(E)$. 
Then by induction $\mu $ is called $n$ times 
continuously differentiable if $T^{n-1}\mu $ is continuously
differentiable such that 
$T^n\mu :=T(T^{n-1}\mu )$, $(D^n\mu )_{\phi }(E)(X_{1, \phi },
...,X_{n, \phi })$ are the $\sigma $-additive real 
measures by $E\in Af(G, \mu )$
for each $X_1,$...,$X_n\in \Xi (G")$, where $(X_j)_{\phi }=:X_{j, \phi }$
for each $j=1,...,n$ and $\phi \in G"$, $D^n\mu : Af(G, \mu )\otimes
\Xi (G")^n\to \bf R$.
\par {\bf 3.2. Note.} Suppose that in either a
$Y^{\Upsilon ,b}$-Hilbert or $Y^{\Upsilon ,b, d'}$-manifold
$N$ modelled on $l_2({\bf K})$ (see \S 2.1) there exists a dense 
$Y^{\Upsilon ,b'}$- or $Y^{\Upsilon ,b',d"}$-Hilbert 
submanifold $N'$ modelled on $l_{2,\epsilon }
=l_{2,\epsilon }({\bf K})$, ${\bf K}=\bf R$ or $\bf C$ or $\bf H$
(see \S 2.2.2), where
\par $(1)$ $a>b>b'$ and $c>d'$ and either
\par $(2)$  $\infty >\epsilon >1/2$ and $d'\ge d"$  or
\par $(3)$ $\infty >\epsilon \ge 0$ and $d'>d"$
(such that either $d'_1>d"_1$ or $d'_1=d"_1$ and $d'_2>d"_2+1$)
correspondingly.
\par If a manifold $N$ is
finite dimensional let $N'=N$. Evidently, each $Y^{\Upsilon ,b}$-manifold
is a complex $C^{\infty }$-manifold. Certainly we suppose,
that a class of smoothness of a manifold $N'$ is not less than 
that of $N$ and classes of smoothness of $M$ and $N$ are not less 
than that of a given loop group for it
as in \S 2.1.5 and of $G'$ as below.
For a chosen loop group $G=(L^MN)_{\xi }$ let its dense subgroup 
$G':= (L^MN')_{\xi '}$ be characterized by parameters:
\par $(a)$ $\xi '=(\Upsilon ,a")$ such that 
$a">b$ for $\xi ={\sf O}$ or $\xi =\sf H$
and the $Y^{\Upsilon ,b}$-manifolds $M$ and $N$
and the $Y^{\Upsilon ,b'}$-manifold $N'$;
\par $(b)$ $\xi '=(\Upsilon ,a")$ such that
$a>a">b$ for $\xi =(\Upsilon ,a)$; 
\par $(c)$ $\xi '=(\Upsilon ,a",c")$ for $\xi =(\Upsilon ,a,c)$
and $dim_{\bf K}M=\infty $ such that $b<a"<a$ and $d'<c"<c$ and either
$(2)$ $\infty >\epsilon >1$ with $d"\le d'$ or 
$(3)$ $\infty >\epsilon \ge 0$ with $d"<d'$,
such that either $d'_1>d"_1$ or $d'_1=d"_1$ and $d'_2>d"_2+1$,
where $M$ and $N$ are $Y^{\Upsilon ,b,d'}$-manifolds,
$N'$ is the $Y^{\Upsilon ,b',d"}$-manifold,
$1\le dim_{\bf K}M=:m<\infty $ in the cases $(a-b)$,
where either $a_1>a"_1$ or $a_1=a"_1$ with $a_2>a"_2+1$,
analogously for $c$ and $c"$, $b$ and $b'$ instead of $a$ and $a"$.
For the corresponding pair 
$G':=(L^M_{\bf R}N')_{\xi '}$ and $G:=(L^M_{\bf R}N)_{\xi }$
let indices in $(1-3)$ and $(a-c)$ be the same with substitution
of $\xi =\sf O$ on $\xi =(\infty ,H)$.
\par For a diffeomorphism group $Diff^t_{\beta ,\gamma }(\tilde M)$
of a Banach manifold $\tilde M$ let $M$ be a dense Hilbert submanifold
in $\tilde M$ as in \cite{lurim1,lurim2}.
For $Diff^{\xi }_{y_0}(N)$ and $(L^MN)_{\xi }$ consider
in $G:=Diff^{\xi }_{y_0}(N) \otimes ^s(L^MN)_{\xi }$ a dense
subgroup $G':=Diff^{\xi '}_{y_0}(N') \otimes ^s(L^MN')_{\xi '}$,
where pairs $(\xi ,\xi ')$ are described above.
\par {\bf 3.3. Theorem.} {\it Let $G$ be either a loop group
or a diffeomorphism group or their semidirect product
for real or complex or quaternion
separable metrizable $C^{\infty }$-manifolds $M$ and $N$,
then there exist a Wiener process on $G$ which induces
quasi-invariant infinite differentiable measures $\mu $
relative to a dense subgroup $G'$.
\par For a given pair
$(G,G')$ there exists a family of nonequivalent Wiener processes on $G$
of the cardinality ${\sf c}=card ({\bf R})$
and a family of the cardinality $\sf c$ of pairwise
orthogonal quasi-invariant $C^{\infty }$-differentiable
measures on $G$ relative to $G'$.}
\par {\bf Proof.} These topological groups also have structures 
of $C^{\infty }$-manifolds, but they do not satisfy the 
Campbell-Hausdorff formula (see Theorems 2.1.7, 2.8, 2.10)
in any open local subgroup.
Their manifold structures and actions of $G'$ on $G$ will be sufficient
for the construction of desired measures. Manifolds over $\bf C$ or
$\bf H$ naturally have structures of manifolds over $\bf R$ also.
\par We take $G=\bar G$ and $Y=\bar Y$ for each loop group
$(L^MN)_{\xi }$ outlined in 3.2.$(b,c)$, for each diffeomorphism group 
$Diff^{\xi }(M)$ of a complex or quaternion manifold $M$ respectively
given above, for each diffeomorphism group 
$G := Di^{ \{ l \} }_{ \{ \gamma \} , \delta ,\eta }(M)$ 
for a real manifold $M$, for corresponding semidirect products
of loop and diffeomorphism groups, since such $G$
has a Hilbert manifold structure (see Theorems 2.1.7, 2.8, 2.10).
For ${\bar G}:=Diff^t_{\beta ,\gamma }(\tilde M)$ there exists
a Hilbert dense submanifold $M$ in a Banach manifold $\tilde M$ 
(see \S 2.6) and a subgroup
$G:=Di^{ \{ l \} }_{ \{ \gamma \} , \delta ,\eta }(M)$  
dense in $\bar G$ and a diffeomorphism subgroup
$G'$ dense in $G$ (see the proof of Theorem 3.10 \cite{lurim2}
and Lemma 2.6.2 above), analogously for loop groups of such classes
of smoothness for real manifolds and corresponding semidirect
products of these groups. This $G'$ can be chosen as in Lemma 2.7.
\par For the chosen loop group $G=(L^MN)_{\xi }$ let its dense subgroup 
$G':= (L^MN')_{\xi '}$ be the same as in \S 3.2 Cases $(b,c)$.
In case $3.2.(a)$ let $\bar G=(L^MN)_{\xi }$ and $G=(L^MN')_{\hat \xi }$
with $\hat \xi =(\Upsilon ,\hat a)$ such that $\hat a>a"$,
then $G'$ let be as in $3.2.(a)$, then also for $G := Diff^{\xi }_{y_0}
(N)\otimes ^s(L^MN)_{\xi }$ take $G' :=
Diff^{\xi '}_{y_0}(N')\otimes ^s (L^MN')_{\xi '}$.
\par Then the embedding $J: T_eG' \hookrightarrow T_eG$
is the Hilbert-Schmidt operator, that follows from \S 2.1 and Lemma 2.7,
Theorem 2.10.
\par On a dense subgroup $G'$ there exists a $1$-parameter group
$\rho : {\bf R}\times G'\to G'$
of diffeomorphisms of $G'$ generated by a $C^{\infty }$-vector field 
$X_{\rho }$ on $G'$ such that $X_{\rho }(p)=(d\rho (s,p)/ds )|_{s=0}$,
where $\rho (s+t,p)=\rho (s,\rho (t,p))$ for each $s, t\in \bf R$,
$\rho (0,p)=p$, $\rho (s,*): G'\to G'$ is the diffeomorphism
for each $s\in \bf R$ (about $\rho $ see \S 1.10.8 \cite{kling}).
Then each measure $\mu $ on $G$ and $\rho $ produce a $1$-parameter
family of measures $\mu _s(W):=\mu (\rho (-s,W))$.
Let $\tau _{G}: TG\to G$ be a tangent bundle on $G$. 
Let also $\theta : Z_{G}\to G$ be a trivial bundle on $G$
with a fibre $Z$ such that $Z_{G}=Z\times G$. We suppose also, that
$L_{1,2}(\theta , \tau _{G})$ is an operator bundle with a fibre
$L_{1,2}(Z,Y)$, where $Z, Z_1,...,Z_n$ are Hilbert spaces, 
$L_{n,2}(Z_1,...,Z_n;Z)$ is
a subspace of a space of all Hilbert-Schmidt $n$ times multilinear 
operators from $Z_1\times ...\times Z_n$ into $Z$
(see \cite{beldal,pietsch}). Then
$L_{n,2}(Z_1,...,Z_n;Z)$ has the structure of the Hilbert space
with the scalar product denoted by 
$$\sigma _2(\phi ,\psi ):=
\sum_{j_1,...,j_n=1}^{\infty }(\phi (e^{(1)}_{j_1},...,e^{(n)}_{j_n}),
\psi (e^{(1)}_{j_1},...,e^{(n)}_{j_n}))$$
for each pair of its elements $\phi , \psi .$ It
does not depend on a choice of the orthonormal bases
$\{ e{(k)}_j: j \} $ in $Z_k$. Let $\Pi :=
\tau _{G}\oplus L_{1,2}(\theta ,\tau _{G})$ be a Whitney sum
of bundles $\tau $ and $L_{1,2}(\theta ,\tau _{G})$. 
If $(U_j,\phi _j)$ 
and $(U_l,\phi _l)$ are two charts of $G$ with an open non-void intersection
$U_j\cap U_l$, then to a connecting mapping $f_{\phi _l,\phi _j}=
\phi _l\circ \phi _j^{-1}$ there corresponds a connecting mapping
$f_{\phi _l,\phi _j}\times {f'}_{\phi _l,\phi _j}$ for the bundle
$\Pi $ and its charts $U_j\times (Y\oplus L_{1,2}(Z,Y))$ for $j=1$ 
or $j=2$, where $f'$ denotes the strong derivative of $f$,
${f'}_{\phi _l,\phi _j}: (a^{\phi _j},A^{\phi _j})\mapsto
({f'}_{\phi _l,\phi _j}a^{\phi _j}, {f'}_{\phi _l,\phi _j}
\circ A^{\phi _j})$, $a^{\phi }\in Y$ and $A^{\phi }\in L_{1,2}(Z,Y)$ 
for the chart $(U,\phi )$,
${f'}_{\phi _l,\phi _j}\circ A^{\phi _j}:=
{f'}_{\phi _l,\phi _j}A^{\phi _j}{f'}_{\phi _l,\phi _j}^{-1}$. 
Such bundles are called quadratic.
Then there exists a new bundle $J$ on $G$ with the same fibre as 
for $\Pi $, but with new connecting mappings: 
$J(f_{\phi _l,\phi _j}): (a^{\phi _j},A^{\phi _j})\mapsto
({f'}_{\phi _l,\phi _j}a^{\phi _j}+tr ({f"}_{\phi _l,\phi _j}(
A^{\phi _j}, A^{\phi _j}))/2, {f'}_{\phi _l,\phi _j}\circ A^{\phi _j})$,
where $tr (A)$ denotes a trace of an operator $A$.
Then using sheafs one gets the It$\hat o$ functor  $I: I(G)\to G$
from the category of manifolds to the category of quadratic bundles.
\par For the construction of differentiable measures on the
$C^{\infty }$-manifold we shall use the following statement:
if $a\in C^{\infty }(TG',TG)$ and $A\in C^{\infty }
(TG',L_{1,2}(TG',TG))$ and $a_x\in T_xG$ and $A_x\in L_{1,2}(T_xG',T_xG)$
for each $x\in G'$, each derivative by $x\in G'$:
$a^{(k)}_x$ and $A^{(k)}_x$ is a
Hilbert-Schmidt mappings into $Y=T_eG$ for each 
$k\in \bf N$ and $\sup_{\eta \in G} \| A_{\eta }(t)A_{\eta }^*(t) \| ^{-1}
\le C$, where $C>0$ is a constant,
then the transition probability $P(\tau ,x,t,W):=
P \{ \omega : \xi (t,\omega )=x, \xi (t,\omega )\in W \} $ 
is continuously stronlgy $C^{\infty }$-differentiable 
along vector fields on $G'$, where $G'$ is  
a dense $C^{\infty }$-submanifold on a space $Y'$,
$Y'$ is a separable real Hilbert space having embedding into $Y$
as a dense linear subspace (see Theorem 3.3 and the Remark
after it in Chapter 4 \cite{beldal}
as  well as Theorems 4.2.1, 4.3.1 and 5.3.3 \cite{beldal}, 
Definitions 3.1 above), $W\in {\sf F}_t$.
\par Now let $G$ be a loop or a diffeomorphism group or their semidirect
product of the corresponding manifolds 
over the field ${\bf K}=\bf R$ or $\bf C$ or $\bf H$.
Then $G$ has the manifold structure.
If $\exp^N: {\tilde T}N\to N$
is an exponential mapping of the manifold $N$, then
it induces the exponential $C^{\infty }$-mapping
${\tilde E}: {\tilde T}(L^MN)_{\xi }\to (L^MN)_{\xi }$
defined by ${\tilde E}_{\eta }(v)=\exp^N_{\eta }\circ v_{\eta }$
(see Theorems 2.1.3.9, 2.1.7 and 2.10),
where ${\tilde T}N$ is a neighbourhood of $N$ in a tangent
bundle $TN$, $\eta \in (L^MN)_{\xi }=:G$, $W_e$ is a neighbourhood 
of $e$ in $G$, $W_{\eta }=W_e\circ \eta .$
At first this mapping is defined on classes of equivalent 
mappings of the loop monoid $(S^MN)_{\xi }$ 
and then on elements of the group, since
$\exp^N_{f(x)}$ is defined for each $x\in M$ and
$f\in \eta \in (S^MN)_{\xi }$.
The manifolds $G$ and $G'$ are of class $C^{\infty }$
and the exponential mappings $\tilde E$ and $\bar E$ for $G$ and $G'$
correspondingly are of class (strongly) $C^{\infty }$.
The analogous connection there exists in the diffeomorphism group
of the manifold $M$ satisfying the corresponding conditions
(see Theorem 3.3 \cite{lurim1}, \S 2.3 and Theorem 2.8)
for which: ${\tilde E}_{\eta }(v)=\exp_{\eta (x)}\circ v_{\eta }$ for each 
$x\in M$ and $\eta \in G$. 
We can choose the uniform atlases $At_u(G)$ such that Christoffel symbols
$\Gamma _{\eta }$ are bounded on each chart (see Theorem 2.11).
This mapping $\tilde E$ is for $G$ as the manifold and has not relations with
its group structure such as given by the Campbell-Hausdorff formula
for some Lie group, for example, finite dimensional Lie group. 
For the case of manifolds $M$ and $N$ over $\bf C$
we consider $G$ and others appearing manifolds with their structure 
over $\bf R$, since ${\bf C}={\bf R}\oplus i \bf R$
as the Banach space over $\bf R$.
\par The exponential mapping
$\exp^{G}: {\tilde T}G\to G$ is defined by the formula $X\mapsto c_X(1).$
The restriction $\exp^G|_{\tilde TG\cap T_pG}$ will also be denoted
by $\exp^{G}_p$. Then there is defined the mapping
$I(\exp^{G}): I({\tilde T} G)\to I(G)$ such that for each chart
$(U,\phi )$ the mapping
$I(\exp^{\phi }): Y\oplus L_{1,2}(Z,Y)\to Y\oplus L_{1,2}(Z,Y)$
is given by the following formula: 
$$I(\exp^{\phi })(a^{\phi },
A^{\phi })=(a^{\phi }- \quad tr ( \Gamma ^{\phi }(A^{\phi },
A^{\phi }))/2,A^{\phi }),$$ 
where $\Gamma $ denotes the Christoffel symbol.
\par Therefore, if ${\sf R}_{x,0}(a,A)$ is a germ of diffusion processes
at a point $y=0$ of the tangent space $T_xG$, then
${\tilde {\exp }}_x{\sf R}_{x,0}(a,A):={\sf R}_x(I( \exp_x )(a,A))$ is a germ
of stochastic processes at a point $x$ of the manifold $G$.
The germs ${\tilde {\exp }}_x{\sf R}_{x,0}(a,A)$ are called stochastic 
differentials and the It$\hat o$ bundle is called the bundle of 
stochastic differentials such that ${\sf R}_{x,0}(a,A)=:
a_xdt+A_xdw$. A section $\sf U$ of the vector bundle
$\Pi =\tau _{Y}\oplus L_{1,2}(\theta ,\tau _{Y})$ is called an 
It$\hat o$ field on the manifold $G$ and it defines a field of
stochastic differentials
${\sf R}_x(I(exp_x)(a,A))={\tilde exp}_x(a_xdt+A_xdw)$.
A random process $\xi $ has a stochastic differential defined by the 
It$\hat o$ field $\sf U:$ $d\xi (s,\omega )={\tilde {\exp }}_{\xi (s,\omega )}
{\sf R}(a_{\xi (s,\omega )},A_{\xi (s,\omega )})$ 
if the following conditions are 
satisfied: for $\nu  _{\xi (s)}$-almost every $x\in Y$
there exists a neighbourhood $V_x$ of a point $x$ and a diffusion process
$\eta _x(t,\omega )$ belonging to the germ ${\sf R}_x(I( \exp_x ))(a,A)$
such that $P_{s,x} \{ \xi (t,\omega )=\eta _x(t,\omega ): 
\xi (t,\omega )\in V_x,
t\ge s \} =1$ $\nu _{\xi (s)}$-almost everywhere, where
$P_{s,x}(S):=P \{ S: \xi (s,\omega )=x \} $, $S$ is a $P$-measurable
subset of $\Omega $, $\nu _{\xi (s)}(F):= P \{ \omega :
\xi (s,\omega )\in  F \} $ (see Chapter 4 in \cite{beldal}).                     
\par If ${\sf U}(t)=(a(t), A(t))$ is a time dependent It$\hat o$ field,
then a random process $\xi (t,\omega )$ 
having for each $t\in [0,T]$ a stochastic 
differential $d\xi =\exp_{\xi (t,\omega )}(a_{\xi (t,\omega )}dt+
A_{\xi (t,\omega )}dw)$
is called a stochastic differential equation on the manifold
$G$, the process $\xi (t,\omega )$ is called its solution (see Chapter 
VII in \cite{dalf}).
As usually a flow of $\sigma $-algebras consistent with
the Wiener process $w(t,\omega )$ is a monotone set
of $\sigma $-algebras ${\sf F}_t$ such that $w(s,\omega )$ is
${\sf F}_t$-measurable for each $0\le s\le t$ 
and $w(\tau ,\omega )-w(s,\omega )$ is independent from ${\sf F}_t$
for each $\tau >s\ge t$, where ${\sf F}_s\supset {\sf F}_t$ 
for each $0\le t\le s$.
\par Then we consider for a manifold $G$ its It$\hat o$ bundle
for which an It$\hat o$ field $\sf U$ has a principal part
$(a_{\eta },A_{\eta })$, where $a_{\eta }\in T_{\eta }G$ and
$A_{\eta } \in L_{1,2}(H,T_{\eta }G)$ and $ker (A_{\eta })= \{ 0 \} $, 
$\theta : H_G\to G$ is a trivial bundle 
with a Hilbert fiber $H$ and $H_G:=G\times H$,
$L_{1,2}(\theta ,\tau _{\eta })$
is an operator bundle with a fibre $L_{1,2}(H,T_{\eta }G)$. 
To satisfy conditions of the theorem about quasi-invariance and
differentiability of a transition probability  
we choose $A$ also such that $\sup_{\eta \in G} \| A_{\eta }(t)
A_{\eta }^*(t) \| ^{-1} \le C$, where $C>0$ is a constant.
If an operator $B$ is selfadjoint, then
$A_{\eta }^{\phi }B{A_{\eta }^{\phi }}^*$ is also selfadjoint, where
$A_{\eta }(t)=:A_{\eta }^{\phi _j}(t)$ is on a chart $(U_j,\phi _j)$.
If $\mu _B$ is a Gaussian measure on $T_{\eta }G$ with the correlation 
operator $B$, then $\mu _{A_{\eta }^{\phi }B{A_{\eta }^{\phi }}^*}$
is the Gaussian measure on $X_{1,\eta },$
where $B$ is selfadjoint and $ker (B)= \{ 0 \} $,
$A_{\eta }: T_{\eta }G\to X_{1,\eta }$, $X_{1,\eta }$ is a Hilbert space. 
We can take initially $\mu _B$ a cylindrical measure on
a Hilbert space $X'$ such that
$T_{\eta }G'\subset X'\subset T_{\eta }G$. 
If $A_{\eta }$ is the Hilbert-Schmidt operator
with $ker (A_{\eta })= \{ 0 \} $, then
$A_{\eta }^{\phi }B{A_{\eta }^{\phi }}^*$ is the nondegenerate selfadjoint
linear operator of trace class and the so called Radonifying 
operator $A_{\eta }^{\phi }$ induces the $\sigma $-additive measure 
$\mu _{A_{\eta }^{\phi }B{A_{\eta }^{\phi }}^*}$
in the completion $X'_{1,\eta }$ of $X'$
with respect to the norm $\| x \| _1:= \| A_{\eta }x \| $
(see \S II.2.4 \cite{dalf}, \S I.1.1 \cite{sko}, 
\S II.2.4 \cite{oksb}). Then
using cylinder subsets we get a new Gaussian $\sigma $-additive
measure on $T_{\eta }G$, which we denote also by
$\mu _{A_{\eta }^{\phi }B{A_{\eta }^{\phi }}^*}$
(see also Theorems I.6.1 and III.1.1 \cite{kuo}).
\par If $U_j\cap U_l\ne \emptyset $, then
$A_{\eta }^{\phi _l}(t)={f_{\phi _l,\phi _j}}'
A_{\eta }^{\phi _l}(t){{f_{\phi _l,\phi _j}}'}^{-1}$, hence
the correlation operator $A_{\eta }^{\phi }B{A_{\eta }^{\phi }}^*$
is selfadjoint on each chart of $G$, that produces the Wiener 
process correctly.
Therefore, we can consider a stochastic process
$\mbox{ }d\xi (t,\omega )={\tilde E}_{\xi (t,\omega )}
[a_{\xi (t,\omega )}dt+A_{\xi (t,\omega )}dw],$
where $w$ is a Wiener process in $T_{\eta }G$
defined with the help of a nuclear nondegenerate 
selfadjoint positive definite operator $B$.
The corresponding Gaussian measures 
$\mu _{tA_{\eta }^{\phi }B{A_{\eta }^{\phi }}^*}$ for $t>0$
(for the Wiener process)  are defined on the Borel
$\sigma $-algebra of $T_{\eta }G$ and 
$\mu _{tA_{\eta }^{\phi }B{A_{\eta }^{\phi }}^*}$
for such Hilbert-Schmidt nondegenerate linear operators
$A_{\eta }$ with $ker (A_{\eta })= \{ 0 \} $
are $\sigma $-additive (see Theorem II.2.1 \cite{dalf}).
If the embedding operator $T_{\eta }G'\hookrightarrow T_{\eta }G$
is of the Hilbert-Schmidt class, then there exist $A_{\eta }$ and $B$
such that $\mu _{tA_{\eta }^{\phi }B{A_{\eta }^{\phi }}^*}$
is the quasi-invariant and $C^{\infty }$-differentiable
measure on $T_{\eta }G$ relative to shifts on vectors from $T_{\eta }G'$
(see Theorem 26.2 \cite{sko} using Carleman-Fredholm determinant
and Chapter IV \cite{dalf} and \S 5.3 \cite{ustzak}).
Henceforth we impose such conditions
on $B$ and $A_{\eta }$ for each $\eta \in G'$.
\par Consider left shifts $L_h: G\to G$
such that $L_h\eta :=h\circ \eta $. 
Let us take $a_e\in T_eG$, $A_e\in L_{1,2}(T_eG',T_eG)$, 
then we put $a_x=(DL_x)a_e$ and $A_x=(DL_x)\circ A_e$
for each $x\in G$, hence $a_x\in T_eG$ and 
$A_x\in L_{1,2}(H_x,(DL_x)T_eG)$, where $(DL_x)T_eG=T_xG$
and $T_eG'\subset T_eG$, $H_x:=(DL_x)T_eG'$.
Operators $L_h$ are (strongly) $C^{\infty }$-differentiable
diffeomorphisms of $G$ such that $D_hL_h: T_{\eta }G\to 
T_{h\eta }G$ is correctly defined, since $D_hL_h=h_*$
is the differential of $h$ \cite{ebi,eichh}.
In view of the choice of $G'$ in $G$ each 
covariant derivative $\nabla _{X_1}...\nabla _{X_n}(D_hL_h)Y$
is of class $L_{n+2,2}({TG'}^{n+1}\times G',TG)$ 
for each vector fields $X_1,...,X_n,Y$ on $G'$
and  $h\in G'$, since for each $0\le l\in \bf Z$ the embedding of $T^lG'$ 
into $T^lG$ is of Hilbert-Schmidt 
class, where $T^0G:=G$ (above and in \cite{beldal}
mappings of trace and Hilbert-Schmidt classes were defined for
linear mappings on Banach and Hilbert spaces and then 
for mappings on vector bundles). Take a dense subgroup $G'$
as it was otlined above and consider left shifts $L_h$ for $h\in G'$.
\par The considered here groups $G$ are separable, 
hence the minimal $\sigma $-algebra generated by all cylindrical 
subalgebras $f^{-1}({\sf B}_n)$, n=1,2,..., coincides with 
the $\sigma $-algebra $\sf B$ of all Borel subsets of $G$, where 
$f: G\to \bf R^n$ is a continuous function, ${\sf B}_n$ 
is the Borel $\sigma $-algebra  of $\bf R^n$. Moreover, $G$ 
is the topological Radon space (see Theorem I.1.2 and 
Proposition I.1.7 \cite{dalf}).
Let $P(t_0,\psi ,t,W):=P( \{ \omega : \xi (t_0,\omega )=\psi ,
\xi (t,\omega )\in W \} )$ be a transition probability of 
a stochastic process $\xi $ for $0\le t_0<t$, which is defined on a 
$\sigma $-algebra $\sf B$ of Borel subsets in $G$, $W\in \sf B$,
since each measure $\mu _{A_{\eta }^{\phi }B{A_{\eta }^{\phi }}^*}$
is defined on the $\sigma $-algebra
of Borel subsets of $T_{\eta }G$ (see above).
\par If $G$ is a manifold with an uniform atlas (see \S 2.1) such that
an It$\hat o$ field $(a,A)$ and Christoffel symbols are bounded, then
there exists a unique up to stochastic equivalence 
random evolution family $S(t,\tau )$ consistent
with the flow of $\sigma $-algebras ${\sf F}_t$ generated by a 
solution $\xi (t,\omega )$ of the stochastic differential equation 
$d\xi =\exp_{\xi (t,\omega )}(a_{\xi (t,\omega )}dt+A_{\xi (t,\omega )}dw)$
on $G$, that is, $\xi (\tau ,\omega )=x$, $\xi (t,\omega )=S(t,\tau ,
\omega )x$
for each $t_0\le \tau <t<\infty $ (see Theorem 4.2.1 \cite{beldal}).
\par On the other hand, $S(t,\tau ;gx)=gS(t,\tau ;x)$
is the stochastic evolution family of operators for each
$0\le t_0\le \tau <t$. There exists $\mu (W):=P(t_0,\psi , t,W)$ such that
it is a $\sigma $-additive
quasi-invariant strongly $C^{\infty }$-differentiable relative to 
the action of $G'$ by the left shifts $L_h$ on $\mu $
measure on $G$, for example, $t_0=0$ and $\psi =e$
with $t_0<t$, that is, $\mu _h(W):=\mu (h^{-1}W)$
is equivalent to $\mu $ and it is strongly infinitely differentiable 
by $h\in G'$.
\par The proof in cases $G=\bar G$ is thus obtained.
In cases $G\subset \bar G$ and $G\ne \bar G$ the use of the standard 
procedure of cylinder subsets induce a Weiener process and a transition 
probability from $G$ on $\bar G$ which is quasi-invariant and 
$C^{\infty }$-differentiable relative to $G'$ (see aslo 
\cite{lurim2}).
Evidently, considering different $(a,A)$ we see that there exists
a family of nonequivalent Wiener processes on $G$
of the cardinality ${\sf c}=card ({\bf R})$.
In view of the Kakutani theorem in \cite{dalf}
there exists a family of the cardinality $\sf c$ of pairwise
orthogonal quasi-invariant $C^{\infty }$-differentiable
measures on $G$ relative to $G'$.
\par {\bf 3.4. Note.} This proof also shows,
that $\mu $ is infinitely differentiable 
relative to each $1$-parameter group $\rho : {\bf R}\times G'\to G'$
of diffeomorphisms of $G'$ generated by a $C^{\infty }$-vector field 
$X_{\rho }$ on $G'$. 

\section{Unitary representations associated with quasi-invariant measures.}
\par {\bf 4.1.1. Note.} A transition probability $P=:\nu $ on $G$
induces strongly continuous unitary regular representation
of $G'$ given by the following formula:
$T_h^{\nu }f(g):=(\nu ^h(dg)/\nu (dg))^{(1+bi)/2}f(h^{-1}g)$
for $f\in L^2(G,\nu ,{\bf C})=:H$, $T_h^{\nu }\in U(H)$,
$U(H)$ denotes the unitary group of the Hilbert space $H$,
where $b\in \bf R$, $i = (-1)^{1/2}$.
For the strong continuity of $T_h^{\nu }$ conditions of
the continuity of the mapping $G'\ni h\mapsto \rho _{\nu }(h,g)\in
L^1(G,\nu ,{\bf C})$ and that $\nu $ is the Borel measure
are sufficient, where $g\in G$, since $\nu $ is the 
Radon measure (see its definition in Chapter I \cite{dalf}). 
On the other hand, the continuity of $\rho _{\nu }(h,g)
=\nu ^h(dg)/\nu (dg)$ by $h$
from a Polish group $G'$ into $L^1(G,\nu ,{\bf C})$ follows from
the inclusion $\rho _{\nu }(h,g)\in L^1(G,\nu ,{\bf C})$  for each 
$h\in G'$ and that $G'$ is a topological subgroup of $G$.
In section 3 mostly Polish groups $\bar G$ and $G'$ were considered.
When $\bar G$ was not Polish it was used an embedding into $\bar G$
of a Polish subgroup $G$ such that $G'\subset G\subset \bar G$
and a measure on $G$ induces a measure on $\bar G$ with the help 
of an algebra of cylindrical subsets. So the considered 
cases of representations reduce to the case of Polish groups $(G,G')$.
\par More generally it is possible to consider
instead of a topological group $G$ a Polish topological space
$X$ on which $G'$ acts jointly continuously: $\phi : (G'\times X)\ni
(h,x)\mapsto hx=:\phi (h,x)\in X$, $\phi (e,x)=x$
for each $x\in X$, $\phi (v,\phi (h,x))=\phi (vh,x)$
for each $v$ and $h\in G'$ and each $x\in X$.
If $\phi $ is a Borel function, then
it is jointly continuous \cite{fidal}.
\par A representation $T: G'\to U(H)$ is called topologically irreducible,
if there is not any unitary operator (homeomorphism) $S$ on $H$ and a 
closed (Hilbert) subspace $H'$ in $H$ such that 
$H'$ is invariant relative to
$ST_hS^*$ for each $h\in G'$, that is, $ST_hS^*(H')\subset H'$.
\par For a topological space $S$ let $S^d$ denotes the derivative
set of $S$,
that is, of all limit points $x\in cl(S\setminus \{ x \} ) $, $x\in S$, 
where $cl(A)$ denotes the closure of a subset $A$ in $S$
(see \S 1.3 \cite{eng}).
A topological space $S$ is called dense in itself if
$S\subset S^d$.
\par A measure $\nu $ on $X$ is called ergodic, if for
each $U\in Af(X,\nu )$ and $F\in Af(X,\nu )$ with $\nu (U)
\times \nu (F)\ne 0$
there exists $h\in G'$ such that $\nu ((h\circ E)\cap F)\ne 0$.  
\par {\bf 4.1.2. Theorem.} {\it Let $X$ be an infinite
Polish topological space 
with a $\sigma $-additive $\sigma $-finite nonnegative nonzero 
ergodic Borel measure $\nu $ with  $supp (\nu )=X$ and
quasi-invariant relative to an infinite dense in itself 
Polish topological group $G'$ acting on $X$ by a Borel function
$\phi $. If
\par $(i)$ $sp_{\bf C} \{\psi |
\quad \psi (g):=(\nu ^h(dg)/\nu (dg))^{(1+bi)/2}, h\in G' \} $
is dense in $H$, where $b\in \bf R$ is fixed, and
\par $(ii)$ for each $f_{1,j}$ and $f_{2,j}$ in $H$, $j=1,...,n,$ 
$n\in \bf N$ and each $\epsilon >0$ there exists $h\in G'$ such that
$|(T_hf_{1,j},f_{2,j})| \le \epsilon |(f_{1,j},f_{2,j})|$,
when $|(f_{1,j},f_{2,j})|>0$.
Then the regular representation $T: G'\to U(H)$ is topologically
irreducible.}
\par {\bf Proof.}  From Condition $(i)$ it follows, that
the vector $f_0$ is cyclic, where $f_0\in H$ 
and $f_0(g)=1$ for each $g\in X$. 
In view of $card (X)\ge \aleph _0$ and an ergodicity of $\nu $
for each $n\in \bf N$ there are subsets $U_j\in Bf(X)$ and
elements $g_j\in G'$
such that $\nu ((g_jU_j)\cap (\bigcup_{i=1,...,j-1,j+1,...,n}U_i))=0$
and $\prod_{j=1}^n\nu _j(U_j)>0$.
Together with Condition $(ii)$ this implies, 
that there is not any finite dimensional 
$G'$-invariant subspace $H'$ in $H$ such that
$T_hH'\subset H'$ for each $h\in G'$ and $H'\ne \{ 0 \}$.
Hence if there is a $G'$-invariant closed subspace $H'\ne 0$
in $H$ it is isomorphic with the subspace
$L^2(V,\nu ,{\bf C})$, where $V\in Bf(X)$ with $\nu (V)>0$. 
\par Let ${\sf A}_G$ denotes a $*$-subalgebra of 
an algebra ${\sf L}(H)$ of bounded linear operators on $H$
generated by the family of unitary operators 
$\{ T_h: h\in G' \} $. In view of the von Neumann
double commuter Theorem (see \S VI.24.2 \cite{fell})
${{\sf A}_G}"$ coincides with the weak and strong operator closures of
${\sf A}_G$ in ${\sf L}(H)$, where ${{\sf A}_G}'$
denotes the commuting algebra of ${\sf A}_G$ and ${{\sf A}_G}"=
({{\sf A}_G}')'$. 
\par Each Polish space is \v{C}ech-complete. By the 
Baire category theorem in a \v{C}ech-complete space $X$
the union $A=\bigcup_{i=1}^{\infty }A_i$ of a sequence of
nowhere dense subsets $A_i$ is a codense subset (see
Theorem 3.9.3 \cite{eng}). On the other hand, in view of Theorem
5.8 \cite{hew} a subgroup of a topological group
is discrete if and only if  it contains an isolated point.
Therefore, we can choose
\par $(i)$ a probability 
Radon measure $\lambda $ on $G'$ such that $\lambda $ has not any atoms and
$supp (\lambda )=G'$. In view of the strong continuity of
the regular representation there exists the S. Bochner integral
$\int_XT_hf(g)\nu (dg)$ for each $f\in H$, which implies its existence 
in the weak (B. Pettis) sence. The measures $\nu $ and $\lambda $
are non-negative and bounded, hence $H\subset L^1(X,\nu ,{\bf C})$
and $L^2(G',\lambda ,{\bf C})\subset L^1(G',\lambda ,{\bf C})$
due to the Cauchy inequality. Therefore, we can apply below 
the Fubini Theorem (see \S II.16.3 \cite{fell}).
Let $f\in H$, then there exists a countable orthonormal base
$\{ f^j: j\in {\bf N} \} $ in $H\ominus {\bf C}f$. Then for each
$n\in \bf N$ the following set $B_n:=\{ q\in L^2(G',\lambda ,{\bf C} ):$
$(f^j,f)_H=\int_{G'}q(h)(f^j,T_hf_0)_H\lambda (dh)$ for $j=0,...,n \} $
is non-empty, since the vector $f_0$ is cyclic, where $f^0:=f$. 
There exists $\infty >R>\| f\|_H$ such that $B_n\cap B^R=:B^R_n$
is non-empty and weakly compact for each $n\in \bf N$, 
since $B^R$ is weakly compact, where
$B^R:=\{ q\in L^2(G',\lambda ,{\bf C} ): \| q\| \le R \} $
(see the Alaoglu-Bourbaki Theorem in \S (9.3.3) \cite{nari}).
Therefore, $B_n^R$ is a centered system of closed subsets
of $B^R$, that is, $\cap_{n=1}^mB^R_n\ne \emptyset $
for each $m\in \bf N$, hence it has a non-empty intersection, consequently,
there exists $q\in L^2(G',\lambda ,{\bf C})$ such that
$$(ii)\mbox{ }f(g)=\int_{G'}q(h)T_hf_0(g)\lambda (dh)$$ for $\nu $-a.e.
$g\in X$.
If $F\in L^{\infty }(X,\nu ,{\bf C})$, $f_1$ and $f_2\in H$,
then there exist $q_1$ and $q_2\in L^2(G',\lambda ,{\bf C})$
satisfying Equation $(ii)$. Therefore, 
$$(iii)\mbox{ }(f_1,Ff_2)_H=:c=$$
$$\int_X\int_{G'}\int_{G'}{\bar q}_1(h_1)q_2(h_2)
\rho _{\nu }^{(1+bi)/2}(h_1,g)
\rho _{\nu }^{(1+bi)/2}(h_2,g)F(g)\lambda (dh_1)\lambda (dh_2)\nu (dg).$$
$$\mbox{Let }\xi (h):=\int_X\int_{G'}\int_{G'}{\bar q_1}(h_1)q_2(h_2)
\rho _{\nu }^{(1+bi)/2}(h_1,g) \rho _{\nu }^{(1+bi)/2}(hh_2,g)
\lambda (dh_1)\lambda (dh_2) \nu (dg).$$
Then there exists $\beta (h)\in L^2(G',\lambda ,{\bf C})$
such that 
\par $(iv)$ $\int_{G'}\beta (h)\xi (h)\lambda (dh)=c$.\\
To prove this we consider two cases. If $c=0$ it is sufficient
to take $\beta $ orthogonal to $\xi $ in $L^2(G',\lambda ,{\bf C})$. 
Each function $q\in L^2(G',\lambda ,{\bf C})$ 
can be written as $q=q^1-q^2+iq^3-iq^4$,
where $q^j(h)\ge 0$ for each $h\in G'$ and $j=1,...,4$,
hence we obtain the corresponding decomposition for $\xi $,
$\xi =\sum_{j,k}b^{j,k}\xi ^{j,k}$, where $\xi ^{j,k}$ corresponds to
$q_1^j$ and $q_2^k$, where $b^{j,k}\in \{ 1,-1,i,-i \}$. 
If $c\ne 0$ we can choose $(j_0,k_0)$ for which $\xi ^{j_0,k_0}\ne 0$
and 
\par $(v)$ $\beta $ is orthogonal to others $\xi ^{j,k}$ with 
$(j,k)\ne (j_0,k_0)$.\\ 
Otherwise, if $\xi ^{j,k}=0$ for each
$(j,k)$, then $q_l^j(h)=0$ for each $(l,j)$ and $\lambda $-a.e. $h\in G'$,
since 
$$\xi (0)=\int_X\nu (dg)(\int_{G'}{\bar q_1}(h_1)\rho _{\nu }^{(1+bi)/2}
(h_1,g)\lambda (dh_1))(\int _{G'}q_2(h_2)
\rho _{\nu }^{(1+bi)/2}(h_2,g)\lambda (dh_2))=0$$ 
and this implies $c=0$, which 
is the contradiction with the assumption $c\ne 0$.
Hence there exists $\beta $ satisfying conditions $(iv, v)$.
\par Let $a(x)\in L^{\infty }(X,\nu ,{\bf C})$, $f$ and $g\in H$, 
$\beta (h)\in L^2(G',\lambda ,{\bf C})$. Since $L^2(G',\lambda ,{\bf C})$ 
is infinite dimensional, then for each finite family of 
$a\in \{ a_1,...,a_m \} \subset L^{\infty }(X,\nu ,{\bf C})$,
$f\in \{ f_1,...,f_m \} \subset H$ there exists
$\beta (h)\in L^2(G',\lambda ,{\bf C})$, $h\in G'$, such that
$\beta $ is orthogonal to $\int_X{\bar f}_s(g)
[f_j(h^{-1}g)
(\rho _{\nu }(h,g))^{(1+bi)/2}-f_j(g)]\nu (dg)$ for each $s,j=1,...,m$. Hence
each operator of multiplication on $a_j(g)$
belongs to ${{\sf A}_G}"$, since due to Formula $(iv)$
and Condition $(v)$ there exists $\beta (h)$ such that 
$$(f_s,a_jf_l)=\int_X\int_{G'}{\bar f}_s(g)\beta (h)(\rho _{\nu }
(h,g))^{(1+bi)/2}f_l(h^{-1}g)\lambda (dh) \nu (dg)=$$
$$=\int_X\int_{G'} {\bar f}_s(g)
\beta (h)(T_hf_l(g))\lambda (dh)\nu (dg)\mbox{, }
\int_X{\bar f}_s(g)a_j(g)f_l(g)\nu (dg)=$$
$$=\int_X \int_{G'}{\bar f}_s(g)\beta (h)f_l(g)\lambda (dh)\nu (dg)=
(f_s,a_jf_l).$$
Hence ${{\sf A}_G}"$ contains 
subalgebra of all operators of multiplication on functions from
$L^{\infty }(X,\nu ,{\bf C})$.
With $G'$ and a Banach algebra $\sf A$ 
the trivial Banach bundle ${\sf B}={\sf A}\times G'$ is associative, in 
particular let ${\sf A}=\bf C$ (see \S VIII.2.7 \cite{fell}).
\par The regular representation $T$ of $G'$ gives rise to a canonical regular
$H$-projection-valued measure $\bar P$:
$\bar P(W)f=Ch_Wf$, where $f\in H$, $W\in Bf(X)$, $Ch_W$ 
is the characteristic function of $W$. Therefore, $T_h\bar P(W)=\bar P
(h\circ W)T_h$ for each $h\in G'$ and $W\in Bf(X)$, since
$\rho _{\nu }(h,h^{-1}\circ g)\rho _{\nu }(h,g)=1=\rho _{\nu }(e,g)$ 
for each $(h,g)\in G'\times X$, 
$Ch_W(h^{-1}\circ g)=Ch_{h\circ W}(g)$ and $T_h(\bar P(W)f(g))
=\rho _{\nu }(
h,g)^{(1+bi)/2}\bar P(h\circ W)f(h^{-1}\circ g)$. Thus $<T,\bar P>$ is 
a system of imprimitivity for $G'$ over $X$, which is denoted 
${\sf T}^{\nu }$. This means that conditions
$SI(i-iii)$ are satisfied: 
\par $SI(i)$ $T$ is a unitary representation
of $G'$; 
\par $SI(ii)$ $\bar P$ is a regular 
$H$-projection-valued Borel measure on $X$ and 
\par $SI(iii)$ $T_h\bar P(W)=\bar P(h\circ W)T_h$ for all $h\in G'$ 
and $W\in Bf(X)$. 
\par For each $F\in L^{\infty }(X,\nu ,{\bf C})$ let $\bar \alpha _F$
be the operator in ${\sf L}(H)$ consisting
of multiplication by $F$: $\bar \alpha _F(f)=Ff$ for each $f\in H$. 
The map $F\to \bar \alpha _F$ is  an isometric $*$-isomorphism
of $L^{\infty }(X,\nu ,{\bf C})$ into ${\sf L}(H)$
(see \S VIII.19.2\cite{fell}). 
\par If $\bar p$ is a projection onto a closed ${\sf T}^{\nu }$-stable
subspace of $H$, then $\bar p$ commutes with all
$\bar P(W)$. Hence $\bar p$ commutes with multiplication by all
$F\in L^{\infty }(X,\nu ,{\bf C})$, so by \S VIII.19.2 \cite{fell}
$\bar p=\bar P(V)$, where $V\in Bf(X)$. Also $\bar p$ commutes with all
$T_h$, $h\in G'$, consequently, $(h\circ V)\setminus V$ and 
$(h^{-1}\circ V)\setminus V$ are $\nu $-null for each $h\in G'$, 
hence $\nu ((h\circ V)\bigtriangleup V)=0$ for all $h\in G'$. In view 
of ergodicity of $\nu $ and Proposition VIII.19.5 \cite{fell}
either $\nu (V)=0$ or $\nu (X\setminus V)=0$, hence
either $\bar p=0$ or $\bar p=I$, where $I$ is the unit operator.
Hence $T$ is the irreducible unitary representation.
\par {\bf 4.2. Theorem.} {\it On a loop group or a diffeomorphism group
or on their semidirect product $G$
there exists a stochastic process, which generates a quasi-invariant 
measure $\mu $ relative to a dense subgroup $G'$ such that an associated 
regular unitary representation $T^{\mu }: G' \to U(L^2(G,\mu ,{\bf C}))$ 
is irreducible. The family $\Psi $ of such pairwise nonequivalent
irreducible unitary representations has the cardinality
$card (\Psi ) = card ({\bf R})$.}
\par {\bf Proof.} From the construction of $G'$ and $\mu $ in
\S 3.2 and Theorem 3.3
it follows that, if a function $f\in L^1(G,\mu ,{\bf C})$ 
satisfies the following condition
$f^h(g)=f(g)$ $(mod $ $\mu )$ by $g\in G$ for each $h\in G'$, 
then $f(x)=const $
$( mod $ $\mu )$, where $f^h(g):=f(hg)$, $g\in G$.  
\par Let $f(g)=Ch_U(g)$ be
the characteristic function of a subset $U$, $U\subset G$, $U\in Af(G,\mu )$,
then $f(hg)=1 $ $\Leftrightarrow g\in h^{-1}U$.  If $f^h(g)=f(g)$ is
accomplished by
$g\in G$ $\mu $-almost everywhere, then $\mu (\{ g\in G:  f^h(g)\ne f(g) \}
)=0$, that is $\mu ( (h^{-1}U)\bigtriangleup U)=0$, consequently, 
the measure $\mu $ 
satisfies the condition $(P)$ from \S VIII.19.5 \cite{fell}, where
$A\bigtriangleup B:=(A\setminus B)\cup (B\setminus A)$ 
for each $A, B\subset G$.
For each subset $E\subset G$ the outer measure is bounded,
$\mu ^*(E)\le 1$, since $\mu (G)=1$ and $\mu $ is non-negative, 
consequently, there exists $F\in
Bf(G)$ such that $F\supset E$ and $\mu (F)=\mu ^*(E)$.  
This $F$ may be interpreted as the representative
of the least upper bound in $Bf(G)$ relative to
the latter equality.
In view of the Proposition VIII.19.5 \cite{fell} the 
measure $\mu $ is ergodic.
\par In view of Theorems 2.1.7, 2.8, 2.10 the Wiener process on
the Hilbert manifold $G$ induces the Wiener process on 
the Hilbert space $T_eG$ with the help of the manifold
exponential mapping. Then the left action $L_h$ of $G'$ on $G$
induces the local left action of $G'$ on a neighbourhood 
$V$ of $0$ in $T_eG$ with $\nu (V)>0$, where $\nu $ is induced by $\mu $.
A class of compact subsets approximates from below each measure $\mu ^f$, 
$\mu ^f(dg):=|f(g)|\mu (dg)$, where $f\in L^2(G,\mu ,{\bf C})=:H$.
Due to the Egorov Theorem II.1.11 \cite{fell} for each $\epsilon >0$
and for each sequence $f_n(g)$ converging to $f(g)$ for $\mu $-almost
every $g\in G$, when $n\to \infty $, there exists a compact subset $\sf K$
in $G$ such that $\mu (G\setminus {\sf K})<\epsilon $ and
$f_n(g)$ converges on $\sf K$ uniformly by $g\in \sf K$, when $n\to \infty $.
\par In view of Lemma IV.4.8 \cite{oksb}
the set of random variables $ \{ \phi (B_{t_1},...,B_{t_n}) :\quad
t_i \in [t_0,T], \phi \in C^{\infty }_0({\bf R^n}), n\in {\bf N} \} $
is dense in $L^2({\sf F}_T,\mu )$, where $T>t_0$. In accordance with
Lemma IV.4.9 \cite{oksb} the linear span of random variables
of the type $\{ \exp \{ \int_0^Th(t)dB_t(\omega )-\int_0^Th^2(t)dt/2 \} :\quad
h\in L^2[t_0,T]$ (deterministic) $ \} $ is dense in $L^2({\sf F}_T,\mu )$,
where $T>t_0$. 
Therefore, in view of Girsanov Theorem 2.1.1 and Theorem 5.4.2 \cite{ustzak}
the following space 
$sp_{\bf C}\{ \psi (g):=(\rho (h,g))^{(1+bi)/2}: h\in G' \}=:Q$ is dense in
$H$, since $\rho _{\mu }(e,g)=1$ for each $g\in G$
and $L_h: G\to G$ are diffeomorphisms of the manifold $G$, $L_h(g)=hg$,
where $b\in \bf R$ is fixed.
Finally we get from Theorem 3.3 above that there exists
$\mu $, which is ergodic and Conditions $(i,ii)$ of Theorem 4.1.2
are satisfied. Evidently $G'$ and $G$ are infinite and dense in themselves.
Hence from Theorem 4.1.2 the statement of this theorem, follows.
\par In view of Theorems 3.10 and 3.13 \cite{lupom} a pair of
representations $T^{\mu }$
and $T^{\nu }$ generated by quasi-invariant measures  $\mu $
and $\nu $ is equivalent if and only if measures $\mu $ and
$\nu $ are equivalent. Considering different Wiener processes
on $G$, their transition probabilities and using the Kakutani theorem
\cite{dalf} it is possible to construct a family $\Psi $ of pairwise
nonequivalent measures and representations such that
$card (\Psi ) = card ({\bf R})$.
\par {\bf 4.3. Note.} Analogously to \S 3.3 there can be constructed 
quasi-invariant and differentiable measures on the manifold
$M$ relative to the action of the diffeomorphism group 
$G_M$ such that $G'\subset G_M$. Then Poisson measures on 
configuration spaces associated with either $G$ or $M$ can be
constructed and producing new unitary representations 
including irreducible as in \cite{lupom}.
\par Having a restriction of a transition measure $\mu $ from \S 3.3
on a proper open neighbourhood of $e$ in $G$ it is possible
to construct a quasi-invariant $\sigma $-finite nonnegative
measure $m$ on $G$ such that $m(G)=\infty $ using left shifts
$L_h$ on the paracompact $G$.
Analogously such measure can be constructed on the
manifold $M$ in the case of the diffeomorphism group using Wiener 
processes on $M$. For definite $\mu $ in view of Theorems 2.9 \cite{lupom}
and 4.2 the corresponding Poisson measure $P_m$ is ergodic.
Therefore, Theorems 3.4, 3.6, 3.9, 3.10, 3.13 and 3.14 \cite{lupom}
also encompass the corresponding class of measures $m$ and $P_m$
arising from the constructed in \S 3.3 transition measures.
\par {\bf 4.4. Theorem.} {\it Let $G$ be an infinite dimensional
Lindel\"of $C^{\infty }$-Lie group
and $G'$ be its dense subgroup such that relative
to their own uniformities $G$ and $G'$ have structures of
Banach manifolds with the Hilbert-Schmidt operator of embedding
$A: T_eG' \hookrightarrow T_eG$. Suppose that $T: G'\to U(H)$
is a strongly continuous injective unitary representation of $G'$
such that $T(G')$ is a complete uniform subspace in $U(H)$
supplied with strong topology.
Then there exists on $G$ a quasi-invariant probability measure
$\mu $ relative to $G'$. If $T$ is topologically reducible,
then $\mu $ is isomorphic to product of measures $\mu _k$
on $G$ quasi-invariant relative to $G'$.}
\par {\bf Proof.} Since $G$ and $G'$ are infinite dimensional and
Lindel\"of, then $G$ and $G'$ have countable bases of neighbourhoods
of $e$ relative to their topologies $\tau $ and $\tau '$ respectively,
hence they are mertizable \cite{hew}. Since $G$ and $G'$ are Lindel\"of
and metrizable, then they are separable \cite{eng}.
Therefore, $T_eG$ and $T_eG'$ are Lindel\"of and separable.
\par In view of Proposition II.1 \cite{neeb} for the separable
Hilbert space $H$ the unitary group endowed with the strong 
operator topology $U(H)_s$ is the Polish group.
Let $U(H)_n$ be the unitary group with the metric
induced by the operator norm. In view of the Pickrell's theorem 
(see \S II.2 \cite{neeb}): if $\pi : U(H)_n\to U(V)_s$
is a continuous representation of $U(H)_n$ on  
the separable Hilbert space $V$, then $\pi $ is also
continuous as a homomorphism from $U(H)_s$ into
$U(V)_s$. Therefore, if $T: G'\to U(H)_s$ is
a continous representation, then there are new representations
$\pi \circ T: G'\to U(V)_s$. On the other hand, the 
unitary representation theory of $U(H)_n$ is the same as that of
$U_{\infty }(H) := U(H)\cap (1+L_0(H))$, since the group $U_{\infty }(H)$
is dense in $U(H)_s$, where $L_0(H)$ denotes the Banach space
of $\bf R$-linear compact operators from $H$ into $H$.
If $H_1$ is an invariant subspace of a representation $T$,
then $H_1$ is separable, since $G'$ is Lindel\"of and separable.
That is there exists a unitary operator $S\in U(H)$ such that
$ \{ ST_hS^*:$  $h\in G' \} $ leaves invariant subspaces $H_1$ and
$H\ominus H_1$, since each operator $ST_hS^*$ is unitary.
Therefore, there exists a representation
$\mbox{ }_1T: G'\to U(H_1)$, which is strongly continuous
and injective. For a Hilbert space $H_1$, the tangent space
$T_eU_{\infty }(H_1)$ can be supplied with the natural
Banach space structure.
\par The norm topology and strong operator topologies induce the same
algebra of Borel subsets of $U(H_1)$, since
relative to these topologies $U(H_1)$ is Lindel\"of.
Consider a rigged Hilbert space $X_+ \hookrightarrow
X_0\hookrightarrow X_-$ with a nondegenerate positive definite
nuclear operator of embedding $W: X_+ \hookrightarrow X_-$.
Choose $X_0\subset T_eU(H_1)$ and $W<A^2$. Then take on $X_-$
a Gaussian measure $\nu $ induced by a cylindrical Gaussian measure
$\lambda _I$ on $X_0$ with the unit correlation operator.
The unitary group $U_{\infty }(H_1)$ satisfies the Campbell-Hausdorff
formula, hence a measure $\nu $ induces a measure $\psi $
on $U_{\infty }(H_1)$. Consider a space $L_A(H_1)$ of all $\bf R$-linear
operators $K$ from $H_1$ into $H_1$ such that $KA$ is a bounded
operator and put $\| K \|_A := \| KA \| $, then $L_A(H_1)$
is the Banach space. Then the completion of $U(H_1)$ relative
to the uniformity induced by $L_A(H_1)$ gives the uniform space
denoted by ${\bar U}(H_1)$. Supply $L_A(H_1)$ by the strong topology
with a base of neighbourhoods of zero $W_{\epsilon }
(x_1,...,x_n):= \{ K\in L_A(H_1):$ $ \| KAx_j \|_{H_1}<\epsilon ,
j=1,...,n \} $, where $x_1,...,x_n\in H_1$, which generates a uniformity.
It induces the strong topology ${\bar s}$ in ${\bar U}(H_1)$.
Relative to such strong topology we denote it by ${\bar U}(H_1)_{\bar s}$.
\par Using cylindrical subsets generated
by projections on finite dimensional subspaces in $T_eU_{\infty }(H_1)$
with the help of finite dimensional subalgebras $T_eU(n)$ embedded
into $T_eU_{\infty }(H_1)$ induce a measure
$\psi $ from $U_{\infty }(H_1)$ on ${\bar U}(H_1)$, which is
the $C^{\infty }$-manifold. Consider a $C^{\infty }$-vector field $X$ in
${\bar U}(H_1)_{\bar s}$, then it induces a local one-parameter
group of diffeomorhisms
$g^t_X$ acting from the left on $U(H_1)$ and hence on ${\bar U}(H_1)$
such that $\partial g^t_X/\partial t|_{t=0}=X$, where $t\in (- \epsilon ,
\epsilon )\subset \bf R$, $\epsilon >0$. On the other hand,
$U_{\infty }(H_1)$ satisfies the Campbell-Hausdorff formula.
Thus $\psi $ can be chosen $\sigma $-additive and quasi-invariant
on the manifold ${\bar U}(H_1)$ relative to the left action
of $U(H_1)$ (see \cite{dalf}).
\par Since $T(G')$ is complete relative to the strong unifomity
inherited from $U(H)_s$, then $T(G')$ is closed in $U(H)_s$
(see about complete uniform spaces in \cite{eng}).
Thus $A: T_eG' \hookrightarrow T_eG$, $T: G'\to U(H)$ and
exponential mappings of $G'$ and $G$ as manifolds induce
an embedding of $G$ into ${\bar U}(H)_{\bar s}$.
Denote images of $G'$ and $G$ in $U(H)_s$ and in ${\bar U}(H)_{\bar s}$
under embeddings by the same letters $G'$ and $G$.
There exists a retraction $r: {\bar U}(H_1)_{\bar s}\to G$, since $G'$
is closed in $U(H_1)_s$ relative to the topology $s$ in $U(H_1)_s$
induced by the strong operator topology, where $r|_{G'}=id$,
$r(U(H_1))=G'$, $r({\bar U}(H_1))=G$,
$r: {\bar U}(H_1)_{\bar s} \to G$ and $r|_{U(H_1)}: U(H_1)_s
\to G'$ are continuous (see about retractions in \cite{isbell}).
Therefore, $\nu $ induces a Gaussian $\sigma $-additive measure $\zeta $
on $G$, where $\zeta (V):= \psi (r^{-1}(V))$ for each Borel subset $V$
in $G$. Since $G$ and $G'$ are $C^{\infty }$-Lie groups,
then they have $C^{\infty }$-manifold
structures such that $\exp : {\tilde T}G\to G$ and $\exp :
{\tilde T}G'\to G'$ are $C^{\infty }$-mappings, where
${\tilde T}G$ denotes a neighbourhood of $G$ in $TG$.
Thus $(\exp \circ A\circ \exp ^{-1})' -I$ is the Hilbert-Schmidt operator
and by Theorem II.4.4 \cite{dalf}
$\psi $ induces a $\sigma $-additive quasi-invariant measure
$\zeta $ on the Borel algebra of $G$ relative to the left action of $G'$,
since $T_{h_1}T_{h_2}=T_{h_1h_2}$ for each $h_1$ and $h_2$ in $G'$.
Instead of a concrete Gaussian measure it is possible to induce
a quasi-invariant measure on ${\bar U}(H_1)$ with the help of a
positive definite functional on $T_eU(H_1)$ satisfying the Sazonov theorem.
\par If $T$ is topologically reducible, then there exists at least
two invariant subspaces $H_1$ and $H_2$ in $H$ relative to
$ST_hS^*$ for each $h\in G'$, where $S$ is a fixed unitary operator.
Then a product measure on ${\bar U}(H_1)\times {\bar U}(H_2)$
induces a product measure on $G$. Since $H$ is separable, then
such product can contain only countable product of probability
measures. In view of the Kakutani theorem it can be chosen
quasi-invariant relative to the left action of $G'$.
\par {\bf 4.5. Proposition.} {\it Consider the semidirect product
$G := Diff^{\xi }_{y_0}(N) \otimes ^s (L^MN)_{\xi }$
of a group of diffeomorphisms and a group of loops, where
$\xi $ is such that $Y^{\xi }(M,N)\subset C^{\infty }(M,N)$. Then the
tangent space $T_eG =: \sf g$ can be supplied with the algebra structure.}
\par {\bf Proof.} Since $TC^{\infty }(M,N) = C^{\infty }(M,TN)$,
then $T_eDiff^{\xi }(N)$ is isomorphic to the algebra $\Sigma (N)$ of
$Y^{\xi }$-vector fields on $N$. For the proof consider foliated
structure of $M$ with foliated submanifolds $M_m$ in $M$
such that $\bigcup_m M_m$ is dense in $M$. Consider subgroups and
subspaces corresponding to restrictions on $M_m$ and then use completion
of the strict inductive limits of subgroups and subalgebars to get general
statement.
\par For each $f\in Y^{\xi }(M_m,s_{0,m};N,y_0)$
the Riemannian volume element $\nu _m$ in $M_m$, where
$dim_{\bf K} M_m =: m$, induces due to the Morse theorem natural
coordinates in $f(M_m)$ defined almost everywhere in $f(M_m)$ relative
to the measure $\mu _m$ on $f(M_m)$ such that $\mu _m(V):=
\nu _m(f^{-1}(V))$ for each Borel subset $V$ in $f(M_m)$.
Let $x_1,...,x_{km}$ be natural coordinates in $f(M_m)$,
where $k=dim_{\bf R}{\bf K}$, ${\bf K}=\bf R$ or $\bf C$ or $\bf H$.
Consider $g\in TY^{\xi }(M_m,s_{0,m};N,y_0)$, then $\lim_{x\to {\bar 1}}
pr_2(g(x))=:z_g\in T_{y_0}N$, where ${\bar 1}:=(1,...,1)\in \bf R^{km}$,
$x_l\in [0,1]$ for each $l=1,...,km$, $pr_2: TN\to Z$ is the natural
projection, where $\{ y \} \times Z=T_yN$ for each $y\in N$,
$Z$ is the vector space over $\bf K$. Therefore, $\lim_{x\to {\bar 1}}
pr_2(g(x)) - x_1...x_{km}z_g =0$, consequently,
$T_{w_0}Y^{\xi }(M_m,s_{0,m};N,y_0)$ is isomorphic with
$Y^{\xi }(M_m,s_{0,m};T_{y_0}N,y_0\times 0)\otimes Z$.
The latter space has $\bf K$-vector structure and the wedge
combination $g\vee f$ of mappings $g$ and $f$ and the equivalence
relation $R_{\xi }$ induce the monoid structure, hence $T_e(S^MN)_{\xi }$
is isomorphic with $(S^MZ)_{\xi }\otimes Z$ and inevitably
$T_e(L^MN)_{\xi }$ is isomorphic with $(L^MZ)_{\xi }\otimes Z$
which is the $\bf K$-vector space and $(f,v)\circ (g,w):=
(f\circ g,v+w)$ gives the algebra structure in $T_e(L^MN)_{\xi }$,
where $f, g\in (L^MZ)$, $v, w \in Z$. If $X, P\in \Sigma (N)$,
then there exists a Lie algebra structure $[X,P]$ in $\Sigma (N)$.
For $P\in Y^{\xi }(M_m,TN)$ there exists $\nabla _XP$.
Thus $T_e(Diff^{\xi }_{y_0}(N) \otimes ^s (L^MN)_{\xi })$
is isomorphic with the semidirect product of algebras
$\Sigma _0(N)\otimes ^s[(L^MZ)\otimes Z]$, where $\Sigma _0(N)$
is a subalgebra of all $X\in \Sigma (N)$ such that $\pi (X(y_0))=y_0$,
where $\pi : TN\to N$ is the natural projection.
\par {\bf 4.6. Definition.} We call the algebra $\sf g$ from
Proposition 4.5 by the vector field loop algebra.
The algebra $T_e Diff^{\xi }(N) =: {\sf g}(N)$ is called the algebra
of vector fields (in $N$). The algbera $T_e(L^MN)_{\xi }$ we call
the loop algebra.
\par For $N=S^4$ supplied with the quaternion manifold
structure (see \S 2.1.3.6) consider the semidirect product
of groups $Diff^{\sf H}(S^4)\otimes ^s\bf H$, where $\bf H$ is considered
as the additive group. Then we call $T_e (Diff^{\sf H}(S^4)
\otimes {\bf H})_{\bf O} =: ({\sf g}(S^4)\otimes {\sf h})_{\bf O}$
the quaternion Virasoro algebra, where ${\sf g}_{\bf O}$ denotes
the octonified algebra $\sf g$, $\bf O$ denotes the octonion division
algebra (over $\bf R$).
\par {\bf 4.7. Theorem.} {\it Let ${\sf g}'$ be a vector field
loop algebra or an algebra of vector fields or a loop algebra.
Then there exists a family $\Phi $ of the cardinality
$card (\Phi ) = card ({\bf R})$ of infinite
dimensional pairwise nonequivalent representations
${\sf t}: {\sf g}'\to {\sf gl}(H)$, where ${\sf gl}(H)$ denotes
the general linear algebra of the Hilbert space $H$ over $\bf R$.}
\par {\bf Proof.} Let $\cal L$ be the category of Lie groups
with differentiable morphisms, let also $\cal A$ be the category
of algebras over $\bf R$.
Use the tangent covariant functor $\cal T$ from the category
$\cal L$ into $\cal A$ such that for each object $L\in \cal L$
we have ${\cal T}(L)=A\in \cal A$ and ${\cal T}: Mor (L,L')\to Mor (A,A')$
such that ${\cal T}(\alpha \alpha ')={\cal T}(\alpha )
{\cal T}(\alpha ')$ for each $\alpha \in Mor (L,L')$
and $\alpha '\in Mor (L',L")$, ${\cal T}(1_A)=1_{{\cal T}(A)}$.
In particular, for a differentiable unitary representation
$\alpha \in Mor (G, U(H))$ we get ${\cal T}(\alpha )\in Mor
({\sf g}, {\sf u}(H))$.
\par Consider at first a unitary representation
of the group $G$ given by Theorem 4.2 such that it is induced
by the differentiable measure. Since $G$ has the $C^{\infty }$-manifold
structure, then the operator $L_h: G\to G$ such that $L_h(g):=hg$
for each $h, g\in G$ is strongly differentiable.
Relative to a strongly differentiable operator $S$
on a Banach space such that $S'-I$ is the Hilbert-Schmidt operator
the Gaussian transition measure transforms in accordance with
Theorem II.4.4 \cite{dalf}. In accordance with Theorem 3.4.3 \cite{beldal}
for each $k\in \bf N$ the solution $\xi $ of the stochastic equation
from \S 3.3 possesses $k$ bounded Fr\'echet derivatives relative to the
action of $L_h$ up to stochastic equivalence of the solution.
The exponential mapping $\exp $ of $G$ as the manifold is also of
class $C^{\infty }$. The measures considered in Theorem 3.3
are infinite differentiable relative to left shifts $L_h$ of a dense
subgroup $G'$ in $G$. Thus the quasi-invariance factor $\rho _{\mu }
(h,g)$ is also strongly differentiable by $h\in G'$. Therefore,
the irreducible unitary representation of $G'$ in $U(H)$
is strongly differentiable and induces the representation
${\sf t}: {\sf g}'\to {\sf u}(H)$, where ${\sf g}'=T_eG'$
and ${\sf u}(H)=T_eU(H)$. On the other hand, each unitary group
$U_{\infty }(H)$ of a complex separable Hilbert space $H$
is isomorphic with the general linear group $GL_{\infty }(H_{\bf R})$
of compact $\bf R$-linear operators from $H_{\bf R}$ into $H_{\bf R}$,
where $H_{\bf R}$ is the space $H$ considered over $\bf R$ which is
induced by the isomorphism of $\bf C$ as the $\bf R$-linear space with
$\bf R^2$. Thus ${\sf u}(H)$ is isomorphic with the general
$\bf R$-linear algebra ${\sf gl}_{\infty }(H_{\bf R})$ of all compact
$\bf R$-linear operators $w: H_{\bf R}\to H_{\bf R}$.
The group $U_{\infty }(H)$ is dense in $U(H)_s$,
consequently, strongly continuous irreducible representation
$T: G'\to U(H)$ induces the irreducible representation
${\sf t}: {\sf g}'\to {\sf gl}(H_{\bf R})$, where
${\sf gl}(H_{\bf R})$ denotes the algebra of all continuous
$\bf R$-linear operators from $H_{\bf R}$ into $H_{\bf R}$.
Since $card (\Psi _0)=card ({\bf R})$ in Theorem 4.2 for the
subfamily $\Psi _0$ of $\Psi $ of strongly differentiable unitary
representations, hence $card (\Phi ) = card ({\bf R})$.
\par {\bf 4.8. Theorem.} {\it The algebra $({\sf g}(S^4)\otimes ^s
{\sf h})_{\bf O}$ is the algebra over octonion division algebra
$\bf O$ such that there exists an embedding into it of the standard
Virasoro algebra of $S^1$. A set of generators of $({\sf g}(S^4)\otimes ^s
{\sf h})_{\bf O}$ is $ \{ \exp (lmz):$ $m\in {\bf Z} \}
\cup \{ y \} $, where ${\bf O}={\bf H}\oplus {\bf H}l$, $l$ is the
doubling generator of $\bf O$ over $\bf H$, $y, z \in \bf H$.}
\par {\bf Proof.} Mention that $\bf H$ is the Abelian Lie group
when $\bf H$ is considered as the additive group. Then to it there
corresponds the Lie algebra $\sf h$ over $\bf R$.
The unit sphere $S^4$ is homeomorphic with the one-point (Alexandroff)
compactification of $\bf H$. There exists the epimorphism
$\exp : {\cal I}_3 \to S({\bf O},0,1)$, where ${\cal I}_3:=
\{ y\in {\bf O}:$ $Re (y)=0 \} $, $S({\bf O},0,1) := \{ z\in {\bf O}:$
$|z|=1 \} $ (see Corollary 3.5 \cite{luoyst2}).
Therefore, $\exp (l{\bf H})$ is the four dimensional unit sphere
$S^4$ embedded into $S({\bf O},0,1)$.
In view of Corollary 3.4 \cite{luoyst2} $\exp (z(1 + 2 \pi k/ |z|))=
\exp (z)$ for each $0\ne z\in {\cal I}_3$, $\exp (0)=1$.
In the $\bf O$-vector space $C^0(B(l{\bf H},0,1),{\bf O})$ is dense
the subspace of polynomials $P_n(lz)=\sum_{|v|\le n}
\{ (a_v,(lz)^v) \} _{q(2|v|)}$, where $z\in \bf H$, $a_{v_k}\in \bf O$
for each $k$, $(a_v,x^v):=a_{v_1}x^{v_1}...a_{v_p}x^{v_p}$, $p\in \bf N$,
$v=(v_1,...,v_p)$, $0\le v_k\in \bf N$, $|v|:=v_1+...+v_p$,
$\{ b_1...b_p \} _{q(p)}$ denotes the product of octonions
$b_1,...,b_p$ in an association order prescribed by the vector
$q(p)$, $B(X,z,r) := \{ y\in X:$ $\rho (z,y)\le r \} $
denotes the ball in the metric space $(X,\rho )$, $r>0$,
$z\in X$ (see \S 2.1 \cite{luoyst2}).
From this it follows, that the $\bf O$-vector space $C^0(S^4,{\bf O})$ is
$\bf R$-linearly isomorphic with the $\bf O$-vector space of continuous
functions $f: {\bf H}\to \bf O$ periodic in the following manner:
$f(z+2\pi k z/|z|)=f(z)$ for each $0\ne z\in \bf H$ and each
$k\in \bf Z$. Each function $\exp (lmz)$ has the decomposition into
the series $\exp (lmz)=\sum_{p=0}^{\infty }(lmz)^p/p!$ converging
on $\bf H$, where $z\in \bf H$, $m\in \bf Z$. Therefore, the
system of equations
\par $(i)$ $(lz)^k=\sum_ma_{m,k}\exp (lmz)$ \\
is equivalent to $\sum_ma_{m,k}m^p/p!=\delta _{k,p}$
and the latter has the real solution
$a_{m,k}\in \bf R$ for each $m, k$, where $\delta _{m,k}=1$ for $m=k$
and $\delta _{m,k}=0$ for $m\ne k$ is the Kronerek delta.
These expansion coefficients $a_{m,k}$ can be expressed in the form:
\par $(ii)$ $a_{m,k} = \int_B (lz)^k\exp (lmz)\lambda (dz) / \lambda (B)$, \\
where $\lambda \ne 0$ denotes the measure on $\bf H$ induced by the
Lebesgue measure on $\bf R^4$, $B := B({\bf H},0,2\pi )$. 
\par Thus $K_m:=\exp (lmz) l^* \partial /\partial z$ is the basis of
generators of the algebra of vector fields on
$S^4 = \exp (l {\bf H})$, $z\in \bf H$,
where $\partial /\partial z$ denotes the superdifferentiation by
the quaternion variable $z$.
Consider the subgroup $Diff^{{\sf H},p}(B)$ of $Diff^{\sf H}(B)$
consisting of periodic diffeomorphisms $f(z(1+2\pi k/|z|))=f(z)$
for each $0\ne z\in \bf H$. Then put
\par $(iii)$  $c(f,g) := \int_B Ln (f'(g(z)).1) d Ln (g'(z).1)$ for each
$f, g\in Diff^{{\sf H},p}(B)$ \\
and define the semidirect product $Diff^{{\sf H},p}(B) \otimes ^s
{\bf H}$ such that
\par $(iv)$ $(g(z),y_1) (f(z),y_2) := (f(g(z)),y_1+y_2+c(f,g))$, \\
where $Ln $ denotes the logarithmic function for $\bf O$
(see about $Ln$ in \cite{luoyst2}).
This induces the semidirect product $Diff^{\sf H}(S^4)\otimes ^s{\bf H}$.
Therefore, $T_e(Diff^{\sf H}(S^4)\otimes ^s{\bf H})_{\bf O}$
is the algebra over $\bf O$ denoted by  $({\sf g}(S^4)\otimes ^s
{\sf h})_{\bf O}$ which is the octonification of
$T_e(Diff^{\sf H}(S^4)\otimes ^s{\bf H})$, that is obtained by extension
of scalars (expansion coefficients) from $\bf H$ to $\bf O$.
In view Formulas $(i-iv)$  $({\sf g}(S^4)\otimes ^s
{\sf h})_{\bf O}$ has the basis of generators $A_m = K_m + w_m y$, where
$w_m\in \bf R$ is a constant for each $m\in \bf Z$, $y\in \bf H$.
Thus $[A_m,y]=0$ for each $m\in \bf Z$ and $y$ belongs to the center
of this algebra, $y\in Z(({\sf g}(S^4)\otimes ^s{\sf h})_{\bf O})$.
Let $f(z)$ be a quaternion holomorphic function from $\bf H$ into $\bf O$.
Then $[K_n, K_m] f(z) = \{ (K_n\exp (lmz)) - (K_m\exp (lnz)) \} 
l^* \partial f(z)/ \partial z$. On the other hand,
$\partial f(z)/\partial z_p=(\partial f(z)/\partial z).i_p$, where
$z=\sum_{p=0}^3z_pi_p$, $z_p\in \bf R$, $i_p \in \{ 1, i, j, k \} $,
$\partial f(z)/\partial z$ is generally neither right nor left linear
operator in $\bf H$. It remains to calculate
$(K_n \exp (lmz)).i_p$ to verify, that $[A_n,A_m]\in
({\sf g}(S^4)\otimes ^s{\sf h})_{\bf O}$ and embed into it the Virasoro
algebra of $S^1$, where $i_p\in \{ 1, i, j, k \} $,
$p=0,1,2,3$. Evidently, $[K_n,K_m].1 f(z) := (m-n)K_{m+n}.1 f(z)$,
where $(\partial /\partial z).1 =: \partial /\partial z_1$.
Then $(\partial \exp (lmz)/\partial z).i_p
=m(\sum_{n=1}^{\infty }\sum_{k=0}^{n-1}((lmz)^k(li_p)) (lmz)^{n-k-1}/n!)$,
since $\bf O$ is power-associative.
There are identites: $(lz)(lw)=-w{\tilde z}$, $((lz)(li_p))(lz)=
z(i_p{\tilde z})$, $(lz)(li_p)=-i_p{\tilde z}$, $(li_p)(lz)=zi_p$
for each $z, w\in \bf H$, $p=0,...,3$, where ${\tilde z}=z^*$
is the conjugated quaternion $z$. Therefore,
\par $(\partial \exp (lmz)/ \partial z).i_p= $ \\
$m(\sum_{n=1}^{\infty }\sum_{k=0}^{n-1}(-|mz|)^{[k/2]+[(n-k-1)/2]}
((lmz)^{k-2[k/2]}(li_p)) (lmz)^{n-k-1-2[(n-k-1)/2]}/n!)$, \\
since $\bf R$ is the centre of $\bf O$, where $[a]$ denotes the integer
part of $a\in \bf R$, $[a]\le a$, $p=1,2,3$. Hence
\par $(v)$ $(\partial \exp (lmz)/ \partial z).i_p=m (e^{lmz} + e^{-lmz})
((lz)^{-1}(i_pz-((lz)l^*)i_p)) $ \\
$+ m (e^{lmz}-e^{-lmz}) ((lz)^{-1}(zi_p - i_p((lz)l^*)) -
((e^{lmz}-e^{-lmz})(lz)^{-1}) ((lz)^{-1}
(i_pz+((lz)l^*)i_p))/4$, \\
where $(lz)^{-1}(i_pz-((lz)l^*)i_p)$ and $((lz)^{-1}(zi_p - i_p((lz)l^*))$ 
and $(lz)^{-1}(i_pz+((lz)l^*)i_p)$  do not depend on $|z|$, $p=1,2,3$.
The equations $(lz)^{-1}i_pz=\sum_{m=0}^{\infty }e^{lmz}v_m$
and $(lz)^{-1}((lz)l^*)i_p=\sum_{m=0}^{\infty }e^{lmz}w_m$
are equivalent to $i_pz=(lz)\sum_{m=0}^{\infty }\sum_{n=0}^{\infty }
(lmz)^nv_m/n!$ and
$((lz)l^*)i_p=(lz)\sum_{m=0}^{\infty }\sum_{n=0}^{\infty }
(lmz)^nw_m/n!$ respectively, which evidently has the solutions
with expansion coefficients $v_m, w_m \in {\bf H}\cup l{\bf H}$,
since $(lz)^n=(-|lz|^2)^{[n/2]}
(lz)^{n-2[n/2]}$ for each $n\in \bf N$, where $z\in \bf H$,
$p=1,2,3$. Using decomposition of
$(e^{lmz}-e^{-lmz}) (lz)^{-1}$ into the power series by $(lz)^k$
and Equation $(i)$ we get the expansion 
$(e^{lmz}-e^{-lmz}) (lz)^{-1}=\sum_{k=0}^{\infty }b_{m,k}e^{kmz}$, where
$b_{m,k}$ are real coefficients. Using $f(z)=\sum_{p=0}^7f_p(z)i_p$
with $f_p(z)\in \bf R$ for each $p$ we get, that
\par $(vi)$ $[K_n, K_m] f(z) = \sum_{s\in \bf Z}\sum_{p=0}^7
(K_s a_{s,p})(f(z)b_{s,p})$, \\
where $i_p\in \{ 1, i, j, k, l, li, lj, lk \} $, $p=0,...,7$,
$a_{s,p}$ and $b_{s,p}\in {\bf O}$ are octonion constants
for each $s, p$.
\par The usual Virasoro algebra $Vir$ of $S^1$ is the complexification
of the algebra of vector fields on $S^1$ centrally extended and it has
generators $L_n:= Y_n + v_nx$, where $v_n\in \bf R$ for each
$n\in \bf Z$, $x\in \bf R$, $Y_n:=\exp (in\phi )i^*\partial /\partial \phi $,
where $\phi $ is the polar angle parameter on $S^1$, $\phi \in [0,2\pi ]$,
$i=(-1)^{1/2}$. The commutation relations are: $[L_n,L_m]:=(m-n)L_{m+n}+
s (n^3-n)\delta_{n,-m}x/12$, $[L_n,x]=0$ for each $n, m\in \bf Z$,
where $s$ is a real constant.
Thus $Vir$ has the embedding into $({\sf g}(S^4)\otimes ^s{\sf h})_{\bf O}$.
\par {\bf 4.9. Corollary.} {\it The algebra
${\sf g}_{\bf O}:=({\sf g}(S^4)\otimes ^s  {\sf h})_{\bf O}$
has a family $\Phi $ of the cardinality $card (\Phi ) = card ({\bf R})$
of infinite dimensional pairwise nonequivalent representations
${\sf t}: {\sf g}_{\bf O}\to {\sf gl}(H)_{\bf O}$.}
\par {\bf Proof.} There exists a group $G=Diff^{\xi }
(S^4)\otimes ^s\bf H$ such that $T_eG'$ has an embedding into $T_eG$
of the Hilbert-Schmidt class, where $G'=Diff^{\sf H}(S^4)\otimes ^s\bf H$
(see \S 2.1.5).
In view of Theorems 4.2 and 4.7 ${\sf g}:=T_eG'$ has a family $\Phi $ of
representations into $T_e{\sf u}(H_{\bf R})$. Each ${\sf t}\in \Phi$
has the natural extension on the octonification: ${\sf t}: {\sf g}_{\bf O}
\to {\sf u}(H_{\bf R})_{\bf O}$.

\par Address: Mathematical Department, Brussels University, V.U.B.,
\par Pleinlaan 2, Brussels 1050, Belgium
\end{document}